\documentclass[british]{scrartcl}
\usepackage[T1]{fontenc}
\usepackage[latin9]{inputenc}
\usepackage{float}
\usepackage{bm}
\usepackage{amsmath}
\usepackage{accents}
\usepackage{amssymb}
\usepackage{graphicx}
\usepackage{pdfpages}
\usepackage{esint}
\usepackage[authoryear]{natbib}
\usepackage{amsthm}
\usepackage{multirow}
\usepackage{url}
\usepackage{color}
\usepackage{enumitem}
\usepackage{hyperref}
\usepackage{caption}
\usepackage{tikz}
\usepackage{pdfsync}
\usepackage{subcaption}
\usepackage{todonotes}
\usetikzlibrary{positioning}
\usetikzlibrary{patterns}
\usetikzlibrary{shapes}
\usetikzlibrary{plotmarks}

\makeatletter

\floatstyle{ruled}
\newfloat{algorithm}{tbp}{loa}
\providecommand{\algorithmname}{Algorithm}
\floatname{algorithm}{\protect\algorithmname}

\usepackage{dsfont}\usepackage{algpseudocode}
\usepackage{mathtools}

\newcounter{rmq}[section]
\newcommand{\R}{\mathbb{R}}

\newcommand{\setX}{\mathcal{X}}
\newcommand{\setY}{\mathcal{Y}}

\renewcommand{\P}{\mathbb{P}}
\newcommand{\E}{\mathbb{E}}

\newcommand{\ind}{\mathds{1}}
\newcommand{\dd}{\mathrm{d}}

\newcommand{\bigO}{\mathcal{O}} 

\newcommand{\comment}[1]{ \ifthenelse{ \equal{\showcomment}{true} }{ {\bf #1} }{} }
\newcommand{\showcomment}{true}

\newtheorem{theorem}{Theorem}
\newtheorem{proposition}{Proposition}
\newtheorem{corollary}{Corollary}
\newtheorem{lemma}{Lemma}

\newtheorem{definition}{Definition}
\newtheorem{remark}{Remark}

\theoremstyle{plain}
\newtheorem{assumption}{Assumption}

\theoremstyle{plain}

\theoremstyle{plain}

\let\oldReturn\Return
\renewcommand{\Return}{\State\oldReturn}

\newcommand{\reqstart}{
    \begin{list}{\thereqcount}{\usecounter{reqcount}}
    \setcounter{reqcount}{\value{reqcountbackup}}
}

\newcommand{\reqend}{
    \setcounter{reqcountbackup}{\value{reqcount}}
    \end{list}
}

  \newcommand\iid{\stackrel{\mathclap{\normalfont\mbox{\tiny{iid}}}}{\sim}}

\DeclareMathOperator*{\argmin}{argmin}

\makeatletter
\newcommand*{\vneq}{%
  \mathrel{%
    \mathpalette\@vneq{=}%
  }%
}
\newcommand*{\@vneq}[2]{%
  \sbox0{\raisebox{\depth}{$#1\neq$}}%
  \sbox2{\raisebox{\depth}{$#1|\m@th$}}%
  \ifdim\ht2>\ht0 %
    \sbox2{\resizebox{\vneqxscale\width}{\vneqyscale\ht0}{\unhbox2}}%
  \fi
  \sbox2{$\m@th#1\vcenter{\copy2}$}%
  \ooalign{%
    \hfil\phantom{\copy2}\hfil\cr
    \hfil$#1#2\m@th$\hfil\cr
    \hfil\copy2\hfil\cr
  }%
}
\newcommand*{\vneqxscale}{1}
\newcommand*{\vneqyscale}{1}

\makeatother

\usepackage{babel}

\usepackage{babel}
\title{Universal Robust Regression via Maximum Mean Discrepancy}

\author{Pierre Alquier$^{(1)}$, Mathieu Gerber$^{(2)}$
\\
\small{(1) ESSEC Business School, Asia-Pacific campus, Singapore}
\\
\small{(2) School of Mathematics, University of Bristol, UK }
}
\date{}

\begin{document}
\maketitle
 
\begin{abstract}
Many modern datasets are collected automatically and are thus easily contaminated by outliers. This led to a regain of interest in robust estimation, including new notions of robustness such as robustness to adversarial contamination of the data. However, most robust estimation methods are designed for a specific model. Notably, many methods were proposed recently to obtain robust estimators in linear models (or generalized linear models), and a few were developed for very specific settings, for example beta regression or sample selection models.
In this paper we develop a new approach for robust estimation in arbitrary regression models, based on Maximum Mean Discrepancy minimization. We build two estimators which are both proven to be robust to Huber-type contamination. We obtain a  non-asymptotic error bound for one them and show that it is also robust to adversarial contamination,  but this estimator is computationally more expensive to  use in practice than the other one. As a by-product of our theoretical analysis of the proposed estimators we derive new results on kernel conditional mean embedding of distributions which are of independent interest.
\end{abstract}

\section{Introduction}
\label{section:intro}

Robustness  is a fundamental problem in statistics, which aims at using statistical procedures that remain stable in presence of outliers. Historically, outliers were mistakes in data collection, or observations of individuals belonging to a different population than the population of interest. Robustness became even more important in the modern context where automatically collected datasets are often heterogeneous. Moreover, some strategic datasets are susceptible of malevolent manipulations.

In the statistical literature, the development of robust estimation methods for regression models generally focusses on the construction of $Z$-estimators   \cite[Chapter 5]{Vaart2000} for which each individual  observation can only have a bounded impact on the  estimating equations, and which therefore have a bounded influence function \citep{hampelrobust,hampelrobust2}. This   strategy has been successfully applied for the  robust estimation of the regression coefficients in generalized linear models    \citep{Kunsch1989,cantoni2001,cantoni2006}  as well as e.g.\ robust inference in the negative binomial regression model with unknown overdispersion parameter \citep{Aeberhard2014} and in the Heckman sample selection model \citep{Zhelonkin2016}. Based on other approaches, robust estimators for    mixtures of linear regression models  \citep[see e.g.][]{bai2012}, for the Beta regression model with unknown precision parameter  \citep{Ghosh2019} and for  robust linear least squares regression \citep{audibert2011robust} have been developed.

In the past ten years there was a renewed interest for robust methods in the machine learning community. \citet{catoni2012challenging} developed a loss function whose minimization leads to robust estimators of the expectation of a random variable, a technique which was then adapted to many situations including linear regression~\citep{catoni2017dimension}. More generally, Lipschitz loss functions, such as the absolute loss  or Huber's loss~\citep{huber1964robust}, lead to robustness of the empirical risk minimization procedure, a fact that was used in~\citet{chinot2018statistical,alquier2019estimation,chinot2019robust,holland2019pac} to study robust procedures of classification and regression. The Median-of-Means  approach of~\citet{nemi1983problem,devroye2016sub} was also adapted to various machine learning problems~\citep{lugosi2016risk,lecueSPA,lugosi2019regularization,depersin,lecueML,lecue2020robusta}, including least squares regression, logistic regression, quantile regression and classification under various losses. 

In the discussion by Sture Holm in \cite{bickel1976another}, as well as in \cite{parr1980minimum}, minimum distance estimation is identified as a way to obtain robust estimators.  Building on this idea, \cite{Basu1998} introduced a density power divergence minimization  approach for robust inference in parametric models for i.i.d.\ observations. This procedure  is extended to  regression models in \cite{Ghosh2013} but suffers from two limitations. Firstly, the optimization of the objective function is, in general, a computationally challenging problem. Secondly, there is no general result which guarantees that the resulting $M$-estimator is robust. Its influence function is however known to be bounded   for the Gaussian linear regression model  \citep{Ghosh2013},  for the Poisson  and logistic regression models \citep{Ghosh2016} and  for the Beta regression model with unknown precision parameter \citep{Ghosh2013}.

In this paper we introduce a new minimum distance estimation strategy for parameter inference in regression models which (a) is proven to be robust to outliers, both in the fixed and random design setting, under general conditions on the statistical model and (b) only requires to be able to sample from the model and to compute the gradient of its log-likelihood function to be  applicable. In this sense, the approach proposed in this work defines a  universal robust regression method. More specifically, we present in this paper a minimum distance estimation procedure for regression models based on the Maximum Mean Discrepancy (MMD) distance.  

The use of  the MMD distance based on bounded kernels for robust minimum distance estimation   was proposed in \citet{barp2019minimum} and in \citet{cherief2019finite} \citep[see also][for a Bayesian type estimator]{cherief2019mmd}. When unbounded kernels are used  the ``automatic'' robustness induced by the MMD metric is  lost, in which case    \citet{lerasle2019monk}  propose  a Median-of-Mean procedure to robustify the MMD based estimator. However,  all these references focus on the simple case where we have a fully parametric model for the distribution of the data.

In this paper, we first extend the MMD based minimum distance approach of \citet{barp2019minimum,cherief2019finite} to the regression setting. This task is non-trivial, especially  in the random design scenario where we have a statistical model only for the distribution of $Y|X$,  and not  for the distribution of $X$.  If the distribution $P^0_X$ of $X$ is known then the method in the latter two references can be used for robust inference since, in this case, if we have a model  $\{P^\theta_{Y|X},\,\theta\in\Theta\}$  for the distribution of $Y|X$ then we have the model  $\{ P^0_X P^{\theta}_{Y|X},\,\theta\in\Theta\}$   for the distribution of  the pair $(X,Y)$. In practice, $P^0_X$  is generally unknown but we can use the observations $\{X_i\}_{i=1}^n$ to compute an empirical estimate $\hat{P}^{n}_X$ of this distribution. In this work, we consider the natural ideal of using $\{ \hat{P}^{n}_X P^{\theta}_{Y|X},\,\theta\in\Theta\}$ as the model for the distribution of  $(X,Y)$  and then to  estimate  $\theta$ using the  approach introduced   in \citet{barp2019minimum,cherief2019finite}. As   shown below, it turns out  that replacing  $P^0_X$ by the non-parametric  estimator $ \hat{P}^{n}_X$ preserves the convergence and robustness properties of the estimator. 

 More precisely, we prove that the resulting estimator $\hat{\theta}_n$ of the model parameter is,   in the fixed and random design setting, (i) universally consistent, in the sense that it will always converge to the best approximation,  in the sense of the MMD distance,  of the truth in the model  without  any  assumption on the distribution generating the observations, and (ii) robust to adversarial contaminations. As a by-product of our theoretical analysis of $\hat{\theta}_n$ in the random design setting  we derive results on kernel conditional  mean embedding of distributions that are of independent interest. 
 
Computing $\hat{\theta}_n$ requires to optimize a function involving a sum over $n^2$ terms. We introduce a stochastic gradient algorithm which allows to efficiently compute this estimator from several thousands data points, but the use of $\hat{\theta}_n$ remains  computationally expensive for large datasets. For this reason, we introduce an alternative estimator $\tilde{\theta}_n$  which, as argued below, is expected to have a similar behaviour to  $\hat{\theta}_n$ in practice    while having the advantage  to be defined by an objective function involving only a sum over $n$ terms. We also establish   that this estimator  is itself   robust to outliers,  in the sense that its  influence function is bounded.

Throughout this work   $\mathcal{X}$ and  $\mathcal{Y}$ are two topological   spaces, equipped respectively with the Borel $\sigma$-algebra $\mathfrak{S}_\mathcal{X}$ and $\mathfrak{S}_\mathcal{Y}$, and, letting  $\mathcal{Z}= \mathcal{X}\times \mathcal{Y}$ and $\mathfrak{S}_\mathcal{Z}=\mathfrak{S}_\mathcal{X}\otimes \mathfrak{S}_\mathcal{Y}$, we denote by $\mathcal{P}(\mathcal{Z})$   the set of  probability distributions on $(\mathcal{Z},\mathfrak{S}_\mathcal{Z})$.  We assume below that any distribution $P\in\mathcal{P}(\mathcal{Z})$  for $(X,Y)$ admits a regular conditional probability for the distribution of $Y$ given $X$\footnote{This is for instance the case if $\setX$ and $\setY$ are two Polish spaces.} and that all the random variables are defined on the same  probability space $(\Omega,\mathcal{F},\mathbb{P})$.

\section{Kernel mean embedding of distributions: Background and new results\label{sec:kernel}}

\subsection{Notation and convention\label{sub:notation}}

We let $k:\mathcal{Z}^2 \rightarrow \mathbb{R}$ be a kernel on $\mathcal{Z}$, i.e.~$k$ is
symmetric and positive definite, and  denote by $(\mathcal{H}, <\cdot,\cdot>_{\mathcal{H}})$   the reproducing kernel Hilbert space (RKHS) over $\mathcal{Z}$ having $k$ as reproducing kernel  \citep[see][for a comprehensive introduction to RKHSs]{muandet2017kernel}. In addition, we let $k_{\mathcal{X}}$ be a kernel on $\mathcal{X}$, $k_{\mathcal{Y}}$ be a kernel  on $\setY$  and   denote by $(\mathcal{H}_{\mathcal{X}}, <\cdot,\cdot>_{\mathcal{H}_{\mathcal{X}}})$ and by $(\mathcal{H}_{\mathcal{Y}}, <\cdot,\cdot>_{\mathcal{H}_{\mathcal{Y}}})$ the   RKHS induced by $k_{\mathcal{X}}$ and by $k_{\mathcal{Y}}$, respectively.   When there is no ambiguity,   with a slight abuse of language we will refer to $\mathcal{H}$    as the RKHS on $\mathcal{Z}$  induced by $k$, although the full characterization of an RKHS requires to specify both a function space and an inner product. The same abuse of language will be used for the RKHSs induced by $k_{\mathcal{X}}$ and by $k_{\mathcal{Y}}$. In the following we denote by $k_{\mathcal{X}}\otimes k_{\mathcal{Y}}$  the product kernel on $\mathcal{Z}$ such that $k_{\mathcal{X}}\otimes k_{\mathcal{Y}}\{(x,y),(x',y')\} = k_{\mathcal{X}}(x,x') k_{\mathcal{Y}} (y,y')$ for all $(x,y),(x',y')\in \mathcal{Z}$, and  by $K_{\alpha,\gamma}$ the Mat\'ern kernel on $\R^d$ with bandwidth parameter $\gamma>0$ and smoothness parameter $\alpha\in (0,\infty]$.  

We refer  to  Example 2.2 in \cite{kanagawa2018gaussian} for the general definition $K_{\alpha,\gamma}$ but we mention here a few useful properties of this  kernel. Firstly,  $K_{\alpha,\gamma}$ reduces to the exponential kernel when $\alpha=1/2$, i.e.~$K_{1/2,\gamma}(x,x')=\exp(-\|x-x'\|/\gamma)$, and to the Gaussian kernel when $\alpha=\infty$, i.e.~$K_{\infty,\gamma}(x,x')=\exp(-\|x-x'\|^2/\gamma)$. Secondly,  for all $\alpha\in(0,\infty]$ there exists a function $\tilde{K}_\alpha:[0,\infty)\rightarrow\R$ such that $K_{\alpha,\gamma}(x,x')=\tilde{K}_\alpha(\|x-x'\|/\gamma)$ for all $(x,x')\in\mathcal{X}^2$ and all $\gamma>0$, implying that the Mat\'ern kernel is a translation invariant kernel. Lastly, for all $\alpha\in(0,\infty]$ and $\gamma>0$ the Mat\'ern kernel is bounded and continuous.

\subsection{Kernel mean embeddings and the maximum mean discrepancy metric}

Assume that the following condition on $k$ holds:

\begin{assumption}\label{asm:k}
 The kernel $k$ is $\mathfrak{S}_\mathcal{Z}$-measurable and such that $|k|\leq 1$.
\end{assumption}

Then, for any probability distribution $P\in\mathcal{P}(\mathcal{Z})$, the quantity  $\mu(P) := \mathbb{E}_{Z\sim P}[k(Z,\cdot)] $, called the  mean embedding of $P$, is well defined in $\mathcal{H}$. If in addition $k$ is such that the mapping $P\mapsto \mu(P)$ is one-to-one, in which case    $k$ is said to be  characteristic,  the mapping  $\mathbb{D}_k:\mathcal{P}(\mathcal{Z})^2\rightarrow [0,2]$  defined by
$$ 
\mathbb{D}_k(P,Q) = \left\| \mu(P) - \mu(Q) \right\|_{\mathcal{H}}, \quad P,Q\in \mathcal{P}(\mathcal{Z})^2,
$$
is a metric on $\mathcal{P}(\mathcal{Z})$, known  as the maximum mean discrepancy metric. We stress that none of the results presented in this work   requires $k$ to be characteristic but it is only under this assumption on $k$ that they provide useful and interpretable convergence guarantees for the proposed estimators.

To better understand the properties of the MMD distance it is useful to express it as an integral probability metric as follows  \citep[see, e.g.][]{muandet2017kernel}
$$ 
\mathbb{D}_k(P,Q) =\sup_{f\in\mathcal{H}:\,\|f\|_{\mathcal{H}}\leq 1}\big|\E_{Z\sim P}\big\{f(Z)\big\}-\E_{Z\sim Q}\big\{f(Z)\big\}\big|, \quad P,Q\in \mathcal{P}(\mathcal{Z})^2.
$$
Two distributions $P$ and $Q$ are therefore close to each other in the sense of the MMD distance if  $\E_{Z\sim P}\big\{f(Z)\big\}\approx\E_{Z\sim Q}\big\{f(Z)\big\}$ for all functions  $f\in\mathcal{H}$ with norm at most 1. Both the norm and the set $\mathcal{H}$ depends on $k$ but, under Assumption \ref{asm:k} and if $k$ is characteristic,   $\mathcal{H}$ is dense in $\mathcal{C}_b(\mathcal{Z})$. In the special case where $\mathcal{Z}\subseteq\R^{d_z}$ for some $d_z\in\mathbb{N}$ and $k$ is  characteristic and translation invariant,   $\mathbb{D}_k(P,Q)$ is the distance in $L_2(\R^d, \eta_k)$ between the characteristic functions of $P$ and $Q$, for some probability measure $\eta_k\in\mathcal{P}(\R^{d_z})$. For instance, $\eta_k$ is a Gaussian distribution if $k$ is a Gaussian kernel and is a  Cauchy distribution if $k$ is an exponential kernel. We refer to  Table 2.1 of \citet{muandet2017kernel} for the expression of $\eta_k$ for various popular translation invariant kernels.

The rest of this section is relevant only for Section \ref{subsec:random_hat}, and can be skipped in a first reading.

\subsection{Kernel conditional mean embeddings}

In the rest of this section we let $(P_{Y|x})_{x\in\mathcal{X}}$ be a collection of distributions on $(\mathcal{X},\mathfrak{S}_\mathcal{X})$, $P_X\in\mathcal{P}(\mathcal{X})$ and $P=P_X P_{Y|\cdot}\in\mathcal{P}(\mathcal{Z})$. In addition, we  assume that the following assumption on $k$    holds:

\begin{assumption}
 \label{asm:kxy}
$k=k_{\mathcal{X}}\otimes k_{\mathcal{Y}}$ where $k_{\mathcal{X}}$ and $k_{\mathcal{Y}}$ are   continuous   on $\mathcal{X}^2$ and  on $\mathcal{Y}^2$, respectively, and such that $|k_{\mathcal{X}}|\leq 1$ and  $|k_{\mathcal{Y}}|\leq 1$. 
\end{assumption}

This  assumption imposes that $\mathcal{H}=\mathcal{H}_{\mathcal{X}}\otimes \mathcal{H}_{\mathcal{Y}}$  and we note   that the kernel $k$ is characteristic if  $k_{\mathcal{X}}$ and $k_{\mathcal{Y}}$  have the additional properties to be translation invariant and characteristic \citep{szabo2018characteristic}.   Under Assumptions \ref{asm:k}-\ref{asm:kxy} we can define a conditional mean embedding operator of   $P$. Such an operator is a mapping $\mathcal{C}_{Y|X}:\mathcal{H}_{\mathcal{X}}\rightarrow \mathcal{H}_{\mathcal{Y}}$ such that
\begin{align}\label{eq:CMB2}
\mathcal{C}_{Y|X} k_{\mathcal{X}}(x,\cdot) =\mu(P_{Y|x}),\quad\forall x\in\mathcal{X}.
\end{align}
Since \eqref{eq:CMB2}  depends on $P$ only through $(P_{Y|x})_{x\in\mathcal{X}}$   in what follows we will often refer to $\mathcal{C}_{Y|X}$ as the conditional mean embedding operator of $(P_{Y|x})_{x\in\mathcal{X}}$.

A valid definition of a conditional mean embedding operator was only recently proposed by \citet{klebanov2020rigorous}, see also \citet{mollenhauer2020nonparametric, li2022optimal},  and is  expressed in term of the   uncentred covariance operator 
 $ \mathcal{C}_P : \mathcal{H}_{\mathcal{X}} \rightarrow \mathcal{H}_{\mathcal{Y}} $ and in term of the  uncentred cross-covariance operator  $ \mathcal{C}_{P_X} : \mathcal{H}_{\mathcal{X}} \rightarrow \mathcal{H}_{\mathcal{X}} $, which are such that, for all $f_1,f_2\in \mathcal{H}_{\mathcal{X}}$ and $g \in \mathcal{H}_{\mathcal{Y}}$,
\begin{equation}\label{eq:cov_op}
 \left< g, \mathcal{C}_P f_1 \right>_{\mathcal{H}_{ \mathcal{Y}}}   =  \E_{(X,Y)\sim P} \left\{ g(Y) f_1(X) \right\},\quad
\left< f_1, \mathcal{C}_{P_X} f_2 \right>_{\mathcal{H}_{\mathcal{X}}}   =  \E_{X\sim P_X} \left\{ f_1(X)f_2(X) \right\}.
\end{equation}
The boundedness of $k_{\mathcal{X}}$ and $k_{\mathcal{Y}}$, imposed by Assumption~\ref{asm:kxy}, implies that $\mathcal{C}_{P}$ and $\mathcal{C}_{P_X}$ exist, are unique, and that they are bounded  linear operators \citep[see][Section 3]{fukumizu2004}.

 With this notation in place, we obtain the following slight variation of \citet[][Theorem 5.3]{klebanov2020rigorous}, whose proof is given in the Appendix \ref{p-lemma:CBM_New} for sake of completeness.
\begin{lemma}\label{lemma:CBM_New}
Assume that Assumptions \ref{asm:k}-\ref{asm:kxy} hold, that the function $f\in\mathcal{H}_{\mathcal{X}}$ is such that $\E_{X\sim P_X}\{f^2(X)\}=0$ if and only $f=0$, and that
\begin{align}\label{eq:cond_embed}
\E_{Y\sim P_{Y|\cdot}}\big\{g(Y)\}\in\mathcal{H}_{\mathcal{X}},\quad\forall g\in\mathcal{H}_{\mathcal{Y}}.
\end{align}
Then,  \eqref{eq:CMB2} holds for the bounded linear operator $\mathcal{C}_{Y|X}=\big(\mathcal{C}_{P_X}^{\dagger}\mathcal{C}_P^*\big)^*$, where $\mathcal{C}_{P_X}^{\dagger}$ denotes the Moore-Penrose pseudo-inverse of  $\mathcal{C}_{P_X}$. 
\end{lemma}

In this work, we use the following property of conditional mean embedding operators:

\begin{lemma}
 \label{lemma:constant:kernel}
For  $P_X',P_X''\in\mathcal{P}(\mathcal{X})$, let  $P' = P_X' P_{Y|\cdot}$ and $P''= P_X'' P_{Y|\cdot}$, and assume that   $(P_{Y|x})_{x\in\mathcal{X}}$ admits a   bounded linear conditional mean embedding operator $\mathcal{C}_{Y|X}$. Then,
\begin{align}\label{eq:bound_d}
\mathbb{D}_k(P',P'')
  \leq
 \|\mathcal{C}_{Y|X}\|_{\mathrm{o}}\,   \mathbb{D}_{k^2_{\mathcal{X}}}(P_X',P_X'').
 \end{align}

\end{lemma}

In words,    conditional mean embedding operators allow to quantify  the MMD distance between two joint distributions $P'$ and $P''$ on $(\mathcal{Z},\mathfrak{S}_\mathcal{Z})$, having the same  regular   conditional distribution $(P_{Y|x})_{x\in\mathcal{X}}$, in term of the MMD distance, induced by the kernel $k_{\mathcal{X}}^2$, between  their marginals $P'_X$ and $P''_X$.

\subsection{Some clarifications on kernel conditional mean embeddings \label{sub:cond_op}}

The main difficulty of the theory of conditional mean embedding operators is that condition \eqref{eq:cond_embed}, introduced by \citet{song2009hilbert}, is hard to interpret. This condition is slightly weakened in    \citet{klebanov2020rigorous} but the  alternative assumptions   proposed in this reference  remain hard to interpret since, as in \eqref{eq:cond_embed}, they require that all functions in $\mathcal{H}_{\mathcal{Y}}$ satisfy a given property.  In fact, to the best of our knowledge,  there is no explicit examples of   distributions $(P_{Y|x})_{x\in\mathcal{X}}$ for which a conditional mean embedding operator exists, beyond the trivial case where $P_{Y|x}$ does not depend on $x$. Below, we  fill this important gap assuming that, for all $x\in\mathcal{X}$, the distribution $P_{Y|x}$ is dominated by some $\sigma$-finite measure $\mu(\dd y)$. 

We start by stating a theorem  that provides sufficient conditions for \eqref{eq:cond_embed} to hold which only involve the Radon-Nikodym derivatives of $(P_{Y|x})_{x\in\mathcal{X}}$  w.r.t.~$\mu(\dd y)$.
 
\begin{theorem}\label{thm:Kernel_Gen2}
Assume that Assumptions \ref{asm:k}-\ref{asm:kxy} hold  and that there exists a $\sigma$-finite measure $\mu(\dd y)$   on $(\setY,\mathfrak{S}_{\mathcal{Y}})$ such that  $P_{Y|x}=p(y|x) \mu(\dd y)$ for all $x\in\mathcal{X}$, where $p(\cdot|\cdot)$ is such that the following conditions hold:
\begin{enumerate}
\item\label{CC1} We have $p(y|\cdot)\in\mathcal{H}_{\mathcal{X}}$ for all $y\in\mathcal{Y}$.
\item\label{CC2} The function $\mathcal{Y}\ni y\mapsto p(y|\cdot)$ is Borel measurable.
\item\label{CC3} For all $y'\in\mathcal{Y}$  the set $\big\{k_{\mathcal{Y}}(y',y)p(y|,\cdot),\,y\in\mathcal{Y}\}$ is separable.
\item\label{CC4}  We have $\int_{\mathcal{Y}} \|p(y|\cdot)\|_{\mathcal{H}_{\mathcal{X}}}\mu(\dd y)<\infty$.
\end{enumerate}
Then, condition \eqref{eq:cond_embed} holds.  
\end{theorem}

If $\mathcal{Y}$ is a finite set then Conditions \ref{CC2}-\ref{CC3} of   Theorem \ref{thm:Kernel_Gen2} always hold while Condition \ref{CC4} is implied by Condition \ref{CC1}. Hence, in this case, assuming Conditions \ref{CC1}-\ref{CC4} reduces to assuming  Condition  \ref{CC1}, which is both sufficient and necessary   for \eqref{eq:cond_embed} to hold when $\mathcal{Y}$ is a finite set. When  $\mathcal{Y}$ is not finite the additional  Conditions \ref{CC2}-\ref{CC4} are used to show that, for all $g\in\mathcal{H}_{\mathcal{Y}}$, the function $y\mapsto g(y)p(y|\cdot)$ is Bochner integrable and thus that $\E_{Y\sim P_{Y|\cdot}}\{g(Y)\}$ is a well-defined function  on $\mathcal{X}$.

The conclusions of Lemma \ref{lemma:CBM_New} and Theorem \ref{thm:Kernel_Gen2} are summarized in the following result:

\begin{corollary}\label{cor:Kernel_Gen2}
Consider the set-up of Theorem \ref{thm:Kernel_Gen2} and assume that $f\in\mathcal{H}_{\mathcal{X}}$ is such that $\E_{X\sim P_X}\{f^2(X)\}=0$ if and only $f=0$. Then, there exists a bounded conditional mean embedding operator $\mathcal{C}_{Y|X}$ for $(P_{Y|x})_{x\in\mathcal{X}}$ which is such that    $\|\mathcal{C}_{Y|X}\|_{\mathrm{o}}\leq \int_{\mathcal{Y}} \|p(y|\cdot)\|_{\mathcal{H}_{\mathcal{X}}}\mu(\dd y)$. 
\end{corollary}

In general, RKHS norms are hard to interpret.  However, if $k_{\mathcal{X}}$ is a Mat\'ern kernel and   $\mathcal{X}\subset\R^d$ is a bounded set with Lipschitz boundary,  e.g.~$\mathcal{X}$ is a hypercube,  then the RKHS $(\mathcal{H}_{\mathcal{X}}, <\cdot,\cdot>_{\mathcal{H}_{\mathcal{X}}})$  is norm equivalent to a Sobolev space \citep[see e.g.][Example 2.6]{kanagawa2018gaussian}.  As   shown  in     Appendix \ref{p_corr_usefuke}, together with Theorem \ref{thm:Kernel_Gen2}  this property of Mat\'ern kernels allows   to obtain  explicit conditions on the Radon-Nikodym derivatives $(p(\cdot| x))_{x\in\mathcal{X}}$ of $(P_{Y|x})_{x\in\mathcal{X}}$ which are sufficient to ensure that, under mild assumptions on $k_{\mathcal{Y}}$, condition \eqref{eq:cond_embed} holds. This allows us to establish   the following proposition, which provides for various definitions of $\mathcal{Y}$ non-trivial  examples of distributions $(P_{Y|x})_{x\in\mathcal{X}}$  that admit a bounded conditional mean embedding operator with respect to a characteristic kernel $k$. 
 
\begin{proposition}\label{prop:CME:example}
Let $\mathcal{X}\subset\R^d$ be a bounded set with Lipschitz boundary and strictly positive Lebesgue measure,  $k_{\mathcal{X}}$ be the restriction of the Mat\'ern kernel $K_{\frac{m}{2},\gamma}$ on $\mathcal{X}\times \mathcal{X}$, for some $m\in\mathbb{N}$ and $\gamma>0$, and let $k_{\mathcal{Y}}$ be a continuous, translation invariant, bounded and characteristic kernel on $\mathcal{Y}$. Then, $k=k_{\mathcal{X}}\times k_{\mathcal{Y}}$ is characteristic and there exists a bounded conditional mean embedding operator for $(P_{Y|x})_{x\in\mathcal{X}}$ when, for all $x\in\mathcal{X}$,
\begin{itemize}
\item  $P_{Y|x}=\sum_{m=1}^M w_m\mathcal{N}_1(\beta_m^\top x,\sigma^2_m)$ for some $M\in\mathbb{N}$, non-negative real numbers $\{w_m\}_{m=1}^M$ such that $\sum_{m=1}^M w_m=1$, vectors  $\{\beta_m\}_{m=1}^M$ in $\R^d$ and strictly positive real numbers $\{\sigma_m\}_{m=1}^M$, so that $\mathcal{Y}=\R$.

\item $P_{Y|x}=\mathcal{P}ois\{\exp(\beta^\top x)\}$ for some $\beta\in\R^d$, so that  $\mathcal{Y}=\mathbb{N}_0$.

\item $P_{Y|x}=\mathcal{B}er[1/\{1+\exp(\beta^\top x)\}]$ for some $\beta\in\R^d$, so that $\mathcal{Y}=\{0,1\}$.

\item $P_{Y|x}=\mathcal{G}amma\{\nu,\nu\exp(-\beta^\top x)\}$ for some $\beta\in\R^d$ and $\nu\in(0,\infty)$, so that $\mathcal{Y}=(0,\infty)$.

\item $P_{Y|x}$ is the distribution of $(Y_{x,1},Y_{x,2})$, where $Y_{x,2}=\ind_{(0,\infty)}(Y_{x,2}^*)$ and $Y_{x,1}=Y_{x,2}Y_{x,1}^*$ with
\begin{align*}
 \begin{pmatrix}
 Y_{x,1}^*\\
  Y^*_{x, 2}
  \end{pmatrix}
  \sim \mathcal{N}_2 
\left\{\begin{pmatrix}
\beta^\top x\\
\gamma^\top x
\end{pmatrix},
\begin{pmatrix}
\sigma^2&\rho\sigma\\
\rho\sigma&1
\end{pmatrix}\right\}
\end{align*}
for some $\beta,\gamma\in\R^d$, $\sigma>0$ and $\rho\in(-1,1)$, so that $\mathcal{Y}=\R\times\{0,1\}$.
\end{itemize}
The above assumptions on $k_{\mathcal{Y}}$ are satisfied if $k_{\mathcal{Y}}$ is the restriction on $\setY\times\setY$ of a  Mat\'ern kernel on $\R^{d_y}$, with $d_y=1$ for the first four definitions of $(P_{Y|x})_{x\in\mathcal{X}}$ and with $d_y=2$ for the last one.
\end{proposition}

\begin{remark}
The assumptions on $\mathcal{X}$ are satisfied e.g.~when this set is a non-empty open-hypercube.
\end{remark}

\section{MMD-based regression}\label{section:notations}

\subsection{ Set-up}

We let  $\{P_{\lambda},\lambda \in \Lambda\}$ be a set of probability distributions on $(\mathcal{Y},\mathfrak{S}_\mathcal{Y})$, $\Theta$ be a Polish space  and $g:\Theta\times \mathcal{X} \rightarrow \Lambda$ be such that the mapping $x\mapsto P_{g(\theta,x)}(A)$ is $\mathfrak{S}_\mathcal{X}$-measurable for all $A\in \mathfrak{S}_\mathcal{Y}$ and all $\theta\in\Theta$. Then, given a $\mathcal{Z}$-valued random variable $(X,Y)$, with $Y$ taking values in $\setY$,  we consider the statistical model $\{ ( P_{g(\theta,x)})_{x\in\setX},\,\theta\in\Theta\}$ for   the conditional distribution of $Y$ given $X$. For example, the Gaussian linear regression model with known variance is obtained by taking $P_\lambda = \mathcal{N}_1(\lambda,\sigma^2)$ and $g(\theta,x)=\theta^\top x$, the logistic regression model by taking $P_\lambda = \mathcal{B}er(\lambda)$  and $g(\theta,x)=1/\{1+\exp(-\theta^\top x)\}$, and the Poisson regression model  by taking $P_\lambda=\mathcal{P}ois(\lambda)$ and $ g(\theta,x) = \exp(\theta^\top x)$. Other classical examples include binomial, exponential, gamma and inverse-Gaussian regression models. 

In the following, $k$ is a kernel on $\mathcal{Z}$ satisfying Assumption \ref{asm:k}, stated in Section \ref{sec:kernel},  and  $D_n=\{(X_i,Y_i)\}_{i=1}^n$  is a   set of $n$ random variables   taking values in $\mathcal{Z}$. Below we assume that a realization of $D_n$ is used to fit the regression model $\{ ( P_{g(\theta,x)})_{x\in\setX},\,\theta\in\Theta\}$  and, for this reason,  we will often refer to $D_n$ as the observations.

\subsection{Definition of the estimators  \texorpdfstring{$\hat{\theta}_n$}{Lg} and \texorpdfstring{$\tilde{\theta}_n$}{Lg}}\label{sec:est}

Let   $\hat{P}^{n} = (1/n)\sum_{i=1}^{n} \delta_{(X_i,Y_i)}$ be the empirical distribution of the observations  and,   for every $\theta\in\Theta$,   let $\hat{P}_\theta^{n}$ be the (random) probability distribution on $(\mathcal{Z},\mathfrak{S}_\mathcal{Z})$ defined by  
\begin{align}\label{eq:model}
\hat{P}_\theta^{n}\big(A\times B\big)=\frac{1}{n} \sum_{i=1}^{n} \delta_{X_i }(A)P_{g(\theta,X_i)}(B),\quad  A\times B\in \mathfrak{S}_\mathcal{X}\otimes \mathfrak{S}_\mathcal{Y}.
\end{align}
In other words, if $(X,Y)\sim  \hat{P}_\theta^{n}$ then $X$ is uniformly distributed on the set $\{X_1,\dots,X_n\}$ and $Y|(X=x)\sim  P_{g(\theta,x)}$.

The first  estimator introduced in this work, $\hat{\theta}_n$, is  defined through the minimization of the MMD   between the probability  distributions $\hat{P}_\theta^n$ and $\hat{P}^{n}$, that is\footnote{When such a minimizer does not exist, we can use an $\epsilon$-minimizer instead and that follows can be trivially adapted.  In addition, we implicitly assume that $\hat{\theta}_n$ and $\tilde{\theta}_n$ are measurable, for all $n\geq 1$.} 
\begin{equation}\label{dfn:estim}
\begin{split}
\hat{\theta}_n(D_n) &\in  \argmin_{\theta\in\Theta} \mathbb{D}_k^2\big( \hat{P}^n_\theta, \hat{P}^n)= \argmin_{\theta\in\Theta} \hat{F}_n(\theta),\quad   \hat{F}_n(\theta):=\sum_{ i,j=1}^n \hat{\ell}(\theta,X_i,X_j, Y_j)
\end{split}
\end{equation}
where,  for all $\theta\in\Theta$, $(x,x')\in\setX^2$ and $y\in\setY$,
\begin{equation}\label{dfn:estim_l}
\hat{\ell}(\theta, x,x',y)=\E_{  Y\sim P_{g(\theta,x)},  \,\,Y'\sim P_{g(\theta,x') }} \big[k\big\{(x,Y),(x',Y')\big\} - 2  k\big\{(x,Y),(x',y)\big\}\big].
\end{equation}
Under Assumption \ref{asm:k},  the kernel $k$ is bounded so that each   term in the double sum appearing in the definition of $\hat{F}_n$  is bounded by $3$. Intuitively, this limits the impact that a single observation can have on $\hat{\theta}_n(D_n)$, making this estimator robust to outliers.
 
The number of terms in the definition of the function $\hat{F}_n$ to minimize is $\mathcal{O}(n^2)$. Below we propose an approach that makes possible to efficiently compute  $\hat{\theta}_n$ for moderate values of $n$, that is for $n$ equals to a few thousands, but this feature of   $\hat{F}_n$ limits the applicability of  $\hat{\theta}_n$  in large datasets. For large scale problems we propose the alternative  estimator $\tilde{\theta}_n$, defined by 
\begin{equation}
\label{eq:estim:approx}
 \tilde{\theta}_n(D_n) \in  \argmin_{\theta\in\Theta}\tilde{F}_n(\theta),\quad   \tilde{F}_n(\theta):=\sum_{i=1}^n \tilde{\ell}(\theta,X_i,Y_i)
\end{equation}
where,  for all $\theta\in\Theta$, $x\in\setX$ and $y\in\setY$,
\begin{equation}\label{eq:estim:approx_l}
\tilde{\ell}(\theta, x,y)=\E_{\substack{ Y,Y'\,\,\iid P_{g(\theta,x)}}} \big\{k_{\mathcal{Y}}(Y,Y')  - 2 k_{\mathcal{Y}}(Y,y)\big\}.
\end{equation}

The function $\tilde{F}_n$ defined in \eqref{eq:estim:approx} has the advantage to  involve  only $n$ terms but, on the other hand,     $\tilde{\theta}_n$  cannot be interpreted as the minimizer of a measure of discrepancy between  $\hat{P}_\theta^n$ and $\hat{P}^{n}$. 
Theoretical results regarding the robustness properties of  $\tilde{\theta}_n$ are provided in the next section,   they are weaker than those obtained for $\hat{\theta}_n$.

Henceforth we  use the shorthand  $\hat{\theta}_n = \hat{\theta}_n(D_n)$ and $\tilde{\theta}_n = \tilde{\theta}_n(D_n)$, which is standard in statistics.

\subsection{Link between the two estimators}\label{sub:link}

In this subsection we assume that $k=k_{\mathcal{X}}\otimes k_{\mathcal{Y}}$ where $k_{\mathcal{X}}=k_{\gamma}$ for some   kernel $k_\gamma$ on $\mathcal{X}$ such  that $k_{\gamma}(x,x)=1$ and such that $\lim_{\gamma \rightarrow 0}k_{\gamma}(x,x')=0$ for all $x'\neq x$. When $\mathcal{X}\subseteq\R^d$ these two properties are for instance satisfied   when, for some $\alpha\in(0,\infty]$, the kernel $k_\gamma$ is a Mat\'ern kernel, that is when $k_\gamma=K_{\alpha,\gamma}$ with $K_{\alpha,\gamma}$ as introduced in   Section \ref{sub:notation}. 

Let   $\ell(\theta,x,x',y) =\hat{\ell}(\theta, x,x',y)/k_{\gamma}(x,x')$ if $k_{\gamma}(x,x')\neq 0$ and $\ell(\theta,x,x',y)=0$ otherwise. Under the above assumptions on $k$ the quantity $\ell(\theta,x,x',y)$ does not depend on   $k_{\gamma}$ and is such that $\ell(\theta,x,x,y)=\tilde{\ell}(\theta,x,y)$. Therefore,  letting
\begin{align}\label{eq:Hn}
h_n\big(\gamma,\theta,D_n\big)=2\sum_{i=1}^{n-1}\sum_{j=i+1}^n k_{\gamma}(X_i,X_j)\ell(\theta, X_i, X_j,Y_j),
\end{align}
it follows that the  estimators $\hat{\theta}_n$ and $\tilde{\theta}_n$ are such that
\begin{align*}
&\hat{\theta}_n\in\argmin_{\theta\in\Theta}\bigg\{ \sum_{i=1}^n \tilde{\ell}(\theta, X_i, Y_i) + h_n\big(\gamma,\theta,D_n\big)\bigg\},\quad \tilde{\theta}_n\in\argmin_{\theta\in\Theta} \sum_{i=1}^n\tilde{\ell}(\theta,  X_i,Y_i).
\end{align*}

Consequently, using $\tilde{\theta}_n$ in place of $\hat{\theta}_n$ amounts to discarding,  in the definition of this latter estimator, the  term $h_n(\gamma,\theta, D_n)$ whose computation  requires $\mathcal{O}(n^2)$ operations. Assuming that the $X_i$'s are $\P$-a.s.\ distinct, under the above assumptions on $k_\gamma$ for all $\theta\in\Theta$ we have    $ \lim_{\gamma\rightarrow 0}h_n(\gamma,\theta, D_n)=0$, $\P$-a.s. Therefore, under   suitable continuity assumptions, $\hat{\theta}_n \rightarrow \tilde{\theta}_n$ as $\gamma  \rightarrow 0$, $\P$-a.s. For this reason,  and as illustrated in Section \ref{section:simulation}, for a small value of $\gamma $ we expect the two estimators $\hat{\theta}_n$ and $\tilde{\theta}_n$   to have a very similar behaviour in practice.  

\subsection{Computation of the two estimators}
\label{subsec:algo}

Computing the estimators  $\hat{\theta}_n$ and $\tilde{\theta}_n$  require   optimizing the functions $\hat{F}_n$ and $\tilde{F}_n$, respectively, which are both defined through an expectation. In some models  the expectations   appearing in \eqref{dfn:estim_l} and in \eqref{eq:estim:approx_l}  can be computed explicitly, in which case the functions  $\hat{F}_n$ and $\tilde{F}_n$  has a known expression and standard optimization algorithms can be used to    optimize them. This is for example the case in logistic or multinomial regression, since for these two models the expectations in \eqref{dfn:estim_l} and in  \eqref{eq:estim:approx_l}  are   finite sums.

In general, the functions  $\hat{F}_n$ and  $\tilde{F}_n$ are however intractable, and as a general strategy for computing $\hat{\theta}_n$ and $\tilde{\theta}_n$ we propose the use of a stochastic gradient   algorithm. As shown in Appendix \ref{subsection:proof:gradient1}, under suitable regularity conditions we have
\begin{align*}
&\nabla_\theta   \hat{\ell} (\theta, x,x',y)\\
 &=2 \E_{Y\sim P_{g(\theta,x)}, \,\,Y'\sim P_{g(\theta,x') }}   \Big(\big[k\big\{(x,Y),(x',Y')\big\} -  k\big\{(x,Y),(x',y)\big\}\big] \nabla_{\theta}  \log p_{g(\theta,x)}(Y)  \Big)
 \end{align*}
 and
 $$
\nabla_\theta \tilde{\ell}(\theta, x,y)=2\E_{\substack{ Y,Y'\,\,\iid P_{g(\theta,x)}}}\big[  \big\{k_{\mathcal{Y}}(Y,Y') - k_{\mathcal{Y}}(Y,y) \big\} \nabla_{\theta}  \log p_{g(\theta,x)}(Y)   \big] 
$$
so that we can easily compute an unbiased estimate of $\nabla \hat{F}_n(\theta)$ and of $\nabla \tilde{F}_n(\theta)$ under the mild conditions that (i) we can sample from  $P_\lambda$ for all $\lambda\in\Lambda$ and  (ii) that we can compute the gradient of the log-likelihood function of a single observation. The stability of this procedure obviously depends on the model, but it holds generally at least for compact parameter spaces $\Theta$, as discussed in Appendix \ref{subsection:proof:gradient1}.

It is important to stress that the functions $\hat{F}_n$ and $\tilde{F}_n$ being typically non-convex and potentially multi-modal, notably   in the presence of outliers, the starting value of the stochastic gradient algorithm must be chosen with care. When the maximum likelihood estimator  of the model parameter can be efficiently computed we recommend to use its value to initialize the optimization procedure. Otherwise, the gradient-free algorithm introduced in \citet{gerber2022global}, designed to compute the global optimum of a function  defined through an expectation, can be used to find a good starting value for the stochastic gradient algorithm.

It is important to mention at this stage that the computation of $\hat{\theta}_n$ can be greatly facilitated by taking $k=k_{\mathcal{X}}\otimes k_{\mathcal{Y}}$  where $k_{\mathcal{X}}=k_\gamma$  with $k_{\gamma}$ as in Section \ref{sub:link}. Indeed, for such a kernel $k$  it often true that with high probability we have $k_{\gamma}(X_i,X_j)\approx 0$ for all $i\neq j$, and thus that $h_n\big(\gamma,\theta,D_n\big)\approx 0$,  where $h_n\big(\gamma,\theta,D_n\big)$ is as defined in \eqref{eq:Hn}.  In this case, we can efficiently compute $\hat{\theta}_n$ with a stochastic gradient algorithm whose cost per iteration  is linear in the sample size $n$, as explained  in   Appendix \ref{subsection:proof:gradient2}.

\section{Theoretical analysis }
\label{section:theory}

\subsection{Set-up and summary of the main results}

In this section we introduce theoretical results concerning the robustness of the   estimators $\hat{\theta}_n$ and $\tilde{\theta}_n$ when outliers are present in the dataset  used by the statistician   to fit the regression model. To this aim, we need to interpret $D_n$ as a contaminated version of a dataset $D^0_n:=\{(X^0_i,Y^0_i)\}_{i=1}^n$, and below   we consider two contamination models, namely the Huber contamination model \citep{huber1964robust} and the adversarial  contamination model, defined in Definition \ref{def:adv} and in Definition \ref{def:Huber}, respectively. 

\begin{definition}\label{def:adv}
We say that the observed dataset $D_n$ is an $\epsilon$-adversarial contamination of $D_n^0$ if there exists a set  $I\subset\{1,\dots,n\}$ such that  $|I|/n\leq \epsilon$ and such that $(X_i,Y_i)=(X^0_i,Y^0_i)$ for all $i\not\in I$.
\end{definition}

\begin{definition}\label{def:Huber}
We say that the observed dataset $D_n$ is an $(\epsilon,Q)$-Huber contamination of $D_n^0$ if there exists a $Q\in\mathcal{P}(\mathcal{Z})$ such that $(X_i,Y_i)\sim \epsilon Q+(1-\epsilon)\delta_{(X_i^0,Y_i^0)}$  for all $i\in\{1,\dots,n\}$.
\end{definition}

Below we   derive   non-asymptotic  bounds for  $\hat{\theta}_n$ under both  the fixed design and the  random design  scenarios, which prove the robustness of this estimator to adversarial contaminations of the data. For $\tilde{\theta}_n$  we derive an asymptotic result in the random design case which establishes the robustness of this estimator to Huber contaminations of the sample.

Recall that results in the fixed design case only provide guarantees on the estimation of the distribution of $Y$ when $X$ is equal to one of the observed $X_i$'s, while in practice regression is often used for out-of-sample predictions. In this case, and assuming that the pairs $(Y^0_i,X^0_i)$'s are i.i.d, this means that we want   guarantees on the estimation of the distribution of $Y$ when $X$ is drawn from the same, unknown, distribution as the observed $X_i$'s, and independently from them. This is precisely what  theoretical results in the random design case provide.

In what follows we let $(P^0_{Y|x})_{x\in\setX}$ be a regular conditional probability of $Y$ given $X$, and thus $Y^0_i|(X^0_i=x)\sim P^0_{Y|x}$ for all $x\in\setX$ and all $i=1,\dots,n$.

\subsection{Non-asymptotic bounds for the estimator \texorpdfstring{$\hat{\theta}_n$}{Lg}--Fixed design }

A typical scenario for the fixed design case is when the $X^0_i$'s are experimental settings that are carefully planned in advance. In this case, measurement errors can only affect the $Y^0_i$'s. For this reason, in this subsection we  assume that the contamination of the sample occurs only on the $Y^0_i$'s, so   that $X_i = X^0_i$ for all $i\in\{1,\dots,n\}$.

Letting
$$
\bar{P}^{0}_n(A\times B) = \frac{1}{n}\sum_{i=1}^n \delta_{X^0_i}(A) P_{Y|X^0_i}^0(B),\quad \forall (A\times B)\in \mathfrak{S}_\mathcal{Z},
$$
we set up our objective as the reconstruction   of   $\bar{P}^{0}_n\in\mathcal{P}(\mathcal{Z})$ by a distribution in the set $\{\hat{P}_\theta^n, \theta\in\Theta\}$.

The following lemma  gives  a non-asymptotic bound on the performances of the estimator  $\hat{\theta}_n$ for this task, under   an adversarial contamination of the sample.

\begin{lemma}\label{thm:fixed1}
Assume that $D_n$ is an $\epsilon$-adversarial contamination of $D_n^0$ such that $X_i=X_i^0$ for all $i\in\{1,\dots,n\}$. Then, under Assumption~\ref{asm:k},  
$$
\mathbb{E}\big\{  \mathbb{D}_k(\hat{P}^n_{\hat{\theta}_n} ,\bar{P}^{0}_n ) \big\}  \leq 4\epsilon + \inf_{\theta\in\Theta}  \mathbb{D}_k(\hat{P}^n_{\theta} ,\bar{P}^{0}_n ) +  2/\sqrt{n}
$$
and, for all $\eta\in(0,1)$,
 \begin{align}
 \label{eq:tthm:fixed1:p}
 \mathbb{P}\Big[\mathbb{D}_k(\hat{P}^n_{\hat{\theta}_n} ,\bar{P}^{0}_n )  & < 4\epsilon + \inf_{\theta\in\Theta}  \mathbb{D}_k(\hat{P}^n_{\theta} ,\bar{P}^{0}_n) +  n^{-1/2}\big\{2+\sqrt{ 2\log(1/\eta)}\big\} \Big] \geq 1-\eta.
 \end{align}
\end{lemma}
\begin{remark}
Lemma \ref{thm:fixed1}   does not require  any  assumption on the distribution of  $\{(X_i^0,Y_i^0)\}_{i=1}^n$.
\end{remark}

In statistical theory  we often   assume that the ``truth is in the model'', that is  that there is a $\theta_0\in\Theta$ such that $\hat{P}^n_{\theta_0} = \bar{P}^0_n$. In this case, Lemma \ref{thm:fixed1} shows that
\begin{align}\label{eq:cont1}
\mathbb{E}\big\{  \mathbb{D}_k(\hat{P}^n_{\hat{\theta}_n} ,\hat{P}^n_{\theta_0} )  \big\}\leq 4\epsilon + 2/\sqrt{n}
\end{align}
while, in the non-contaminated case where $D_n = D_n^0$,
\begin{align}\label{eq:cont2}
\mathbb{E} \big\{  \mathbb{D}_k(\hat{P}^n_{\hat{\theta}_n} ,\bar{P}^{0}_n ) \big\}  \leq  2/\sqrt{n}.
\end{align}
In  words, when  computed from the uncontaminated dataset $D_n^0$ the  estimator $\hat{\theta}_n$ is consistent for estimating $\theta_0$, with respect to the MMD distance, and,   provided that  $\epsilon$ is small, an $\epsilon$-adversarial contamination of the sample will have only a negligible  impact on the estimated parameter value. Similar conclusions can    be derived from the inequality in probability given in \eqref{eq:tthm:fixed1:p}.

Lemma~\ref{thm:fixed1}  implies  the convergence of $\hat{\theta}_n$ with respect to the MMD distance. However, under additional assumptions, it is possible to relate this form of convergence to the convergence in the sense of the Euclidean distance $\|\cdot\|$ to the true parameter, or the pseudo-true parameter.

 \begin{theorem}
\label{cor:fixed}
Consider the set-up of Lemma~\ref{thm:fixed1} and assume that there is a unique $\theta_0\in\Theta$ such that $\mathbb{D}_k(\hat{P}^n_{\theta_0} ,\bar{P}^{0}_n ) = \inf_{\theta\in\Theta} \mathbb{D}_k(\hat{P}^n_{\theta} ,\bar{P}^{0}_n )$. In addition,  assume that there exist a neighbourhood $U$ of $\theta_0$ and a constant $\mu>0$ such that
$\mathbb{D}_k(\hat{P}^n_{\theta},\hat{P}^n_{\theta_0} )  \geq \mu\|\theta-\theta_0\|$ for all $\theta\in U$. 
Let $\alpha = \inf_{\theta\in U^c}  \mathbb{D}_k(\hat{P}^n_{\theta}, \hat{P}^n_{\theta_0}) \in (0,2]$ and assume that $\epsilon \in[0, \alpha /32)$.
Then, 
$$
\mathbb{P}\big(\lim\sup_{n\rightarrow\infty } \|\hat{\theta}_n-\theta_0\| \leq 4\epsilon/\mu\big)=1
$$
 and, for all $n\geq 64/\alpha^2$ and all $\eta \in [ 2e^{-n \alpha^2/38},1)$,
$$
\mathbb{P}\Big[  \|\hat{\theta}_n-\theta_0\| < (4\epsilon/\mu)+ n^{-1/2} \big\{2+\sqrt{2\log(2/\eta)}\big\}/\mu \Big] \geq 1-\eta.
$$
 \end{theorem}

\subsection{Non-asymptotic bounds for the estimator   \texorpdfstring{$\hat{\theta}_n$}{Lg}--Random design}\label{subsec:random_hat}

We   assume now that the pairs $(X^0_i,Y^0_i)$'s are i.i.d.~from some probability distribution $P^0\in\mathcal{P}(\mathcal{Z})$, and we denote by $P_X^0$  the marginal distribution of the $X^0_i$'s. In addition, for every $\theta\in\Theta$  we let $P_\theta\in\mathcal{P}(\mathcal{Z})$ be defined by
\begin{align*}
P_\theta(A\times B)=\E_{X\sim P^0_X}\big\{\ind_A(X) P_{g(\theta,X)}(B)\},\quad  A\times B\in \mathfrak{S}_\mathcal{Z}.
\end{align*}
Then, we set up our objective as the reconstruction of $P^0$ by a distribution in $\{P_\theta,\,\theta\in\Theta\}$, and our natural candidate is $P_{\hat{\theta}_n}$.

Since the approximating set $\{P_\theta,\,\theta\in\Theta\}$ is unknown, because it depends on $P_X^0$,   achieving this objective requires  more care  than in the fixed design setting. In particular, it requires  to control the MMD  distance between the distribution $P_\theta$ and its empirical counterpart $\hat{P}_\theta^n$ defined in \eqref{eq:model}, a task that we perform using   Lemma \ref{lemma:constant:kernel}. For this reason, the results presented below assume  that $k=k_{\mathcal{X}}\otimes k_{\mathcal{Y}}$ and require  the following additional assumption:  
\begin{assumption}
 \label{asm:cond:expect}
 For all  $\theta\in\Theta$  there exists a bounded linear conditional mean embedding operator $\mathcal{C}_{\theta}:\mathcal{H}_{\mathcal{X}}\rightarrow\mathcal{H}_{\mathcal{Y}}$ for $(P_{g(\theta,x)})_{x\in\mathcal{X}}$. In addition,   $\mathfrak{C}:=\sup_{\theta\in\Theta}\|\mathcal{C}_{\theta}\|_{\mathrm{o}}<\infty$.
\end{assumption}

Under this additional assumption we have the following result:

\begin{lemma}\label{thm:random1}
Assume that $D_n$ is an $\epsilon$-adversarial contamination of $D_n^0$. Then, under Assumptions~\ref{asm:k}-\ref{asm:cond:expect}, and with $\mathfrak{C}<\infty$ as in Assumption \ref{asm:cond:expect}, we have
$$
   \mathbb{E} \big\{\mathbb{D}_k \big( P_{\hat{\theta}_n} , P^0  \big)  \big\}  \leq 8\epsilon + 
 \inf_{\theta\in\Theta}  \mathbb{D}_k  ( P_{\theta} , P^0  )
 +
 (\mathfrak{C} + 3)/\sqrt{n}
 $$
and, for all $\eta\in (0,1)$,
 \begin{align*}
 \mathbb{P} & \Big[
  \mathbb{D}_k( P_{\hat{\theta}_n} , P^0 )
  <
  8\epsilon +  
  \inf_{\theta\in\Theta}  \mathbb{D}_k( P_{\theta} , P^0)
  + n^{-1/2}( \mathfrak{C}+3)\big\{1+\sqrt{2\log(4/\eta) }\big\} 
 \Big] \geq 1-\eta.
 \end{align*}
\end{lemma}

From Lemma \ref{thm:random1} we can readily obtain the random design counterpart of the inequalities \eqref{eq:cont1}-\eqref{eq:cont2}, obtained in the fixed design setting, to prove the consistency of $\hat{\theta}_n$ in the well-specified case and in the absence of contamination.

The following theorem is the main result of this subsection.
 \begin{theorem}
 \label{cor:random}

Assume that there is a unique $\theta_0\in\Theta$ such that $\mathbb{D}_k(P_{\theta_0},P^0 )=\inf_{\theta\in\Theta}\mathbb{D}_k(P_{\theta},P^0 )$ and  that there exist a neighbourhood $U$ of $\theta_0$ and a constant $\mu>0$ such that
\begin{align}\label{eq:loweDk}
\mathbb{D}_k(P_{\theta},P^0  ) \geq \mathbb{D}_k(P_{\theta_0},P^0  ) +\mu\|\theta-\theta_0\|,\quad\forall \theta\in U.
\end{align}
Let  $\alpha = \inf_{\theta\in U^c}  \mathbb{D}_k(P_{\theta},P^0)- \mathbb{D}_k(P_{\theta_0},P^0  ) \in (0,2]$ and  assume that $D_n$ is an $\epsilon$-adversarial contamination of $D_n^0$ for some $\epsilon\in[0, \alpha/64)$.  Assume also that Assumptions~\ref{asm:k}-\ref{asm:cond:expect} hold and let  $\mathfrak{C}<\infty$ be as in Assumption \ref{asm:cond:expect}. Then, 
$$
\mathbb{P}\big(\limsup_{n\rightarrow\infty } \|\hat{\theta}_n-\theta_0\| \leq 8 \epsilon/\mu\big)=1
$$
and there exist     constants $C_1,C_2\in(0,\infty)^2$, that depend only on $\alpha$ and on $\mathfrak{C}$, such that for all $ n\geq C_1$ and all  $\eta \in [ 8e^{- C_2 n},1)$ we have
$$
\mathbb{P}\Big[ \|\hat{\theta}_n-\theta_0\| < (8\epsilon/\mu)+ n^{-1/2} (\mathfrak{C}+3) \big\{1+\sqrt{2\log(8/\eta)}\big\} /\mu \Big] \geq 1-\eta.
$$

\end{theorem}

\begin{remark}
If $P^0=P_{\theta_0}$  for some $\theta_0\in\Theta$, i.e.~if the model is well-specified, then a sufficient condition for $\theta_0$ to be the unique global minimizer of the mapping $\theta\mapsto \mathbb{D}_k(P_{\theta}, P^0  )$ is that $k$ is characteristic and the model $\{P_\theta,\,\,\theta\in\Theta\}$ is identifiable, in the sense that $\theta_1\neq\theta_2\Rightarrow P_{\theta_1}\neq P_{\theta_2}$.
\end{remark}

In Theorem \ref{cor:random} the distribution $P_{\theta_0}$ should be interpreted as the best approximation of $P^0$ in the sense of the MMD distance $\mathbb{D}_k$. It is worth noting that, unless the model is well-specified, in which case   $P_{\theta_0}=P^0$, both the parameter value $\theta_0$ and the distribution $P_{\theta_0}$ depend on the  choice of $k$. If a small fraction $\epsilon$ of the data is corrupted then the theorem ensures that $\hat{\theta}_n$  still estimates well $\theta_0$.  In particular, the first part of Theorem \ref{cor:random} implies that the influence function of the estimator $\hat{\theta}_n$ is bounded. It is also worth noting that taking $\epsilon=0$ in Theorem \ref{cor:random} establishes  the almost sure convergence of $\hat{\theta}_n$ towards $\theta_0$.

Condition \eqref{eq:loweDk} of Theorem \ref{cor:random} requires that the function
$
\theta\mapsto \mathbb{D}_k(P_{\theta}, P^0   )- \mathbb{D}_k(P_{\theta_0}, P^0  )
$
is strongly convex in a neighbourhood $U$ of $\theta_0$, a condition which is rather weak, as shown in the next proposition.

\begin{proposition}\label{prop:Taylor}
Assume that $\theta_0\in\argmin_{\theta\in\Theta}\mathbb{D}_k(P_{\theta},P^0 )$ is unique  and that the function   $\theta\mapsto \mathbb{D}_k(P_{\theta},P^0)$ is twice continuously differentiable at $\theta_0$. Then, condition \eqref{eq:loweDk} of Theorem \ref{cor:random} holds.
\end{proposition}

Proposition \ref{prop:CME:example}, given in Section \ref{sub:cond_op}, provides examples of  characteristic kernels $k$ for which, for all $\theta\in\Theta$, a conditional mean embedding operator $\mathcal{C}_\theta$ of $(P_{g(\theta,x)})_{x\in\mathcal{X}}$ exists for five popular regression models, namely for (i)  the linear Gaussian regression model, first example of Proposition \ref{prop:CME:example}  with $\theta=(\beta,\sigma)$, as well as for mixtures of such models, (ii) the Poisson regression model,    second example of Proposition \ref{prop:CME:example}  with $\theta=\beta$, (iii) the logistic regression model, third example of Proposition \ref{prop:CME:example}  with $\theta=\beta$, (iv)  the Gamma regression model, fourth example of Proposition \ref{prop:CME:example}  with $\theta=(\beta,\nu)$, and (v)  the Heckman sample selection model, last example of Proposition \ref{prop:CME:example}  with $\theta=(\beta,\gamma,\sigma,\rho)$.

Using the last part of Corollary \ref{cor:Kernel_Gen2}, one can easily check that for each of these models, under a suitable definition of $\Theta$ we have  $\sup_{\theta\in\Theta}\|\mathcal{C}_\theta\|_{\mathrm{o}}<\infty$, as required by Assumption \ref{asm:cond:expect}. In the set-up of Proposition \ref{prop:CME:example}, this is for instance the case if, for some compact set $B\subset\R^d$ and some $\delta>0$, the parameter space $\Theta$ is such that $\beta\in B$ for all models, such that $\gamma\in B$ for the  Heckman sample selection model, such that  $\sigma>\delta$ for the linear Gaussian regression model or mixtures of such models, and for the  Heckman sample selection model, and such that  $\nu\in(\delta,1/\delta)$ for the Gamma regression model.

 \begin{remark}
 Proposition \ref{prop:CME:example} assumes that $\mathcal{X}\subset\R^d$ is a bounded set and, as illustrated with the above examples,  Assumption \ref{asm:cond:expect}   typically holds when $\Theta$ is a compact set.  It is however important to note that the dependence to the outliers of the bounds given in   Lemma \ref{thm:random1} and in Theorem \ref{cor:random}   depend neither on $\Theta$ nor on $\mathcal{X}$. Indeed, the outliers impact  these bounds only through their proportion  $\epsilon$.
 \end{remark}

\subsection{Asymptotic  guarantees for the estimator  \texorpdfstring{$\tilde{\theta}_n$}{Lg}--Random design}

For this estimator we set up our objective as the reconstruction of the regular conditional probability $(P^0_{Y|x})_{x\in\mathcal{X}}$ by a distribution in the set $\{ (P_{g(\theta,x)})_{x\in\setX},\,\theta\in\Theta\}$. As in Section \ref{subsec:random_hat}, in what follows  we denote by $P^0$   the   distribution of the $(X_i^0, Y_i^0)$'s. In addition,  for any distribution $Q\in\mathcal{P}(\mathcal{Z})$ we let $Q_X$ denote the distribution of $X$  and $(Q_{Y|x})_{x\in\setX}$ denote  a regular conditional probability for the distribution of $Y$ given $X$, where $(X,Y)\sim Q$.

It is direct to see that $\tilde{\theta}_n$ is an $M$-estimator   and therefore   sufficient conditions on $k_{\mathcal{Y}}$ and on the statistical model to ensure that $\tilde{\theta}_n$ converges to some value $\tilde{\theta}_{0}\in\Theta$ as $n\rightarrow\infty$ can be obtained from the general theory on $M$-estimators  \citep[][Chapter 5]{Vaart2000}. Using this approach, we easily obtain the following proposition: 

\begin{proposition}\label{prop:random2}
Let $k_{\mathcal{Y}}$ satisfy $|k_{\mathcal{Y}}|\leq 1$,   $D_n=D_n^0$  for all $n$, and assume  that the following conditions hold:  
\begin{itemize}
\item   $\Theta$ is  compact.
\item   The mapping  $
\theta\mapsto \mathbb{E}_{Y, Y'\iid P_{g(\theta,x)} } \big\{k_{\mathcal{Y}}(Y,Y') -2 k_{\mathcal{Y}}(Y,y)\big\}
$ is continuous on $\Theta$  for all $(x,y)\in\setX\times\setY$. 
\item    The mapping $\theta\mapsto \E\big\{ \mathbb{D}_{k_{\mathcal{Y}}}\big(P_{g(\theta,X_1)},  P^0_{Y|X_1}\big)^2\big\}$ has a unique global minimum at $\tilde{\theta}_{0}\in\Theta$.
\end{itemize} Then, $\tilde{\theta}_n\rightarrow \tilde{\theta}_{0}$ in $\mathbb{P}$-probability.
\end{proposition}

If the model is well-specified, that is if $(P^0_{Y|x})_{x\in\mathcal{X}}\in\{ (P_{g(\theta,x)})_{x\in\setX},\,\theta\in\Theta\}$, then a sufficient condition for $\tilde{\theta}_0$ to be well-defined is that $k_{\mathcal{Y}}$ is a characteristic kernel and that the model $\{ (P_{g(\theta,x)})_{x\in\setX},\,\theta\in\Theta\}$ is identifiable, in the sense that $\mathbb{P}\big(\theta_1\neq\theta_2\Rightarrow P_{g(\theta_1,X)}\neq P_{g(\theta_2,X)}\big)=1$.  It is also important to stress that, since   for $\tilde{\theta}_n$ we focus on Huber's type contaminations of the sample, in Proposition \ref{prop:random2} there is no loss of generality to assume that  $D_n=D_n^0$ for all $n$. Indeed, in the random design setting, under an  $(\epsilon,Q)$-Huber contamination of the sample we have that  $(X_i,Y_i)\iid \tilde{P}^0:=\epsilon Q+(1-\epsilon) P^0$ and thus, since   Proposition \ref{prop:random2} requires no assumption on   $P^0$, its result  remains valid if $P^0$ is replaced by $\tilde{P}^0$.  It is also worth mentioning that, unless the model is well-specified,    the parameter value $\tilde{\theta}_{0}$ defined in  Proposition \ref{prop:random2}  will be typically  different from the parameter value $\theta_0$  the estimator $\hat{\theta}_n$ converges to. Finally, we note that using the general theory on $M$-estimators   one can obtain sufficient conditions on $k_{\mathcal{Y}}$ and on the statistical model which ensure that $\tilde{\theta}_n$ is $\sqrt{n}$-consistent and asymptotically Gaussian.

The following theorem provides an asymptotic guarantee regarding  the robustness of  $\tilde{\theta}_n$  to Huber type contaminations of the data.  Notably, a direct implication of this theorem is that the influence function   of $\tilde{\theta}_n$ is bounded.

\begin{theorem}\label{thm:random2}
Let $k_{\mathcal{Y}}$ satisfy $|k_{\mathcal{Y}}|\leq 1$, $Q\in\mathcal{P}(\mathcal{Z})$ and assume that following two conditions hold:
\begin{itemize}
\item  The mapping  $\theta\mapsto \E\big\{ \mathbb{D}_{k_{\mathcal{Y}}}\big(P_{g(\theta,X^0_1)},  P^0_{Y|X^0_1}\big)^2\big\}$ has a unique global minimum at $\tilde{\theta}_0$.

\item  There exist a neighbourhood $U$ of $\tilde{\theta}_0$ and a constant $\mu>0$ such that, for all $\theta\in U$,
\begin{align}\label{eq:cond_tildeU}
\hspace{-0.3cm}\E_{X \sim P_X^0}\big\{ \mathbb{D}_{k_{\mathcal{Y}}}(P_{g(\theta,X)}, P_{Y|X}^0)^2\big\}\geq \E_{X \sim P_X^0}\big\{ \mathbb{D}_{k_{\mathcal{Y}}}(P_{g(\tilde{\theta}_0 ,X)}, P_{Y|X}^0)^2\big\}+\mu\|\theta-\tilde{\theta}_0\|. 
\end{align}
\end{itemize}
Let
$$
\alpha= \inf_{\theta\in U^c} \E_{X \sim P_X^0}\big\{ \mathbb{D}_{k_{\mathcal{Y}}}(P_{g(\theta,X)}, P_{Y|X}^0)^2\big\}-\E_{X \sim P_X^0}\big\{ \mathbb{D}_{k_{\mathcal{Y}}}(P_{g(\tilde{\theta}_0,X)}, P_{Y|X}^0)^2\big\}\in(0,4]
$$
and assume that $\epsilon\in [0, \alpha/(52+\alpha))$. Then,  for all
$$
\tilde{\theta}_{Q,\epsilon}\in \argmin_{\theta\in\Theta}\E_{X_1\sim\epsilon Q_X+(1-\epsilon) P^0_X}\Big[ \mathbb{D}_{k_{\mathcal{Y}}}\big\{P_{g(\theta,X_1)}, \epsilon Q_{Y|X_1}+(1-\epsilon)P^0_{Y|X_1}\big\}^2\Big]
$$
 we have $\|\tilde{\theta}_{Q,\epsilon}-\tilde{\theta}_0\|\leq   52 \epsilon/(\mu -\epsilon\mu)$.
\end{theorem}

Following similar steps as in the proof of Proposition \ref{prop:Taylor} it is easily checked  that condition \eqref{eq:cond_tildeU} holds when the function $\theta\mapsto \E_{X \sim P_X^0}\{ \mathbb{D}_{k_{\mathcal{Y}}}(P_{g(\theta,X)},P_{Y|X}^0)^2\}$ has $\tilde{\theta}_0$ as unique global minimizer and is twice continuously differentiable around  this parameter value.

\section{Numerical experiments}\label{section:simulation}

\subsection{Set-up\label{sub:sim_setup}}

All the results presented in this section are obtained for $\mathcal{X}=\R^d$, $\mathcal{Y}\subset\R^{d_y}$ for some $d_y\in\mathbb{N}$, and  for the kernel $k=k_{\mathcal{Y}} \otimes k_{\mathcal{X}}$ such that $k_{\mathcal{X}}=K_{0.5,0.01}$, i.e.~$k_{\mathcal{X}}$ is the exponential kernel on $\R^d$ with bandwidth parameter $\gamma=0.01$, and such that $k_{\mathcal{Y}}$ is the exponential kernel on $\mathcal{Y}$ with bandwidth parameter equal to 1.

For each  experiment   the  observations  used to fit the model  are obtained from an uncontaminated dataset $d^0_{N}:=\{(x_i^0,y_i^0)\}_{i=1}^N$   that we  contaminate  as follows. We choose an $\epsilon\in [0,1)$ and  randomly select a set $I\subset\{1,\dots,N\}$ such that $|I|=\lfloor \epsilon N \rfloor$. For all $i\in\{1,\dots,N\}$  we then let $(x_i, y_i)=(x^0_i, y^0_i)$ if $i\not\in I$ and  $(x_i, y_i)=(x^c_i, y^c_i)$ for some $(x^c_i, y^c_i)\in\mathcal{Z}$ if $i\in I$. The way we generate the  $(x^c_i, y^c_i)$'s will vary from one example to the next and will therefore be specified in due course. Finally,  for $n\leq N$ we let $d_n=\{(x_i,y_i)\}_{i=1}^n$ be the sample available   to estimate the model parameter.

Below the value of the estimators   $\hat{\theta}_n$ and $\tilde{\theta}_n$ are obtained using AdaGrad  \citep{Duchi}, an adaptive stochastic gradient algorithm, and the strategy exposed in Appendix \ref{subsection:proof:gradient2} for computing the  former estimator is implemented. As  suggested in Section \ref{subsec:algo}, the  algorithms  that   compute    $\hat{\theta}_n$ and $\tilde{\theta}_n$ use the maximum likelihood estimate of the model parameter  as starting value.

\begin{table} 
\begin{center}
\begin{scriptsize}
  \begin{tabular}{|c|c|c|c|c|c|c|c|c|}
   \hline
      
   $\tau$ &type &$n$&  $\beta_{\mathrm{ols},n}$& $\beta_{\mathrm{lad},n}$&$\beta_{\mathrm{rob},n}$& $\hat{\beta}_{\mathrm{mom},n}$ &  $\hat{\beta}_n$
          &  $\tilde{\beta}_n$
   \\

   \hline \hline
   \multirow{3}{*}{0} & & $100$   &0.372&0.353&0.350& 0.409 & 0.355& 0.334\\
    &&1\,000 & 0.116& 0.092&0.104&  0.168 & 0.108&0.107  \\
    &&5\,000  &0.053&0.039&0.046& 0.111 & 0.049&0.047  \\
   \hline
    \multirow{3}{*}{1}  & \multirow{3}{*}{$\mathsf{Y}$} &100 &0.464& 0.339&0.385& 0.707 &0.350&0.342\\
    &&1\,000 &0.181& 0.094&  0.106& 0.921 & 0.105& 0.097  \\
    &&5\,000  & 0.103& 0.043&  0.049& 0.807 & 0.054&0.051 \\
   \hline
     \multirow{3}{*}{2}  & \multirow{3}{*}{$\mathsf{Y}$} &100&0.647&0.351&0.359& 1.315 & 0.337& 0.333 \\
    &&1\,000 &0.241&0.097&0.110& 1.513 &  0.114 &0.115  \\
    &&5\,000  & 0.175&0.039&  0.047& 1.357 & 0.051&0.052 \\
   \hline
       \multirow{3}{*}{3}  & \multirow{3}{*}{$\mathsf{Y}$} &100&0.724&0.331& 0.343& 1.519 & 0.329&0.320 \\
    &&1\,000 & 0.309& 0.100& 0.108& 1.864 & 0.113&0.110 \\
    &&5\,000  & 0.250&0.043&  0.048& 1.759 & 0.053& 0.055  \\
    \hline

          \multirow{3}{*}{1}  & \multirow{3}{*}{$\mathsf{X}$} &100&0.870&0.356&0.374& 0.476 & 0.342&0.338   \\
    &&1\,000 &0.836& 0.111& 0.105& 0.254 & 0.104&0.096    \\
    &&5\,000  & 0.818& 0.065&0.049& 0.174 & 0.054&   0.052    \\
    \hline
          \multirow{3}{*}{2}  & \multirow{3}{*}{$\mathsf{X}$} &100& 1.575&0.400&0.347& 0.655 & 0.337&0.331 \\
    &&1\,000 &1.467&0.160&0.110& 0.319 &  0.112& 0.115  \\
    &&5\,000  & 1.401& 0.119& 0.046& 0.245 & 0.051& 0.052         \\
     \hline
          \multirow{3}{*}{3}  & \multirow{3}{*}{$\mathsf{X}$} &100&1.838&0.442&0.344& 0.743 & 0.331& 0.323  \\
    &&1\,000 & 1.805&  0.216&0.108& 0.377 & 0.113& 0.109\\
    &&5\,000  & 1.771& 0.183&  0.048& 0.293 & 0.054&0.056  \\
    \hline
  \end{tabular}
  \hspace{0.1cm}
      
  \vspace{0.3cm}

 \caption{Results for the Gaussian linear regression model.  For each experimental setting  we report the RMSE over
25 replications.}
 \label{table1}
 \end{scriptsize}
 \end{center}
\end{table}

\subsection{Gaussian linear regression}\label{subsec:expe:linear}

We let $d=8$ and, for every $x\in\R^d$, we let $P_{g(\theta,x)}=\mathcal{N}_1(\beta^\top x,\sigma^2)$ with $\theta=(\beta,\sigma)\in\Theta:=\R^d\times(0,\infty)$. For this example the dataset $d_N^0$ is constructed by simulating  independent observations  using 
$$
Y^0_i=\beta_0^\top X^0_i+\epsilon_i,\quad (X^0_i,\epsilon_i)\iid\mathcal{N}_d(0,I_d)\otimes \mathrm{Laplace}(0,\sigma_0)
$$ 
with $\beta_0=(4,4,3,3,2,2,1,1)$ and $\sigma_0=1$. The model $\{ (P_{g(\theta,x)})_{x\in\R^d},\,\theta\in\Theta\}$ is therefore misspecified and, in what follows, we focus on the estimation of $\beta_0$. We let $N=5\,000$ and, for $n\leq N$ and $\epsilon>0$, the contaminated dataset $d_n$ is generated as explained in Section \ref{sub:sim_setup} where two types of outliers $(x_i^c,y^c_i)$  are considered. More precisely, we say that the outliers are of type   $\mathsf{X}$ when $y_i^c=y_i$ and $x_i^c$ is such that $x^c_{ij}=x_{ij}$ for all $j>1$ and such that $x_{i1}^c$ is a random draw from the from the $\mathcal{N}_1(5,1)$ distribution, and of type  $\mathsf{Y}$ when $y_i^c$  is a random draw from the $\mathcal{N}_1(10,1)$ distribution and $x^c_i=x_i$.

In Table \ref{table1} we report, for different values of $n\leq N$ and of $\epsilon\in[0,0.03]$, the root mean squared error (RMSE)  for the estimation of $\beta_0$ obtained for six estimators, namely the ordinary least squares estimator $\beta_{\mathrm{ols},n}$, the least absolute deviations estimator $\beta_{\mathrm{lad},n}$, the robust estimator of $\beta$ in linear Gaussian regression models proposed by \cite{koller2011sharpening},  computed using the \texttt{R} package \texttt{robustbase},  the robust Median-of-Means estimator $\hat{\beta}_{\mathrm{mom},n}$ introduced by \citet{lecueML},  computed using the \texttt{Python} package \texttt{scikit-learn-extra}\footnote{\url{https://scikit-learn-extra.readthedocs.io/en/stable/index.html}} and  using   three blocks, and the two proposed estimators $\hat{\beta}_n$ and $\tilde{\beta}_n$. 

When the sample is not contaminated, i.e.~when $d_n=d_n^0$, we observe that the ordinary least squares estimator $\beta_{\mathrm{ols},n}$ is the best estimator. However, this estimator is extremely sensitive to the presence of outliers  of both types, a fact that is already well documented in the literature   \citep[see for example][Chapter 1]{Rousseeuw}.  We also observe    that the  estimator $\hat{\beta}_{\mathrm{mom},n}$ performs poorly, with  an RMSE which is in all cases larger than that of $\beta_{\mathrm{ols},n}$. The theory predicts that increasing the number of blocks used by the Median-of-Means procedure  should make $\hat{\beta}_{\mathrm{mom},n}$   more robust  but then the  optimization procedure, as implemented in \texttt{scikit-learn-extra}, becomes less stable. For this reason,  unreported numerical results  have shown that increasing the number of blocks does not improve the empirical performance of $\hat{\beta}_{\mathrm{mom},n}$. As predicted by our theory, the  MMD based estimators $\hat{\beta}_n$ and $\tilde{\beta}_n$ are robust to the two considered  types of outliers. In addition, we observe that their performance is almost identical to that of the robust estimator $\beta_{\mathrm{rob},n}$ proposed \cite{koller2011sharpening}, and that their  RMSE  is in all cases smaller than that of $\beta_{\mathrm{lad},n}$ when the outliers are of type $\mathsf{X}$. However, for outliers of  type $\mathsf{Y}$ this latter estimator tends to be the best one, with an RMSE slightly lower than that of  the estimators   $\hat{\beta}_n$, $\tilde{\beta}_n$  and $\beta_{\mathrm{rob},n}$.

\begin{table}
\begin{center}
\begin{scriptsize}
  \begin{tabular}{|c|c|c|c|c|c|}
   \hline
      
   $\epsilon$ &$n$&  $\theta_{\mathrm{mle},n}$ & $\theta_{\mathrm{rob},n}$&  $\hat{\theta}_n$
          &  $\tilde{\theta}_n$
   \\

   \hline \hline
   \multirow{3}{*}{0\%} & $100$  &1.504&1.737 &2.132 &2.012  \\
   &  $ 1\,000$  &0.565&0.608 & 0.733&0.735  \\
   &  $5\,000$  &0.210& 0.296& 0.290&0.289  \\
    \hline
   \multirow{3}{*}{1\%} & $100$  & 1.767&2.020 &2.017&1.802  \\
   &  $ 1\,000$  &1.325&1.164 & 0.706&0.695  \\
   &  $5\,000$  & 1.293&1.083 &0.283 & 0.283  \\
    \hline
     \multirow{3}{*}{2\%} & $100$  &2.218&4.131 &2.238 &2.186   \\
   &  $ 1\,000$  &1.766&1.551 &0.683 &0.681      \\
   &  $5\,000$  & 1.936& 1.669&0.252 &0.257 \\
    \hline
     \multirow{3}{*}{3\%} & $100$  & 2.657&3.048 &1.786 &1.755   \\
   &  $ 1\,000$  &2.496&2.253 & 0.631& 0.642 \\
   &  $5\,000$  &2.404&2.142 & 0.243& 0.242   \\
    \hline
    
  \end{tabular}
  \hspace{0.5cm}\begin{tabular}{|c|c|c|c|c|c|}
   \hline
      
   $\epsilon$ &$n$& $\beta_{\mathrm{mle},n}$ & $\beta_{\mathrm{rob},n}$&  $\hat{\beta}_n$
          &  $\tilde{\beta}_n$
   \\

   \hline \hline
   \multirow{3}{*}{0\%} & $100$  &1.451&1.696&2.0638&1.944  \\
   &  $ 1\,000$  &0.536&0.574&0.705&0.708  \\
   &  $5\,000$  &0.202&0.238&0.273&0.269   \\
    \hline
   \multirow{3}{*}{1\%} & $100$  &1.598&1.731&1.947&1.727  \\
   &  $ 1\,000$  &1.009&0.584&0.675&0.669  \\
   &  $5\,000$  &0.878&0.273 &0.267&0.266  \\
    \hline
     \multirow{3}{*}{2\%} & $100$  &1.958&3.684 &2.189&2.147 \\
   &  $ 1\,000$  &1.312&0.638&0.654&0.651 \\
   &  $5\,000$  &1.389&0.341&0.239&0.241 \\
    \hline
     \multirow{3}{*}{3\%} & $100$  &2.200&1.989&1.732&1.697   \\
   &  $ 1\,000$  &1.861&0.707&0.604&0.611 \\
   &  $5\,000$  &1.762&0.418&0.231&0.230   \\
    \hline
  \end{tabular}

  \vspace{0.3cm}

 \end{scriptsize}
   \caption{Results for the Heckman sample selection model (synthetic data). The left table is for the estimation of $\theta=(\beta,\gamma,\sigma,\rho)$  while right table is for he estimation of $\beta$ only. For each experimental setting, we report the RMSE  over 25 replications.}
 \label{table2}
 \end{center}

\end{table}

\begin{table}
\begin{center}
\begin{scriptsize}
  \begin{tabular}{|c|c|c|c|c||c|c|c|c||c|c|c|c|}
  \hline
  &\multicolumn{4}{c}{$\epsilon=0\%$}& \multicolumn{4}{c}{$\epsilon=1\%$}&\multicolumn{4}{c}{$\epsilon=3\%$}\\
  \hline
 & $\theta_{\mathrm{mle},n}$&$\theta_{\mathrm{rob},n}$& $\hat{\theta}_n$ &$\tilde{\theta}_n$& $\theta_{\mathrm{mle},n}$&$\theta_{\mathrm{rob},n}$& $\hat{\theta}_n$ &$\tilde{\theta}_n$& $\theta_{\mathrm{mle},n}$&$\theta_{\mathrm{rob},n}$& $\hat{\theta}_n$ &$\tilde{\theta}_n$\\
  \hline
  $\beta_1$&5.037 & 5.401&  5.056      & 5.111&5.445 & 5.441  &5.202 & 5.090& 5.726 &-0.735 & 5.289    & 5.085   \\
   $\beta_2$& 0.212 &  0.201 &0.176    &  0.212&0.180 & 0.190  &0.168 &0.213&0.172  & 0.120 &0.163   & 0.214  \\
   $\beta_3$&0.350 & 0.255 &0.389   &0.318   &0.250 & 0.244  & 0.378 & 0.320& 0.193 & 0.665 &0.369   & 0.319  \\
   $\beta_4$&0.019 & 0.013  & 0.044     &  0.027 &0.013  & 0.015  & 0.038      & 0.026 & 0.006 & 0.064  &0.034    & 0.027  \\
   $\beta_5$&-0.220 &-0.155 &  -0.120     &-0.210 &-0.129 & -0.131  & -0.115   &-0.200 & -0.096 & -0.388 &-0.128  &-0.212 \\
   $\beta_6$& 0.541 & 0.481 &  0.509     &  0.526 &0.489  & 0.474   &
 0.507    & 0.517& 0.475  & 0.753  & 
0.502   & 0.519 \\
   $\beta_7$& -0.029 & -0.067&  -0.032   &  -0.084&-0.056 & -0.064 &-0.048   & -0.075 &-0.090  & 0.177&-0.060    & -0.078    \\
    $\gamma_1$&-0.724&  -0.749 &-1.110  &  -0.827 &-0.613& -0.735 &-1.021  &-0.825&-0.508  & 5.441 & -1.069  & -0.819    \\
    $\gamma_2$& 0.098 &  0.105& 0.159     & 0.116 &0.109 &  0.120 &
 0.150    & 0.125 &0.106  & 0.190 &0.153  & 0.120  \\
    $\gamma_3$& 0.644 & 0.687 & 0.688      & 0.736 &0.584 &  0.665 & 0.713    & 0.752 & 0.505 &  0.244 & 0.715   & 0.742  \\
    $\gamma_4$& 0.070  & 0.070 &0.072       &  0.071&0.060 &  0.064 & 0.069      &  0.069& 0.053 & 0.015  & 0.072  & 0.068 \\
    $\gamma_5$& -0.373 &-0.398&  -0.430     &-0.415 &-0.344& -0.388 &-0.444    & -0.417 & -0.302 & -0.131 &-0.434  & -0.412  \\
    $\gamma_6$& 0.795 & 0.834 & 1.229    & 0.906 &0.559 &  0.753 & 1.219   &  0.888 & 0.344 & 0.474  & 1.171 & 0.869 \\
    $\gamma_7$& 0.182 &  0.183 &  0.205     &  0.212 &0.179 & 0.177 & 0.189   & 0.214& 0.144 & -0.064 &0.198   & 0.207 \\
    $\sigma$& 1.271 & 1.318 & 1.195    & 1.217 &1.351 &  1.357 & 1.205  & 1.229  & 1.454 & 1.357 & 1.220& 1.242 \\
    $\rho$&  -0.124 & -0.051& -0.311    &-0.315  &-0.520& -0.047 &-0.371   & -0.307& -0.720 & -0.047 &-0.381   & -0.308 
\\
\hline
  \end{tabular}

  \vspace{0.3cm}

 \caption{Estimated value of $\theta$ in the Heckman sample selection model for the  2001 Medical Expenditure Panel Survey dataset. \label{table22}}

 \end{scriptsize}
 \end{center}
\end{table}

\subsection{Heckman sample selection model\label{sub:Heck}}

For   $d\in\mathbb{N}$ and $\Theta\subseteq\R^{d}\times(0,\infty)\times(-1,1)$, the Heckman sample selection model   $\big\{ \big(P_{g(\theta,x)}\big)_{x\in\R^d},\,\theta\in\Theta\big\}$ is such that, for all $x\in\R^d$ and $\theta=(\beta,\gamma,\sigma,\rho)\in\Theta$, the distribution $P_{g(\theta,x)}$ is the distribution $P_{Y|x}$  defined in the last example of Proposition \ref{prop:CME:example}. For this model,  in addition to the two MMD based estimators $\hat{\theta}_n$ and $\tilde{\theta}_n$, results are presented for $\theta_{\mathrm{mle},n}$, the maximum likelihood estimator   of $\theta$, computed using the \texttt{R} package \texttt{sampleSelection} \citep{Toomet2008}, and for   $\theta_{\mathrm{rob},n}$, the robust two-step estimator of $\theta$ proposed by \cite{Zhelonkin2016}, computed using the \texttt{R} package \texttt{ssmrob} written by  \citet{zhelonkin2021ssmrob}.  We stress that the  estimator  $\theta_{\mathrm{rob},n}$ is designed specifically for robust estimation in Heckman sample selection models.

We first  let   $d=8$, 
$$
\Theta=\Big\{(\beta,\gamma,\sigma,\rho)\in \R^{d}\times(0,\infty)\times(-1,1) \text{ such that }\beta_{d-i+1}=\gamma_i=0,\,\forall i\in\{1,\dots,d/2\}\Big\}
$$
 and construct the dataset $d_N^0$ by simulating $N=5\,000$ independent   observations  using 
$$
Y^0_i|X^0_i\sim P_{g(\theta_0,X^0_i)},\quad X^0_i \iid \mathcal{N}_{d}(0,I_{d})
$$
 where $\theta_0=(\beta_0 ,\gamma_0 ,\sigma_0,\rho_0)$ with
 $$
 \beta_0=(4,3,2,1,0,0,0,0),\quad \gamma_0=(0,0,0,0,4,3,2,1),\quad \sigma_0=1.5,\quad \rho_0=0.5. 
 $$
The model we consider here assumes that the outcome equation depends only on the first $d/2$ components of $X_i^0$ while the selection equation depends only on the last $d/2$ components of $X_i^0$, so that the total number of parameter to estimate is $d+2=10$. For $n\leq N$ and $\epsilon>0$ the contaminated dataset $d_n$ is constructed as explained in Section \ref{sub:sim_setup}, with $y_i^c=y_i$ and $x_i^c$ such that $x^c_{ij}=x_{ij}$ for all $j>1$ and such that $x_{i1}^c$ is a random draw from the  $\mathcal{N}_1(5,1)$ distribution.

Table \ref{table2} shows, for different values of $n\leq N$ and of $\epsilon\in[0,0.03]$, the  RMSE  obtained with the different estimators for the estimation of $\theta_0$ as well as for the estimation of $\beta_0$, which is often the main parameter of interest in this model. We observe that the maximum likelihood estimator is the best estimator when there are no outliers, as expected from the asymptotic theory. On the other hand, this estimator is sensitive to the presence of outliers. The robust estimator $ \theta_{\mathrm{rob},n}$ of \cite{Zhelonkin2016} improves upon $ \theta_{\mathrm{mle},n}$ when the sample is contaminated. However,  for all the considered values of $\epsilon>0$ and  $n$, this estimator is dominated by $\hat{\theta}_n$ and by $\tilde{\theta}_n$. When restricting our attention to the estimation of $\beta$, we observe in  Table \ref{table2} that the two MMD based estimators have a lower RMSE than $\beta_{\mathrm{rob},n}$ for all $\epsilon>0$ when $n=5\,000$ and for $\epsilon=0.03$ when $n=1\,000$.

We now consider the   real dataset  and  model used in \citet[Table 4]{Zhelonkin2016}. The dataset, available from the \texttt{R} package \texttt{ssmrob}, contains $N=3\,328$ observations $d_N^0:=\{(y_{1i}^0,y_{2i}^0, x_i^0)\}_{i=1}^N$, extracted from the 2001 Medical Expenditure Panel Survey, and an Heckman sample selection model is used to regress $y^0_{i1}$, the  log ambulatory expenses for the $i$th individual,  on the vector $x_i^0$ containing $d=7$ covariates, including an intercept. Both the selection and the outcome equation is assumed to depend on all the components of $x_i^0$, so that $\Theta=\R^{d}\times(0,\infty)\times(-1,1)$ and the number of parameter to estimate is   $2d+2=16$. We let $n=N$ and, for a given $\epsilon>0$, we contaminate  the dataset $d_N^0$, where    $y_{2i}^0=\ind_{(0,\infty)}(y_{1i}^0)$, by applying the approach described in Section \ref{sub:sim_setup}  with $(y_{1i}^c,y_{2i}^c, x_i^c)=(y_{1i}^0,1-y_{2i}^0, x_i^0)$. Below  we denote by $d_{n,\epsilon}$ the resulting contaminated version of $d_n^0$.

The estimated parameters values, obtained for  $\epsilon\in\{0,0.01,0.03\}$ and the four considered estimators,  are presented in Table \ref{table22}. For $\epsilon=0$  the results obtained with $\theta_{\mathrm{mle},n}$ and with $\theta_{\mathrm{rob},n}$ reproduce those given in  \cite{Zhelonkin2016}, and the two MMD based estimators provide   similar estimated values of the model parameters. In order to assess the sensibility of the different estimators to a contamination of the data, in Table \ref{table:SS_real} we give the value of $\Delta_{\epsilon,n}(\theta_n):=\|\theta_n(d_{n,\epsilon})-\theta_n(d^0_{n})\|$ for  all $\theta_n\in\{\theta_{\mathrm{mle},n}, \theta_{\mathrm{rob},n},\hat{\theta}_n,\tilde{\theta}_n\}$. As expected, the maximum likelihood estimator is very sensitive to the presence of outliers and $\theta_{\mathrm{rob},n}$ improves upon $\theta_{\mathrm{mle},n}$. The most striking feature of Table \ref{table:SS_real} is the remarkable performance of $\tilde{\theta}_n$. Notably, the   value  $\Delta_{\epsilon,n}(\tilde{\theta}_n)$ is about 2.675 times smaller than that of $\Delta_{\epsilon,n}(\theta_{\mathrm{rob},n})$ when $\epsilon=0.01$, and about 8.5
times smaller when $\epsilon=0.03$. For this latter value of $\epsilon$  the estimator $\hat{\theta}_n$ outperforms the estimator $\theta_{\mathrm{rob},n}$  since,  in this case, $\Delta_{\epsilon,n}(\hat{\theta}_n)$ is  about 1.8   times smaller than   $\Delta_{\epsilon,n}(\theta_{\mathrm{rob},n})$. On the contrary, when there is a small proportion of outliers, i.e.~when $\epsilon=0.01$, we observe from Table \ref{table:SS_real} that  $\theta_{\mathrm{rob},n}$ dominates  $\hat{\theta}_n$. We however recall that $\hat{\theta}_n$ depends on a   kernel $k_{\mathcal{X}}$  which is kept the same throughout this Section \ref{section:simulation}, and thus which is not optimized in any sense  to the particular problem at hand.

\begin{table}
\begin{center}
\begin{tabular}{c||c|c|c|c}
$\epsilon$&$\theta_{\mathrm{mle},n}$&$\theta_{\mathrm{rob},n}$&$\hat{\theta}_{n}$&$\tilde{\theta}_{n}$\\
\hline
1\%&0.651&0.107&0.186&0.040\\
3\%& 1.090&0.468&0.259&0.055
\end{tabular}
\caption{Value of $\Delta_{\tau,n}(\theta_n)$ for the Heckman sample selection model for the  2001 Medical Expenditure Panel Survey dataset, with $\theta_n\in\{\theta_{\mathrm{mle},n}, \theta_{\mathrm{rob},n},\hat{\theta}_{n},\tilde{\theta}_{n}\}$.\label{table:SS_real}}
\end{center}
\end{table}

\subsection{Gamma regression model\label{sub:gamma}}

With  $d=8$ and  $\Theta=\R^d\times (0,\infty)$, we  consider the  Gamma regression model $\big\{ \big(P_{g(\theta,x)}\big)_{x\in\R^d},\,\theta\in\Theta\big\}$  which is such that, for all $x\in\R^d$ and $\theta=(\beta,\nu)\in\Theta$, the distribution $P_{g(\theta,x)}$ is the distribution $P_{Y|x}$  defined in the fourth example of Proposition \ref{prop:CME:example}. For this example  the dataset $d_N^0$ is obtained by simulating $N=5\,000$  independent   observations using $Y^0_i|X^0_i\sim P_{g(\theta_0,X^0_i)}$ and $X^0_i\iid \mathcal{N}_d(0,I_d)$, with $\theta_0=(1,\dots,1)$. Then, for every $n\leq N$ and $\epsilon\in [0,0.03]$, the contaminated dataset $d_n$ is constructed as described in Section \ref{sub:sim_setup},   with $y_i^c=y_i$ and $x_i^c$ such that $x^c_{ij}=x_{ij}$ for all $j>1$ and such that $x_{i1}^c$ is a random draw from the from $\mathcal{N}_1(-0.5,1)$ distribution.

Table \ref{table3} presents the RMSE for the estimation of $\theta_0$ and $\beta_0$ obtained with     $\hat{\theta}_n$ and $\tilde{\theta}_n$. Results are also given for  $\theta_{\mathrm{mle},n}$, the maximum likelihood estimator of $\theta$, and for  $\theta_{\mathrm{rob},n}$, the estimator proposed by \cite{cantoni2001,cantoni2006} for robust inference in generalized linear model, computed using the \texttt{R} package \texttt{robustbase}. We observe from this table that, as expected, $\theta_{\mathrm{mle},n}$ is not robust to the presence of outliers, while the two proposed MMD based estimators are. In all cases, the RMSE  obtained with the estimator $\theta_{\mathrm{rob},n}$ proposed by  \cite{cantoni2001,cantoni2006} is however slightly smaller than the that obtained with $\hat{\theta}_n$ and $\tilde{\theta}_n$.  

We stress that, in some sense, the better performance of $\theta_{\mathrm{rob},n}$ is reassuring  since this estimator is precisely  designed  for robust inference in generalized linear models and, in particular, has been motivated in  \cite{cantoni2006} for robust inference in Gamma regression models. By contrast, the applicability of the estimators $\hat{\theta}_n$ and $\tilde{\theta}_n$, and their theoretical guarantees   derived in Section \ref{section:theory} regarding their robustness, hold  for a much broader class of regression models.

\begin{table}
\begin{center}
\begin{scriptsize}
  \begin{tabular}{|c|c|c|c|c|c|}
   \hline

   $\epsilon$ &$n$&  $\theta_{\mathrm{mle}}$ & $\theta_{\mathrm{rob}}$&  $\hat{\theta}_n$
          &  $\tilde{\theta}_n$
   \\
 
   \hline \hline
   \multirow{3}{*}{0\%} & $100$  & 0.478 &0.425&0.501 &0.498    \\
   &  $ 1\,000$  &0.123&0.120 &0.142&0.144 \\
   &  $5\,000$  &0.055&0.053&0.067&0.069     \\
    \hline
   \multirow{3}{*}{1\%} & $100$  &0.394&0.398 &0.521&0.521   \\
   &  $ 1\,000$  &0.250& 0.108 &0.149 &0.148\\
   &  $5\,000$  &0.350&0.049 &0.071 &0.070  \\
    \hline
     \multirow{3}{*}{2\%} & $100$  &0.395&0.381&0.539  &0.545  \\
   &  $ 1\,000$  &0.353 & 0.113& 0.146& 0.145 \\
   &  $5\,000$  &0.551 &0.056 &0.073&0.073 \\
    \hline
     \multirow{3}{*}{3\%} & $100$  &0.441&0.407&0.471 &0.475     \\
   &  $ 1\,000$  &  0.349 &0.119& 0.157& 0.158 \\
   &  $5\,000$  & 0.631&0.069 & 0.076 & 0.078 \\
    \hline
  \end{tabular} \hspace{0.5cm} \begin{tabular}{|c|c|c|c|c|c|}
   \hline

   $\epsilon$ &$n$&  $\beta_{\mathrm{mle}}$ & $\beta_{\mathrm{rob}}$&  $\hat{\beta}_n$
          &  $\tilde{\beta}_n$
   \\
 
   \hline \hline
   \multirow{3}{*}{0\%} & $100$  &0.308&0.306 &0.444 &0.443     \\
   &  $ 1\,000$  &0.087&0.093&0.132&0.133     \\
   &  $5\,000$  &0.041&0.043&0.061 &0.063       \\
    \hline
   \multirow{3}{*}{1\%} & $100$  &0.297 &0.313&0.473&0.475    \\
   &  $ 1\,000$  &0.114&0.098& 0.132&0.133  \\
   &  $5\,000$  &  0.070&0.042&0.061&  0.059 \\
    \hline
     \multirow{3}{*}{2\%} & $100$  &0.302 &0.301& 0.490&0.493 \\
   &  $ 1\,000$  & 0.120&0.094&0.127& 0.127  \\
   &  $5\,000$  & 0.105& 0.046&0.066& 0.067\\
    \hline
     \multirow{3}{*}{3\%} & $100$  &0.296 & 0.296&0.403& 0.406     \\
   &  $ 1\,000$  &0.119&  0.091&0.138& 0.136   \\
   &  $5\,000$  & 0.129& 0.044& 0.067& 0.069  \\
    \hline
    \end{tabular}
  \vspace{0.3cm} 
 \caption{Results for the Gamma regression model. For each experimental setting, we report the mean square error over
25 replications.}
 \label{table3}
 \end{scriptsize}
 \end{center}
\end{table}

\section{Conclusion}
\label{section:conclusion}

Some important questions remain open, such as the dependence of the convergence rate of the estimator $\hat{\theta}_n$ to the dimension of the parameter space $\Theta$, and the existence of  conditional mean embedding operators  when  $\mathcal{X}$ is unbounded  which would enable to apply some our theoretical results for this estimator   to problems where the regressors can take  arbitrarily large values. Finally,  further work is needed establish  non-asymptotic guarantees for the estimator $\tilde{\theta}_n$.

\section*{Acknowledgements}
We would like to acknowledge the Heilbronn Institute for Mathematical Research who funded a two-week  visit for P. A.~to the University of Bristol, during which this work was started. We thank Timoth\'ee Mathieu (INRIA) who provided detailed explanations on the MOM implementation in \texttt{scikit-learn-extra}.

\section*{Supplementary material}

Supplementary material available at \textit{Biometrika} includes additional information about the computation of the two estimators and all the proofs.

\bibliographystyle{apalike}
\bibliography{complete}

\appendix

\section{Computation of the estimators}
\label{subsection:proof:gradient}

\subsection{Gradient of the loss\label{subsection:proof:gradient1}}

\begin{proposition}\label{prop:gradient:1}
Assume that each $P_\lambda$ has a density $p_\lambda$ with respect to a measure $\mu$ such that $\lambda\mapsto p_\lambda $ is differentiable, and that  $\theta\mapsto g(\theta,x)$ is differentiable for any $x\in\setX$.
\begin{enumerate}
\item Assume that there is a  function   $\hat{b}:\setY^2\rightarrow\R$ such that  $$
 \int_\setY\int_\setY \hat{b}(y,y') \mu(\mathrm{ d}y)\mu(\mathrm{ d}y') < \infty 
 $$
  and such that,   for all $(\theta,x,x',y,y')$,
  $$
 \big|k((x,y),(x',y')) \nabla_{\theta}p_{g(\theta,x)}(y)p_{g(\theta',x')}(y')  \big| \leq \hat{b}(y,y').
 $$
Then, for all $(\theta,x,x',y)$ we have
 \begin{align*}
&\nabla_\theta   \hat{\ell} (\theta, x,x',y)\\
 &=2 \E_{Y\sim P_{g(\theta,x)}, \,\,Y'\sim P_{g(\theta,x') }}   \bigg[\Big(k\big((x,Y),(x',Y')\big) -  k\big((x,Y),(x',y)\big)\Big) \nabla_{\theta}  \log p_{g(\theta,x)}(Y)  \bigg].
 \end{align*}
\item Assume that there exists  a function $\tilde{b}:\setY^2\rightarrow\R$ such that 
$$
\int_{\setY}\int_{\setY} \tilde{b}(y,y') \mu(\mathrm{ d}y)\mu(\mathrm{ d}y') <  \infty
$$
 and such that,  for all $(\theta,x,y,y')$,
$$
\big|k(y,y') \nabla_{\theta}[p_{g(\theta,x)}(y)p_{g(\theta',x)}(y') ] \big| \leq \tilde{b}(y,y').
$$
 Then, for all $(\theta,x,y)$ we have
$$
\nabla_\theta \tilde{\ell}(\theta, x,y)=2\E_{\substack{ Y,Y'\,\,\iid P_{g(\theta,x)}}}\Big[  \Big(k_{\mathcal{Y}}(Y,Y') - k_{\mathcal{Y}}(Y,y) \Big) \nabla_{\theta}  \log p_{g(\theta,x)}(Y)   \Big].
$$
\end{enumerate} 
\end{proposition}

\begin{remark}
We need more assumption to ensure stability and convergence of the stochastic gradient algorithm. See for example Proposition 5.2 in~\cite{cherief2019finite} (and the references therein), where the authors require the existence of the variance of 
$$
\hat{L}(\theta,x,x'. U,U',y):=2\Big(k\big((x,Y),(x',Y')\big) -  k\big((x,Y),(x',y)\big)\Big) \nabla_{\theta}  \log p_{g(\theta,x)}(Y)
$$  
when $U\sim P_{g(\theta,x)}$ and $U'\sim P_{g(\theta,x')}$. However, under Assumption~\ref{asm:k}, it boils down to the corresponding assumption on $\nabla_{\theta}  \log p_{g(\theta, x)}(U)$. For example, if there is $v>0$ such that for any $(x,\theta)$, $\mathbb{E}_{U\sim P_{g(\theta,x)} }[ \| \nabla_{\theta}  \log p_{g(\theta, x)}(U) \|^2 ]\leq v$, then  
$$
\mathrm{Var}(\hat{L}(\theta,x,x'. U,U',y)) \leq 16v,\quad \forall  (\theta,x,x',y).
$$
\end{remark}

\begin{proof}  We start by the proof of point 2. By definition,
\begin{align*}
 \tilde{\ell}(\theta,X_i,Y_i)
&=  \mathbb{E}_{Y\sim P_{g(\theta,X_i)},Y'\sim P_{g(\theta,X_i)} } \big[k_{\mathcal{Y}}(Y,Y') -2 k_{\mathcal{Y}}(Y,Y_i)\big]
\\
&=  \iint \big[k_{\mathcal{Y}}(y,y') -2 k_{\mathcal{Y}}(y,Y_i)\big]  p_{g(\theta,X_i)}( y) p_{g(\theta,X_i)}(y') \mu(\mathrm{d}y)\mu(\mathrm{d}y')
\\
&=  \iint k_{\mathcal{Y}}(y,y')  p_{g(\theta,X_i)}( y) p_{g(\theta,X_i)}(y') \mu(\mathrm{d}y)\mu(\mathrm{d}y') - 2  \int k_{\mathcal{Y}}(y,Y_i)  p_{g(\theta,X_i)}( y) \mu(\mathrm{d}y),
\end{align*}
so that
\begin{align}
\nabla_\theta  \tilde{\ell}(\theta,X_i,Y_i) \nonumber
&=  \nabla_\theta \iint k_{\mathcal{Y}}(y,y')  p_{g(\theta,X_i)}( y) p_{g(\theta,X_i)}(y') \mu(\mathrm{d}y)\mu(\mathrm{d}y')
\\
& - \nabla_\theta \int k_{\mathcal{Y}}(y,Y_i)  p_{g(\theta,X_i)}( y) \mu(\mathrm{d}y) \nonumber
\\
&= \iint k_{\mathcal{Y}}(y,y')  \nabla_\theta \left[ p_{g(\theta,X_i)}( y) p_{g(\theta,X_i)}(y')\right] \mu(\mathrm{d}y)\mu(\mathrm{d}y') \nonumber
\\
&  - 2 \int k_{\mathcal{Y}}(y,Y_i) \nabla_\theta \left[ p_{g(\theta,X_i)}( y)\right] \mu(\mathrm{d}y)
\label{proof:prop:1:1}
\end{align}
where the inversion of $\int$ and $\nabla$ is jusfified thanks to the existence of the function $\tilde{b}$. Remark that
$$
\nabla_\theta \left[ p_{g(\theta,X_i)}( y)\right] = \nabla_\theta \left[ \log p_{g(\theta,X_i)}( y)\right] p_{g(\theta,X_i)}
$$
and that
\begin{multline*}
\nabla_\theta \left[ p_{g(\theta,X_i)}( y) p_{g(\theta,X_i)}( y')\right] 
\\
= \nabla_\theta \left[ \log p_{g(\theta,X_i)}( y)\right]  p_{g(\theta,X_i)}( y) p_{g(\theta,X_i)}( y') + \nabla_\theta \left[ \log p_{g(\theta,X_i)}( y')\right]  p_{g(\theta,X_i)}( y) p_{g(\theta,X_i)}( y').
\end{multline*}
Plugging this into~\eqref{proof:prop:1:1} gives:
\begin{align*}
\nabla_\theta  \tilde{\ell}(\theta,X_i,Y_i)
&= \iint k_{\mathcal{Y}}(y,y') \nabla_\theta \left[ \log p_{g(\theta,X_i)}( y)\right]  p_{g(\theta,X_i)}( y) p_{g(\theta,X_i)}( y')   \mu(\mathrm{d}y)\mu(\mathrm{d}y')
\\
&  + \iint k_{\mathcal{Y}}(y,y') \nabla_\theta \left[ \log p_{g(\theta,X_i)}( y')\right]  p_{g(\theta,X_i)}( y) p_{g(\theta,X_i)}( y')   \mu(\mathrm{d}y)\mu(\mathrm{d}y')
\\
&  - 2\int k_{\mathcal{Y}}(y,Y_i) \nabla_\theta \left[ \log p_{g(\theta,X_i)}( y)\right] p_{g(\theta,X_i)} \mu(\mathrm{d}y)
\\
& = 2 \iint k_{\mathcal{Y}}(y,y') \nabla_\theta \left[ \log p_{g(\theta,X_i)}( y)\right]  p_{g(\theta,X_i)}( y) p_{g(\theta,X_i)}( y)   \mu(\mathrm{d}y)\mu(\mathrm{d}y')
\\
&  - 2\sum_{i=1}^n \int k_{\mathcal{Y}}(y,Y_i) \nabla_\theta \left[ \log p_{g(\theta,X_i)}( y)\right] p_{g(\theta,X_i)} \mu(\mathrm{d}y)
\end{align*}
by symmetry, and thus,
$$
\nabla_\theta  \tilde{\ell}(\theta,X_i,Y_i)
=  \frac{2}{n} \sum_{i=1}^n \mathbb{E}_{Y\sim P_{g(\theta,X_i)},Y'\sim P_{g(\theta,X_i)} } \biggl\{ \bigl[k_{\mathcal{Y}}(Y,Y') -k_{\mathcal{Y}}(Y,Y_i)\bigr]
\nabla_\theta \left[ \log p_{g(\theta,X_i)}( Y)\right]
\biggr\}.
$$

The proof of point 1, from the expression in~\eqref{dfn:estim}, is exactly similar. \end{proof}

\subsection{A closer look at the computation of   \texorpdfstring{$\hat{\theta}_n$}{Lg} }\label{subsection:proof:gradient2}

Let $k=k_{\gamma}\otimes k_{\mathcal{Y}}$ with $k_{\gamma}$ as in Section \ref{sub:link} and let $L(\theta, x,x',y)$ be a random variable such that $\E[L(\theta, x,x',y)]=\nabla_\Theta \ell(\theta,x,x',y)$, with  $\ell(\theta,x,x',y)$  as defined in Section \ref{sub:link}. Then, given $n$ observations $d_n:=\{(x_i,y_i)\}_{i=1}^n$ in $\mathcal{Z}$, the random variable
$$
H_n\big(\gamma,\theta,d_n\big):=2\sum_{i=1}^{n-1}\sum_{j=i+1}^n k_{\gamma}(x_i,x_j)L(\theta, x_i,x_j,y_j)
$$
is such that $\E[H_n\big(\gamma,\theta,d_n\big)]=\nabla_\theta h_n(\gamma,\theta, d_n)$, with $h_n(\gamma,\theta, d_n)$ as defined in \eqref{eq:Hn}. 

Next, for an integer $M_1\in\{1,\dots,(n-1)n/2-1\}$  we let   
$$
\mathcal{S}_{M_1}\subset \mathcal{S}:=\{(i,j): 1\leq i<j\leq n\}
$$
 be such that the set $\{k_{\gamma}(x_i,x_j)\}_{(i,j)\in \mathcal{S}_{M_1}}$ contains the $M_1$ largest elements of the set $\{k_{\gamma}(x_i,x_j)\}_{(i,j)\in\mathcal{S}}$, and for an integer  $M_2\in\mathbb{N}$ such that $M_1+M_2\leq (n-1)n/2$ we  let $\{ (I_i,J_i)\}_{i=1}^{M_2}$ be a simple random  sample    obtained without replacement from the set $\mathcal{S}\setminus \mathcal{S}_{M_1}$. Then,  the random variable
\begin{align*}
H^{(M_1,M_2)}_n(\gamma,\theta, d_n)&:=2\sum_{(i,j)\in\mathcal{S}_{M_1}} k_{\gamma}(x_i,x_j)L(\theta, x_i,x_j,y_j)\\
&+\frac{(n-1)n-2M_1}{M_2}\sum_{m=1}^{M_2} k_{\gamma}(x_{I_m},x_{J_m})L(\theta, x_{I_m},x_{J_m},y_{J_m})
\end{align*}
is such that $\E[H^{(M_1,M_2)}_n(\gamma,\theta, d_n)]=h_n(\gamma,\theta, d_n)$, and thus
\begin{align}\label{eq:grad_sto}
\E\bigg[\sum_{i=1}^N L(\theta, x_i,y_i)+H^{(M_1,M_2)}_n(\gamma,\theta, d_n)\bigg]=\nabla_\theta \sum_{ i,j=1}^n \hat{\ell}(\theta,X_i,X_j, Y_j).
\end{align}

This approach for computing an unbiased estimate of $\nabla_\theta \sum_{ i,j=1}^n \hat{\ell}(\theta,X_i,X_j, Y_j)$ involves the    construction of the sets $\mathcal{S}$ and $\mathcal{S}_{M_1}$, which requires $\bigO(n^2)$ operations. However,  once these two sets are obtained,  obtaining a realization of $G_n(\theta,d_n):=\sum_{i=1}^N L(\theta, x_i,y_i)+H^{(M_1,M_2)}_n(\gamma,\theta, d_n)$ for a given $\theta$ can be done in only  $\bigO(n+M_1+M_2\log (M_2))$ operations using e.g.\ the simple random sampling without replacement method proposed by \cite{Gupta}.

For this procedure to work well in practice the parameters $M_1$ and $M_2$ must be such that  the variance of $G_n(\theta,d_n)$  is small.  When a  small value for $\gamma$ is chosen it is often true that $k_{\gamma}(x_i,x_j)\approx 0$ for most pairs $(i,j)\in\mathcal{S}$. When this happens, taking $M_1=\bigO(n)$ and $M_2$ such that $M_2\log(M_2)=\bigO(n)$  allows to efficiently compute $\hat{\theta}_n$ using a stochastic gradient algorithm whose cost per iteration  is linear in the sample size $n$. However, the memory requirement the approach we just described is $\bigO(n^2)$, which limits is applicability to moderate values of $n$ (to $n$ equals to a few thousands, say).

\section{Proof of Lemma \ref{lemma:CBM_New}\label{p-lemma:CBM_New}}

\subsection{Preliminaries}

We first recall the following result \citep[see][Proposition 1.6]{da2014stochastic}:

\begin{lemma}\label{lemma:swap} 
Let $\mathcal{A}$ and $\mathcal{B}$ be two Hilbert spaces,  $T:\mathcal{A}\rightarrow \mathcal{B}$ be a bounded linear operator  and $Z$ be a random variable taking values in    $A$ and such that $\mathbb{E}\big[\|Z\|_A\big]<\infty$. Then $\mathbb{E}[T(Z)]=T(\mathbb{E}[Z])$.
\end{lemma}


We recall that, under Assumption \ref{asm:k}-\ref{asm:kxy}, for any probability distribution $P\in\mathcal{P}(\mathcal{Z})$ the mean embedding  $\mu(P)  = \mathbb{E}_{Z\sim P}[k(Z,\cdot)] $  of $P$  is well defined in $\mathcal{H}$,  and that   $\mu(P)$ has the key property to be such that
\begin{align*}
<f,\mu(P)>_{\mathcal{H}}=<f,\mathbb{E}_{Z\sim P}[k(Z,\cdot)]>_{\mathcal{H}}=\mathbb{E}_{Z\sim P}\big[<f,k(Z,\cdot)>_{\mathcal{H}}\big]=\E_{Z\sim P}[f(Z)],\quad\forall f\in\mathcal{H}
\end{align*}
where the second equality holds by Lemma \ref{lemma:range}, noting that for all $f\in\mathcal{H}$ the mapping $g\mapsto <f,g>_{\mathcal{H}}$ is a bounded linear operator on $\mathcal{H}$ while, under Assumption \ref{asm:kxy}, $\E_{Z\sim P}[\|k(Z,\cdot)\|_{\mathcal{H}}]\leq 1$.

Recall also that the boundedness of $k_{\mathcal{X}}$ and $k_{\mathcal{Y}}$ (Assumption~\ref{asm:kxy}) implies that $\mathcal{C}_{P}$ and $\mathcal{C}_{P_X}$ exist, are unique, and that they are bounded, linear operators \citep[see][Section 3]{fukumizu2004}.

We then have the following result  (also proved in the proof of Corollary 3 in \citealp{fukumizu2004} as well as in \citealp{klebanov2020rigorous}, Theorem 4.1).
\begin{lemma}\label{lemma:range}
Assume that Assumption \ref{asm:k}-\ref{asm:kxy} and condition \eqref{eq:cond_embed} hold. Then, $\mathrm{range}(\mathcal{C}^*_P)\subseteq \mathrm{range}(\mathcal{C}_{P_X})$.
\end{lemma}
\begin{proof}
Let $g\in\mathcal{H}_{\mathcal{Y}}$ and $f\in\mathcal{H}_{\mathcal{X}}$. Then,
\begin{align*}
<\mathcal{C}^*_{P}g, f>_{\mathcal{H}_{\mathcal{X}}}&=<g, \mathcal{C}_{P}f>_{\mathcal{H}_{\mathcal{Y}}}\\
&=\E_{(X,Y)\sim P}\big[g(Y)f(X)\big]\\
&=\E_{X\sim P_X}\big[\E\big[g(Y)| X] f(X)\big]\\
&=<\E_{Y\sim P_{Y| \cdot}}[g(Y)], \mathcal{C}_{P_X} f >_{\mathcal{H}_{\mathcal{X}}}\\
&=<\mathcal{C}_{P_X}\E_{Y\sim P_{Y|\cdot}}[g(Y)],  f >_{\mathcal{H}_{\mathcal{X}}}
\end{align*}
where the fourth equality holds under \ref{eq:cond_embed} and the last one uses the fact that $\mathcal{C}_{P_X}$ is self-adjoint. Since $g\in\mathcal{H}_{\mathcal{Y}}$ and $f\in\mathcal{H}_{\mathcal{X}}$ are arbitrary, it follows that
$$
\mathcal{C}^*_{P}g=\mathcal{C}_{P_X}\E_{Y\sim P_{Y|\cdot}}[g(Y)],\quad\forall g\in\mathcal{H}_{\mathcal{Y}}
$$
and the  proof of the lemma  is complete.

\end{proof}

\subsection{Proof of the lemma}

\begin{proof}
Let $I:\mathcal{H}_{\mathcal{H}}\rightarrow \mathcal{H}_{\mathcal{H}}$ be the identity operator on $\mathcal{H}_{\mathcal{X}}$ and $\mathcal{P}_{\mathrm{Ker}(\mathcal{C}_{P_X})}:\mathcal{H}_{\mathcal{H}}\rightarrow \mathcal{H}_{\mathcal{H}}$ be the orthogonal projection on $\mathrm{Ker}(\mathcal{C}_{P_X})$. Recall that  $\mathcal{P}_{\mathrm{Ker}(\mathcal{C}_{P_X})}$ is a linear operator such that $\|\mathcal{P}_{\mathrm{Ker}(\mathcal{C}_{P_X})}\|_{\mathrm{o}}=1$. Therefore,
the linear operator $\mathcal{C}^\dagger_{P_X}\mathcal{C}_{P_X}=\mathcal{I}-\mathcal{P}_{\mathrm{Ker}(\mathcal{C}_{P_X})}$ is   bounded. In addition, by Lemma \ref{lemma:range}, $\mathrm{range}(\mathcal{C}^*_P)\subseteq \mathrm{range}(\mathcal{C}_{P_X})$ and therefore $\mathcal{C}^\dagger_{P_X}\mathcal{C}^*_{P}:\mathcal{H}_{\mathcal{Y}}\rightarrow \mathcal{H}_{\mathcal{X}}$ is a bounded linear operator \citep[][Theorem 2.3]{arias2009reduced}. Hence, recalling that if $A:\mathcal{H}_1\rightarrow\mathcal{H}_2$ is a bounded linear operator between two Hilbert spaces then $\|A^*\|_{\mathcal{H}_2}=\|A\|_{\mathcal{H}_1}$, it follows that $(\mathcal{C}^\dagger_{P_X}\mathcal{C}^*_{P})^*:\mathcal{H}_{\mathcal{X}}\rightarrow \mathcal{H}_{\mathcal{Y}}$   is a bounded linear operator.

To proceed further let    $g\in\mathcal{H}_{\mathcal{Y}}$ and $f\in\mathcal{H}_{\mathcal{X}}$. Then,
\begin{equation}\label{eq:C1}
\begin{split}
<f, \mathcal{C}_{P_X}\E_{Y\sim P_{Y|\cdot}}[g(Y)]>_{\mathcal{H}_{\mathcal{X}}}&=\E_{X\sim P_X}\big[f(X)  \E[g(Y)|X]\big]\\
&=\E_{(X,Y)\sim P_X}\big[g(Y) f(X)\big]\\
&=<g, \mathcal{C}_P f>_{\mathcal{H}_{\mathcal{Y}}}\\
&=<\mathcal{C}_P^* g,  f>_{\mathcal{H}_{\mathcal{X}}}
\end{split}
\end{equation}
while, on the other hand, recalling that $f'=\mathcal{C}_{P_X}\mathcal{C}^\dagger_{P_X} f'$ for all $f'\in\mathrm{range}(\mathcal{C}_{P_X})$, and recalling that $\mathrm{range}(\mathcal{C}^*_P)\subseteq \mathrm{range}(\mathcal{C}_{P_X})$ by Lemma \ref{lemma:range},
\begin{align}\label{eq:C2}
<f, \mathcal{C}_{P_X}\mathcal{C}^\dagger_{P_X}\mathcal{C}^*_{P}g>_{\mathcal{H}_{\mathcal{X}}}&=<f,  \mathcal{C}^*_{P}g>_{\mathcal{H}_{\mathcal{X}}}.
\end{align}
Hence, by \eqref{eq:C1}-\eqref{eq:C2}, it follows that
$$
<f, \mathcal{C}_{P_X}\big(\E_{Y\sim P_{Y|\cdot}}[g(Y)]-  \mathcal{C}^\dagger_{P_X}\mathcal{C}^*_{P}g\big)>_{\mathcal{H}_{\mathcal{X}}}=0 
$$
and thus
\begin{align*}
\E_{X\sim P_X}\Big[f(X)\big(\E_{Y\sim P_{Y|\cdot}}[g(Y)]&-  \mathcal{C}^\dagger_{P_X}\mathcal{C}^*_{P}g\big)(X)\Big]\\
 &=<f, \mathcal{C}_{P_X}\big(\E_{Y\sim P_{Y|\cdot}}[g(Y)]-  \mathcal{C}^\dagger_{P_X}\mathcal{C}^*_{P}g\big)>_{\mathcal{H}_{\mathcal{X}}}\\
&=0. 
\end{align*}
Consequently, since $f\in\mathcal{H}_{\mathcal{X}}$ is arbitrary, it follows that, under the assumptions of the lemma,  
\begin{align}\label{eq:span}
\E_{Y\sim P_{Y|\cdot}}[g(Y)]=  \mathcal{C}^\dagger_{P_X}\mathcal{C}^*_{P}g 
\end{align}

Remark now that for $x\in\mathcal{X}$ and $y\in\mathcal{Y}$ we have
\begin{equation}\label{eq:xy}
\begin{split}
 \mathcal{C}^\dagger_{P_X}\mathcal{C}^*_{P}k_{\mathcal{Y}}(y,\cdot)(x)&=<\mathcal{C}^\dagger_{P_X}\mathcal{C}^*_{P}k_{\mathcal{Y}}(y,\cdot), k_{\mathcal{X}}(x,\cdot)>_{\mathcal{H}_{\mathcal{X}}}\\
 &=<k_{\mathcal{Y}}(y,\cdot), (\mathcal{C}^\dagger_{P_X}\mathcal{C}^*_{P})^*k_{\mathcal{X}}(X,\cdot)>_{\mathcal{H}_{\mathcal{Y}}}\\
 &=\big(\mathcal{C}^\dagger_{P_X}\mathcal{C}^*_{P}\big)^*k_{\mathcal{X}}(x,\cdot)(y)
\end{split}
\end{equation}
where the first equality uses the reproducing property of $k_{\mathcal{X}}$ and the third equality the  reproducing property of $k_{\mathcal{Y}}$.

Let $y\in\mathcal{Y}$ and $x\in\mathcal{X}$. Then, using \eqref{eq:span} with $g=k_{\mathcal{Y}}(y,\cdot)$ and \eqref{eq:xy}, we have
\begin{align*}
\mu(P_{Y|x})(y)&=\E_{Y\sim P_{Y|x}}\big[k_{\mathcal{Y}}(y,Y)\big]\\
&= \mathcal{C}^\dagger_{P_X}\mathcal{C}^*_{P}k_{\mathcal{Y}}(y,\cdot)(x) \\
&=\big(\mathcal{C}^\dagger_{P_X}\mathcal{C}^*_{P}\big)^*k_{\mathcal{X}}(x,\cdot)(y)
\end{align*}
and the proof is complete.
\end{proof}

\section{Proof of Theorem \ref{thm:Kernel_Gen2}\label{p-thm:Kernel_Gen2}}

\subsection{A preliminary result for proving Theorem \ref{thm:Kernel_Gen2} }

\begin{lemma}\label{Lemma:Bochner}
Assume that $|k_{\mathcal{X}}|\leq 1$ and let  $\mu(\dd y)$  be a  $\sigma$-finite measure  on $(\setY,\mathfrak{S}_{\mathcal{Y}})$ and $f:\mathcal{X}\times\setY\rightarrow\R$ be such that
\begin{enumerate}
\item\label{B1} $f(\cdot,y)\in\mathcal{H}_{\mathcal{X}}$ for all $y\in\setY$,
\item\label{B2} The function  $\setY\ni y\mapsto f(\cdot,y)$ is Borel measurable,
\item\label{B3} The set $\{f(\cdot,y):\,\, y\in \setY\}$ is separable,
\item\label{B4}  $\int_\setY \|f( \cdot,y)\|_{\mathcal{H}_{\mathcal{X}}}\mu(\dd y)<\infty$.
\end{enumerate}
Then, $\int_{\setY}f(\cdot,y)\mu(\dd y)\in\mathcal{H}_{\mathcal{X}}$.
\end{lemma}

\begin{proof} Since the set $\{f(\cdot,y):\,\, y\in \setY\}$ is separable and  the mapping $y\mapsto f(\cdot,y)$  is Borel measurable the function $y\mapsto f(\cdot, y)$ is strongly measurable. Therefore, there exist  \citep[][Proposition E.2]{cohn} a sequence $\big( \{E_{i,n}\}_{i=1}^n\big)_{n\geq 1}$ and a sequence  $\big( \{f_{i,n}\}_{i=1}^n\big)_{n\geq 1}$ such that
\begin{enumerate}
\item $E_{i,n}\in \mathfrak{S}_{\mathcal{Y}}$ and $f_{i,n}\in\mathcal{H}_{\mathcal{X}}$ for all   $n\geq i\geq 1$,
\item $\lim_{n\rightarrow 0}\|\sum_{i=1}^n \ind_{E_{i,n}}(y)f_{i,n}-f(y,\cdot)\|_{\mathcal{H}_{\mathcal{X}}}= 0$ for all $y\in\setY$,
\item\label{fn2} $\|\sum_{i=1}^n \ind_{E_{i,n}}(y)f_{i,n}\|_{\mathcal{H}_{\mathcal{X}}}\leq \|f(y,\cdot)\|_{\mathcal{H}_{\mathcal{X}}}$ for all $n\geq 1$ and all $y\in\setY$.
\end{enumerate}

For every $n\geq 1$ let $f_n:\mathcal{X}\times\setY\rightarrow\R$ be defined by
$$
f_n(x,y)=\sum_{i=1}^n \ind_{E_{i,n}}(y)f_{i,n}(x),\quad (x,y)\in\mathcal{X}\times\setY.
$$
Under the assumptions of the lemma we have $\int_{\setY} \|f(\cdot,y)\|_{\mathcal{H}_{\mathcal{X}}}\mu(\dd y)<\infty$, and thus,
$$
\int_{\setY}\|f_n(\cdot,y)\|_{\mathcal{H}_{\mathcal{X}}}\mu(\dd y)\leq \int_{\setY} \|f(\cdot,y)\|_{\mathcal{H}_{\mathcal{X}}}\mu(\dd y)<\infty,\quad\forall n\geq 1,
$$
showing that, for all $n\geq 1$,  the simple function $y\mapsto f_n(\cdot,y)$ is Bochner integrable. Consequently, for all $n\geq 1$ the function
$$
\tilde{f}_n:=\int_{\setY}f_n(\cdot,y)\mu(\dd y)=\sum_{i=1}^n \Big(\int_{E_{i,n}} \mu(\dd y)\Big)f_{i,n}
$$  
is well-defined. Notice that  $\tilde{f}_n\in \mathcal{H}_{\mathcal{X}}$ for all $n\geq 1$.

To proceed further remark that 
$$
|f_n(x,y)|\leq \|f_n(\cdot,y)\|_{\mathcal{H}_{\mathcal{X}}}\leq  \|f(\cdot,y)\|_{\mathcal{H}_{\mathcal{X}}},\quad\forall (x,y)\in\mathcal{X}\times\setY
$$
where the first inequality holds since  $|k_{\mathcal{X}}|\leq 1$ by assumption while the second inequality hods by the third aforementioned properties of  $\big( \{E_{i,n}\}_{i=1}^n\big)_{n\geq 1}$ and   $\big( \{f_{i,n}\}_{i=1}^n\big)_{n\geq 1}$.

By assumption, $\int_{\setY}\|f(\cdot,y)\|_{\mathcal{H}_{\mathcal{X}}}\dd y<\infty$  and thus, by the dominated converge theorem, and using the fact that the convergence in $\|\cdot\|_{\mathcal{H}_{\mathcal{X}}}$ norm implies the point-wise convergence, 
\begin{align}\label{eq:Need0}
\lim_{n\rightarrow\infty}\tilde{f}_n(s)=\int_{\setY}f(s,y)\dd y,\quad\forall s\in\mathcal{X}.
\end{align}
Therefore, recalling that $\tilde{f}_n\in \mathcal{H}_{\mathcal{X}}$ for all $n\geq 1$, to complete the proof it remains to show that the sequence $(\tilde{f}_n)_{n\geq 1}$ is Cauchy w.r.t.\ the $\|\cdot\|_{\mathcal{H}_{\mathcal{X}}}$ norm.

To this aim   remark that, since
$$
\big\| f_n(\cdot,y)-f(\cdot,y) \big\|_{\mathcal{H}_{\mathcal{X}}}\leq 2\big\|f(\cdot,y) \big\|_{\mathcal{H}_{\mathcal{X}}},\quad\forall n\geq 1
$$
while, by assumption, $\int_\setY \|f(\cdot,y)\|_{\mathcal{H}_{\mathcal{X}}}\dd y<\infty$,
the dominated convergence theorem implies that
\begin{align}\label{eq:Need1}
\lim_{n\rightarrow \infty} \int_{\setY} \big\| f_n(\cdot,y)-f(\cdot,y) \big\|_{\mathcal{H}_{\mathcal{X}}}\dd y= 0.
\end{align}
On the other hand, for every $n>m\geq 1$ we have  
\begin{equation}\label{eq:Need2}
\begin{split}
\big\|\tilde{f}_n-\tilde{f}_m\big\|_{\mathcal{H}_{\mathcal{X}}}&=\Big\|\int_{\setY}\big\{f_n(\cdot,y)-f_m(\cdot,y)\big\}\mu(\dd y)\Big\|_{\mathcal{H}_{\mathcal{X}}}\\
&\leq \int_{\setY} \big\| f_n(\cdot,y)-f_m(\cdot,y) \big\|_{\mathcal{H}_{\mathcal{X}}}\mu(\dd y)\\
&\leq \int_{\setY} \big\| f_n(\cdot,y)-f(\cdot,y) \big\|_{\mathcal{H}_{\mathcal{X}}}\mu(\dd y)+\int_{\setY} \big\| f_m(\cdot,y)-f(\cdot,y) \big\|_{\mathcal{H}_{\mathcal{X}}}\mu(\dd y)
\end{split}
\end{equation}
where the first inequality holds by  \citet[][Proposition E.5]{cohn}, since a shown above the function $y\mapsto f_n(\cdot, y)$ is Bochner integrable. Together, \eqref{eq:Need1} and  \eqref{eq:Need2} show that the sequence $(\tilde{f}_n)_{n\geq 1}$ is indeed Cauchy w.r.t.\ the $\|\cdot\|_{\mathcal{H}_{\mathcal{X}}}$ norm, and the proof of the lemma is complete. 
\end{proof}

\subsection{Proof of Theorem \ref{thm:Kernel_Gen2} }
 
\begin{proof} Let $g\in\mathcal{H}_{\mathcal{Y}}$ so that $g=\sum_{i=1}^\infty a_i k_{\mathcal{Y}}(y_i,\cdot)$ for a sequence $(y_i)_{i\geq 1}$ in $\mathcal{Y}$ and a sequence $(a_i)_{i\geq 1}$ in $\R$. For all $n\geq 1$ let  $g_n=\sum_{i=1}^n a_ik_{\mathcal{Y}}(y_i,\cdot)$ and $f_n:\mathcal{X}\times\setY\rightarrow \R$ be defined by 
$f_n(x,y)=g_n(y)p(y|x) $, $(x,y)\in\mathcal{X}\times\setY$. We first show that, for all $n\geq 1$, the function $f_n$ verifies the assumptions of Lemma \ref{Lemma:Bochner}.

By Conditions \ref{CC1} and \ref{CC3} of the theorem,   it readily follows that $f_n$ verifies Conditions  \ref{B1} and \ref{B3} of Lemma \ref{Lemma:Bochner}, for all $n\geq 1$. To show that this is also the case for Condition \ref{B2} of Lemma \ref{Lemma:Bochner} let $\mathcal{B}\big(\mathcal{H}_{\mathcal{X}} \big)$ be the Borel $\sigma$-algebra on $\mathcal{H}_{\mathcal{X}}$. Let $n\geq 1$ and assume first that $\mathcal{H}_{\mathcal{X}}$ contains the non-zero constant functions so that  the function $y\mapsto  g_n(y)$ is  $\mathcal{B}(\mathcal{H}_{\mathcal{X}})$-measurable. Then, since by assumption the function  $y\mapsto p(y|\cdot)$ is    $\mathcal{B}(\mathcal{H}_{\mathcal{X}})$-measurable  and since the product of two Borel measurable functions is a Borel measurable function, it follows that the function $\setY\ni y\mapsto f_n(\cdot,y)$ is $\mathcal{B}(\mathcal{H}_{\mathcal{X}})$-measurable, as required. Assume now that $\mathcal{H}_{\mathcal{X}}$  does not contain the non-zero constant functions. Let $\tilde{H}_{\mathcal{X}}$ be the RKHS on $\mathcal{X}$ having $k_{\mathcal{X}}+1$ as reproducing kernel so that, as shown above,  the function $\setY\ni y\mapsto f_n(\cdot,y)$ is $\mathcal{B}\big(\tilde{\mathcal{H}}_{\mathcal{X}}\big)$-measurable. Consequently,
\begin{align}\label{eq:mes1}
\big\{y\in\setY:\, f_n(\cdot|y)\in A\big\}\in  \mathfrak{S}_\mathcal{Y},\quad\forall A\in  \mathcal{B}\big(\tilde{\mathcal{H}}_{\mathcal{X}}\big).
\end{align}
Recalling that $\tilde{\mathcal{H}}_{\mathcal{X}}=\big\{f+c,\, f\in \mathcal{H}_{\mathcal{X}},\,c\in\R\big\}$ and that $\|f\|_{\tilde{\mathcal{H}}_{\mathcal{X}}}= \|f\|_{\mathcal{H}_{\mathcal{X}}}$ for all $f\in \mathcal{H}_{\mathcal{X}}$ \citep[][Theorem 5.1]{paulsen2016}, it follows that  $\mathcal{B}(\mathcal{H}_{\mathcal{X}})\subset \mathcal{B}(\tilde{\mathcal{H}}_{\mathcal{X}})$ which, together with   \eqref{eq:mes1}, implies that
\begin{align*}
\big\{y\in\setY:\, f_n(\cdot|y)\in A\big\}\in  \mathfrak{S}_\mathcal{Y},\quad\forall A\in  \mathcal{B}\big(\mathcal{H}_{\mathcal{X}}\big).
\end{align*}
This shows that the function $\setY\ni y\mapsto f_n(\cdot,y)$ is  $\mathcal{B}(\mathcal{H}_{\mathcal{X}})$-measurable, and thus, for all $n\geq 1$, $f_n$ satisfies Condition \ref{B2} of  Lemma \ref{Lemma:Bochner}.

Lastly, using the fact that $|k_{\mathcal{Y}}|\leq 1$ and Condition \ref{CC4} of the theorem, for all $n\geq 1$ we have
\begin{align*}
\int_{\setY} \|f_n(\cdot ,y)\|_{\mathcal{H}_{\mathcal{X}}}\mu(\dd y)&\leq \big(\sup_{y\in\setY} |g_n(y)|\big)\int_{\mathcal{\mathcal{Y}}}  \|  p(y|\cdot)\|_{\mathcal{H}_{\mathcal{X}}}\mu(\dd y)\\
&\leq \|g_n\|_{\mathcal{H}_{\mathcal{Y}}}\int_{\mathcal{Y}}  \|  p(y|\cdot)\|_{\mathcal{H}_{\mathcal{X}}}\mu(\dd y)\\
&<\infty
\end{align*}
and thus, for all $n\geq 1$, $f_n$ verifies Condition \ref{B4} of Lemma \ref{Lemma:Bochner}, which concludes to show that, for all $n\geq 1$,   $f_n$ verifies all the assumptions of Lemma \ref{Lemma:Bochner}.

Therefore, by Lemma \ref{Lemma:Bochner}, the function $\tilde{f}_n:= \int_{\setY} f_n(\cdot,y)\mu(\dd y)$ exists and belongs to $\mathcal{H}_{\mathcal{X}}$, for all $n\geq 1$. In addition,  for all $m>n\geq 1$ we have  \citep[see][Proposition E.5, for the first inequality]{cohn}
\begin{align*}
\big\|\tilde{f}_n-\tilde{f}_m\big\|_{\mathcal{H}_{\mathcal{X}}}&=\Big\|\int_{\mathcal{Y}} (g_n-g_m)(y)  p(y|\cdot)\mu(\dd y)\Big\|_{\mathcal{H}_{\mathcal{X}}}\\
&\leq \int_{\mathcal{Y}}|g_n(y)-g_m(y)| \, \|   p(y|\cdot)\|_{\mathcal{H}_{\mathcal{X}}}\mu(\dd y)\\
&\leq \sup_{y\in\setY}|g_n(y)-g_m(y)|\int_{\mathcal{Y}} \big\|  p(y|\cdot)\big\|_{\mathcal{H}_{\mathcal{X}}}\mu(\dd y)
\end{align*}
where, since $|k_{\mathcal{Y}}|\leq 1$ by assumption,
\begin{align}\label{eq:unif}
\limsup_{n\rightarrow\infty}\sup_{m>n} \sup_{y\in\setY}|g_n(y)-g_m(y)|\leq  \limsup_{n\rightarrow\infty} \sup_{m>n}\|g_n-g_m\|_{\mathcal{H}_{\mathcal{Y}}}=0.
\end{align}
Consequently, the sequence $(\tilde{f}_n)_{n\geq 1}$ is Cauchy w.r.t.\ the $\|\cdot\|_{\mathcal{H}_{\mathcal{X}}}$ norm  and   therefore converges point-wise to a function $\tilde{f}\in\mathcal{H}_{\mathcal{X}}$. Thus, to complete the proof it remains to show that 
$$
\lim_{n\rightarrow\infty} \tilde{f}_n(x)=\E_{Y\sim P_{Y|X=x}}[g(Y)],\quad\forall x\in\mathcal{X}.
$$
Since for every $n\geq  1$ and $x\in\mathcal{X}$ we have
\begin{align*}
\big|\tilde{f}_n(x)-\E_{Y\sim P_{Y|X=x}}[g(Y)]\big|\leq \int_{\mathcal{Y}} |g_n(y)-g(y)|\, p(y|x)\mu(\dd y)\leq \sup_{y\in\setY}|g_n(y)-g(y)|,
\end{align*}
it follows, by  \eqref{eq:unif}, that $
\lim_{n\rightarrow\infty}\sup_{x\in\mathcal{X}}|\tilde{f}_n(x)-\E_{Y\sim P_{Y|X=x}}[g(Y)]|=0$,
and the proof of the theorem is complete.
\end{proof}

\section{Proof of Corollary \ref{cor:Kernel_Gen2}}

Corollary \ref{cor:Kernel_Gen2} is a direct consequence of  Lemma \ref{lemma:CBM_New}, Theorem \ref{thm:Kernel_Gen2} and of the following lemma:

\begin{lemma}\label{lemma:bound_norm}
 Assume that Assumptions \ref{asm:k}-\ref{asm:kxy} hold  and that there exists a $\sigma$-finite measure $\mu(\dd y)$   on $(\setY,\mathfrak{S}_{\mathcal{Y}})$ such that  $P_{Y|x}=p(y|x) \mu(\dd y)$ for all $x\in\mathcal{X}$, where $p(\cdot|\cdot)$ satisfies Assumptions \ref{CC1}-\ref{CC4} of Theorem \ref{thm:Kernel_Gen2}. Moreover, assume that there exists a bounded conditional mean embedding operator $\mathcal{C}_{Y|X}$ for $(P_{Y|x})_{x\in\mathcal{X}}$. Then,    $\|\mathcal{C}_{Y|X}\|_{\mathrm{o}}\leq \int_{\mathcal{Y}} \|p(y|\cdot)\|_{\mathcal{H}_{\mathcal{X}}}\mu(\dd y)$.
\end{lemma}
\begin{proof}
Let  $g\in\mathcal{H}_{\mathcal{Y}}$ and remark that 
$$
(\mathcal{C}_{Y|X}^* g)(x)=<\mathcal{C}_{Y|X}^* g, k_{\mathcal{X}}(x,\cdot)>_{\mathcal{H}_{\mathcal{X}}}=<g,\mathcal{C}_{Y|X} k_{\mathcal{X}}(x,\cdot)>_{\mathcal{H}_{\mathcal{Y}}}=\E_{Y\sim P_{Y|X=x}}[g(Y)],\quad\forall x\in\mathcal{X}
$$
where the first equality uses the reproducing property of $k_{\mathcal{X}}$ and the third \eqref{eq:r1}.

Consequently,
\begin{align*}
\|\mathcal{C}_{Y|X}^* g\|_{\mathcal{H}_{\mathcal{X}}}&=\Big\|\int_{\mathcal{Y}} g(y)p(y|\cdot)\dd y\Big\|_{\mathcal{H}_{\mathcal{X}}}\\
&\leq  \int_{\mathcal{Y}} \|g(y)p(y|\cdot)\|_{\mathcal{H}_{\mathcal{X}}}\dd y\\
&\leq \sup_{y\in\mathcal{Y}}|g(y)|\int_{\mathcal{Y}} \|p(y|\cdot)\|_{\mathcal{H}_{\mathcal{X}}}\dd y\\
&\leq \|g\|_{\mathcal{H}_{\mathcal{Y}}}\int_{\mathcal{Y}} \|p(y|\cdot)\|_{\mathcal{H}_{\mathcal{X}}}\dd y
\end{align*}
where, under the assumptions of the lemma, the first inequality holds by \citet[][Proposition E.5]{cohn} and where the last inequality uses the fact that $|k_{\mathcal{Y}}|\leq 1$.

Therefore,
$$
\|\mathcal{C}_{Y|X}\|_{\mathrm{o}}=\|\mathcal{C}^*_{Y|X}\|_{\mathrm{o}}\leq \int_{\mathcal{Y}} \|p(y|\cdot)\|_{\mathcal{H}_{\mathcal{X}}}\dd y
$$
and the proof of the lemma is complete.
\end{proof}

\section{A useful corollary of  Theorem \ref{thm:Kernel_Gen2}\label{p_corr_usefuke}}
In order to state the next result we let $\Lambda_d(\dd x)$ denote  the Lebesgue measure on $\R^d$, $A_s= \big\{\tilde{a}\in\mathbb{N}_0^{d}:\,\sum_{i=1}^{d}\tilde{a}_i\leq s\big\}$ for all $s\in\mathbb{N}_0$ and   $|a|=\sum_{i=1}^d a_i$ for all $a\in\R^d$.

\begin{corollary}\label{cor:CME}
Assume that $\mathcal{X}\subseteq\R^d$   is bounded with Lipschitz boundary and, for  some constants $m\in\mathbb{N}$ and $\gamma>0$, let $k_{\mathcal{X}}$ be the restriction of the Mat\'ern kernel $K_{\frac{m}{2},\gamma}$ on $\mathcal{X}\times \mathcal{X}$. Let  $s=(d+m)/2$ if $(d+m)$ is even and $s=(d+m+1)/2$ if $(d+m)$ is odd, and assume that there exists a $\sigma$-finite measure $\mu(\dd y)$   on $(\setY,\mathfrak{S}_{\mathcal{Y}})$ such that  $P_{Y|x}=p(y|x) \mu(\dd y)$ for all $x\in\mathcal{X}$, where $p(\cdot|\cdot)$ satisfies the following conditions:
\begin{itemize}
\item for all $y\in\mathcal{Y}$, the function $p(y|\cdot)$ is $s$ times continuously differentiable on $\setX$, with
$$
\max_{a\in A_s}\sup_{(x,y)\in\mathcal{X}\times\mathcal{Y}} \Big|\frac{\partial^{\sum_{i=1}^d a_i}}{\partial x_1^{a_1}\dots\partial x_d^{a_d}} p(y| x)\Big|<\infty 
$$
and with
$$
\max_{a\in A_s}\int_{ \mathcal{Y}} \left[\int_{\mathcal{X}}    \Big\{\frac{\partial^{\sum_{i=1}^d a_i}}{\partial x_1^{a_1}\dots\partial x_d^{a_d}} p(y| x)\Big\}^2 \Lambda_d(\dd x)\right]^{\frac{1}{2}}\mu(\dd y)<\infty,\quad\forall a\in A_s,
$$
\item the function $y\mapsto \frac{\partial^{\sum_{i=1}^d a_i}}{\partial x_1^{a_1}\dots\partial x_d^{a_d}} p(y| x)$ is continuous on $\mathcal{Y}$,  for all $x\in\mathcal{X}$ and $a\in A_s$.
\end{itemize}
Assume also that the set $\mathcal{Y}$ is separable and that Assumptions \ref{asm:k}-\ref{asm:kxy} hold. Then, conditions \ref{CC1}-\ref{CC4} of Theorem \ref{thm:Kernel_Gen2} hold  and thus   \eqref{eq:cond_embed} is satisfied.
\end{corollary}

\begin{proof}  Remark first that to prove the result it is enough to consider the case where  $(m+d)$ is  even. Indeed, if $(m+d)$ is odd then in what follows we  can replace
\begin{itemize}
\item  the set $\setX$ by  $\tilde{\setX}=\setX\times \R^{1}$,
\item for all $y\in\mathcal{Y}$, the function $p(y|\cdot):\mathcal{X}\rightarrow\R$ by the function $\tilde{p}(y|\cdot):\tilde{\setX}\rightarrow\R$  defined by $\tilde{p}(y|  (x,v'))=p(y|x)$ for all $(x,u)\in\tilde{\setX}$,
\item  $d$ by $\tilde{d}=d+1$.
\end{itemize}

Recall that, since $(m+d)$ is  even,  the RKHS $\mathcal{H}_{\mathcal{X}}$  is norm-equivalent to the Sobolev space $W_2^s\big(\mathcal{X})$ \citep[see e.g.][Example 2.6]{kanagawa2018gaussian}. In addition, recall that   the norm $\|\cdot\|_{W_2^s ( \mathcal{X} )}$ is defined by
$$
\|f\|_{W_2^s ( \mathcal{X})}=\sum_{a\in A_s}\bigg(\int_{\mathcal{X}}\left|\frac{\partial^{\sum_{i=1}^d a_i}}{\partial u_1^{a_1}\dots\partial u_d^{a_d}} f(x)\right|^2\Lambda_d(\dd x)\bigg)^{\frac{1}{2}},\quad f\in W_2^s\big(X) 
$$
and let
$$
D_a p(y|x)=\frac{\partial^{\sum_{i=1}^d a_i}}{\partial x_1^{a_1}\dots\partial x_d^{a_d}} p(y| x),\quad\forall (a,x,y)\in A_s\times\mathcal{X}\times\mathcal{Y}.
$$

To proof the corollary remark first that, under its assumptions, for all $y\in\mathcal{Y}$ we have $\|p(\cdot|y)\|_{W_2^s ( \mathcal{X})}<\infty$. Thus, for all $y\in\setY$, the function $p(y|\cdot)$ belongs to the Sobolev space $W_2^s\big(\mathcal{X}\big)$, and thus  to the RKHS $\mathcal{H}_{\mathcal{X}}$. This shows that  $p(\cdot|\cdot )$ verifies Condition \ref{CC1}  of Theorem \ref{thm:Kernel_Gen2}. 

In addition,  under the assumptions of the corollary we have 
\begin{equation}\label{eq:norm_p}
\begin{split}
\int_\setY \|p(y|\cdot)\|_{W_2^s(\mathcal{X})} \mu(\dd y)=\sum_{a\in A_s}\int_\setY \left\{\int_{\mathcal{X}}  D_a p(y|x)^2\Lambda_d(\dd x) \right\}^{\frac{1}{2}} \mu(\dd y)
&<\infty
\end{split}
\end{equation}
and thus,   since the norm $\|\cdot\|_{\mathcal{H}_{\mathcal{X}}}$  is equivalent to the norm $\|\cdot\|_{W_2^s(\mathcal{X})}$,  it follows that    
$$
\int_\setY\|p(y|\cdot)\|_{\mathcal{H}_{\mathcal{X}}}\mu(\dd y)<\infty
$$
showing that  $p(\cdot|\cdot )$ verifies   Condition \ref{CC4}  of Theorem \ref{thm:Kernel_Gen2}. 
 
To proceed further  recall  that the image of a separable space by a continuous function is separable. Hence, since $\mathcal{Y}$ is assumed to be separable, to show that Condition \ref{CC3} of Theorem \ref{thm:Kernel_Gen2}  holds it suffices to show that, for every $y'\in\setY$, the function 
\begin{align}\label{eq:map}
\setY\ni y\mapsto k_{\mathcal{Y}}(y',y) p(y|\cdot)\in \mathcal{H}_{\mathcal{X}}
\end{align}
 is continuous. To this aim, let   $y'\in\setY$ and $(y'_i)_{i\geq 1}$ be a sequence in $\mathcal{Y}$ such that $\lim_{i\rightarrow\infty}y'_i=y'$. Then,  since    $k_{\mathcal{Y}}$ is continuous  by assumption, to show that the function defined in \eqref{eq:map} is continuous it is enough to show that
 \begin{align}\label{eq:to_showseq}
 \limsup_{i\rightarrow\infty}\|p(y'_i|\cdot)-p(y'|\cdot)\|_{\mathcal{H}_{\mathcal{X}}}=0.
 \end{align}
The   norm $\|\cdot\|_{\mathcal{H}_{\mathcal{X}}}$ being norm-equivalent to the norm  $\|\cdot\|_{W_2^s(\mathcal{X})}$, there exists a constant $C<\infty$ such that $\|f\|_{\mathcal{H}_{\mathcal{X}}}\leq C\|f \|_{W_2^s (\mathcal{X} )}$ for all $f\in \mathcal{H}_{ \mathcal{X}}$ and thus, for all $i\geq 1$, we have
 \begin{equation}\label{eq:seq_1}
 \begin{split}
 \|p(y'_i|\cdot)-p(y'|\cdot)\|^2_{\mathcal{H}_{\mathcal{X}}}&\leq  C^2\|p(y'_i|\cdot)-p(y'|\cdot)\|^2_{W_2^s(\mathcal{X})}\\
 &\leq C^2 \sum_{a\in A_s}\int_{\mathcal{X}}\big|D_a p(y'_i|x)-D_a p(y'|x)\big|^2\Lambda_d(\dd x).
 \end{split}
 \end{equation}
 By assumption,  for all $x\in \mathcal{X}$ and  all $a\in A_s$, the function $y\mapsto D_a p(y|x)$
is continuous on $\setY$ while,  for all $(a,x)\in A_s\times \mathcal{X}$ we have
\begin{align*}
\sup_{i\geq 1}\big|D_a p(y'_i|x)-D_a p(y'|x)\big|^2&\leq 2\sup_{(x,y)\in\mathcal{X}\times\mathcal{Y}} |D^a p(y|x)|<\infty.
\end{align*}
Consequently, since $\mathcal{X}$ is bounded,  \eqref{eq:to_showseq} follows from  \eqref{eq:seq_1}  and   the dominated convergence theorem, and thus $p(\cdot|\cdot )$ satisfies Condition \ref{CC3} of Theorem \ref{thm:Kernel_Gen2}.

Finally,  since  as shown above  the mapping $\setY\ni y\mapsto p(y|\cdot)$ is continuous, it follows that this mapping is Borel measurable and thus $p(\cdot|\cdot )$ satisfies Condition \ref{CC2} of Theorem \ref{thm:Kernel_Gen2}. Hence, all the conditions of Theorem \ref{thm:Kernel_Gen2} and the proof is complete.
\end{proof}

\section{Proof of Proposition \ref{prop:CME:example}}

\subsection{Preliminary result}

\begin{lemma}
\label{lemma:bij}
Assume that $\mathcal{X} \subseteq \mathbb{R}^{d}$ for some integer $d$ and that $\mathcal{X}$ is path-wise connected and such that $\Lambda_d(\mathcal{X})>0$. Assume also that $k_{\mathcal{X}}$ is continuous on $\mathcal{X}^2$. Then, there exists a distribution  $P_X\in\mathcal{P}(\setX)$ such that
\begin{align}
\label{eq:bij}
 \big\{f \in\mathcal{H}_{\mathcal{X}}:\,\E_{X\sim P_X}\big[f(X)  h(X)\big]=0, \,\forall h\in\mathcal{X}\big\}=\{0\}.
 \end{align}
\end{lemma}

\begin{proof}

Remark first that since $k_{\mathcal{X}}$ is continuous on $\mathcal{X}^2$ any function $f\in \mathcal{H}_{\mathcal{X}}$ is continuous on $\mathcal{X}$ \citep[][Theorem 2.17]{paulsen2016}. Let $P_X$ denote the $\mathcal{N}_d(0,I_d)$ distribution, truncated on $\mathcal{X}$ if $\mathcal{X} \neq \mathbb{R}^d$. Assume that there exists a non-zero function $f\in\mathcal{H}_{\mathcal{X}}$   such that 
$$
\E_{X\sim P_X}\big[f(X)h(X)\big]=0, \quad \,\forall h\in\mathcal{H}_{\mathcal{X}}.
$$
Then, $\E_{X\sim P_X}[f(X)^2]=0$ and, since $P_X$ admits a strictly positive density $p_X$ on $\mathcal{X}$ w.r.t.~$\Lambda_d$, we have $f(x)=0$ for $\Lambda_d$-almost every $x\in\mathcal{X}$. However, as $f$ is assumed to be continuous, and $\mathcal{X}$ is path-wise connected, the function $f$ is zero everywhere.
\end{proof}

\subsection{Proof of the proposition}

\textit{Proof:}

The fact that $k$ is characteristic follows from \citet{szabo2018characteristic} and the properties of the Mat\'ern kernel.
 
Next, remark that, by Lemma \ref{lemma:bij},  under the assumption of the proposition   there exists a distribution $P_X\in\mathcal{P}(\mathcal{X})$ such that the only function $f\in\mathcal{H}_{\mathcal{X}}$ for which we have $\E_{X\sim P_X}[f(X)^2]=0$ is the zero function. In addition,   since $ \mathcal{X}$ is bounded with Lipschitz boundary  we can use Corollary \ref{cor:CME} to check that  there exists a bounded linear conditional mean operator     for $(P_{Y|x})_{x\in \mathcal{X}}$.

To this aim, for all  $x\in\mathcal{X}$ we let $p(y|x)$ be the density of $P_{Y|x}(\dd y)$ w.r.t.~$\mu(\dd y)$. The $\sigma$-finite measure $\mu(\dd y)$ on $\mathcal{Y}$ will be specified below for each example but, for all the considered examples,   for all $y\in\Theta\times \mathcal{Y}$ the mapping $x\mapsto p(y|\cdot)$ is infinitely many times differentiable. Consequently,  letting $s$ and $A_s$ be as defined in Corollary \ref{cor:CME}, we can define
$$
D^a p(y|x)= \frac{\partial^{\sum_{i=1}^d a_i}}{\partial x_1^{a_1}\dots\partial x_d^{a_d}} p(y| x) ,\quad \forall (a,x,y)\in \Theta\times A_s\times\mathcal{X}\times\mathcal{Y}.
$$
Then, by  Corollary \ref{cor:CME},  a  bounded linear conditional mean operator   for $(P_{Y|x})_{x\in \mathcal{X}}$ exists if
\begin{enumerate}
\item\label{continuity}  the mapping $y\mapsto D^a p(y|x)$ is continuous for all $(a,x)\in A_s\times \mathcal{X}$,

\item   the following two conditions hold: 
\begin{align}
&\max_{a\in A_s}\sup_{(x,y)\in\mathcal{X}\times\mathcal{Y}} |D^a p(y|x)|<\infty\label{eq:bound_mod} \\
&\max_{a\in A_s}\int_{ \mathcal{Y}} \left[\int_{\mathcal{X}}    \{D^a p(y|x)\}^2 \Lambda_d(\dd x)\right]^{\frac{1}{2}}\mu(\dd y)<\infty.\label{eq:bound_mod2}
\end{align}
\end{enumerate}

For all the examples considered in the proposition it is trivial to see that the mapping $y\mapsto D^a p(y|x)$ is continuous for all $(a,x)\in  A_s\times \mathcal{X}$. Under the assumptions made on   $\mathcal{X}$,  Conditions  \eqref{eq:bound_mod} and \eqref{eq:bound_mod2} are easily checked from the definition of $p(y|x)$ given below for each examples

\paragraph{Example 1:} For this example $\mathcal{Y}=\R$ and we let  $\mu(\dd y)$ be the Lebesgue measure on $\R$ so that
$$
p(y|x)=\sum_{m=1}^M w_m\frac{1}{\sqrt{2\pi\sigma_m^2}}\exp\Big\{-\frac{(y-\beta_m^\top x)^2}{2\sigma_m^2}\Big\},\quad\forall (x,y)\in\mathcal{X}\times\mathcal{Y}.
$$

\paragraph{Example 2:}  For this example, $\mathcal{Y}=\mathbb{N}_0$ and we let $\mu(\dd y)$ be the counting measure on $\mathbb{N}_0$ so that
$$
p(y|x)= \frac{\exp\big\{y\,\beta^\top x- \exp( \beta^\top x)\big\}}{y!},\quad\forall (x,y)\in \mathcal{X}\times\mathcal{Y}.
$$

\paragraph{Example 3:} For this example, $\mathcal{Y}=\{0,1\}$ and we let $\mu(\dd y)$ be the counting measure on $\{0,1\}$ so that
$$
p(y|x)=\left(\frac{1}{ 1+\exp(-\beta^\top x)}\right)^y  \left(\frac{1}{ 1+\exp(\beta^\top x)}\right)^{1-y},\quad\forall (x,y)\in\mathcal{X}\times\mathcal{Y}.
$$

\paragraph{Example 4:} For this example, $\mathcal{Y}=(0,\infty)$ and we let $\mu(\dd y)$ be the Lebesgue measure on $\R$ so that 
$$
p(y|x)= \frac{1}{\Gamma(\nu)} y^{\nu-1} \exp(-\nu \beta^\top x)\exp\big\{-\nu y\exp(-\beta^\top x)\big\},\quad\forall (x,y)\in\mathcal{X}\times\mathcal{Y}.
$$

\paragraph{Example 5:} For this example, $\mathcal{Y}=\R\times\{0,1\}$ and we let 
$$
 \mu(\dd y)=\big(\Lambda_1(\dd y_1)+\delta_{\{0\}}(\dd y_1))\otimes \delta_{\{0\}}(\dd y_2)
$$
so that
$$
P_\lambda(\dd y)=\check{p}_\lambda(y)\mu(\dd y)
$$
where, denoting by $\phi(\cdot;\mu,\sigma^2)$ the probability density function of the $\mathcal{N}_1(\mu,\sigma^2)$  distribution w.r.t.\ $\Lambda_1$, for all $(x,y)\in\mathcal{X}\times\mathcal{Y}$ we have
\begin{align*}
p(y|x)&=\phi\big(y_1;\beta^\top x, \sigma^2\big)\Phi\Big(\{\gamma^\top x+(\rho/\sigma)\beta^\top x\}/\sqrt{1-\rho^2}\Big)\ind_{\R\setminus\{0\}}(y_1)\big(1-\ind_{\{0\}}(y_2)\big)\\
 &+\Phi(-\gamma^\top x)\ind_{\{0\}}(y_1)\ind_{\{0\}}(y_2).
 \end{align*} 
 \hfill$\square$

\section{Proof of Lemma \ref{lemma:constant:kernel}\label{p-lemma:constant:kernel}}

\begin{proof} Let $\mathcal{C}_{Y|X}:\mathcal{H}_{\mathcal{X}}\rightarrow  \mathcal{H}_{\mathcal{Y}}$ be a bounded linear operator such that
\begin{align}\label{eq:r1}
\mu(P_{Y|x})=\mathcal{C}_{Y|X} k_{ \mathcal{X}}(x,\cdot),\quad\forall x\in\mathcal{X}.
\end{align}
and let $\tilde{\mathcal{C}}_{P_{Y|X}}:\mathcal{H}_{\mathcal{X}}\otimes\mathcal{H}_{\mathcal{X}}\rightarrow \mathcal{H}$ be the (unique) linear operator on $\mathcal{H}_{\mathcal{X}}\otimes\mathcal{H}_{\mathcal{X}}$ such that
$$
\tilde{\mathcal{C}}_{P_{Y|X}}(f_1\otimes f_2)=f_1\otimes \mathcal{C}_{Y|X} f_2,\quad f_1\in\mathcal{H}_{\mathcal{X}},f_2\in\mathcal{H}_{\mathcal{X}}. 
$$
For all $f_1\in \mathcal{H}_{\mathcal{X}}$ and $f_2\in\mathcal{H}_{\star,\mathcal{X}}$ we have 
\begin{align*}
\|\tilde{\mathcal{C}}_{P_{Y|X}}(f_1\otimes f_2)\|_{\mathcal{H}}
&=\|f_1\otimes \mathcal{C}_{Y|X}f_2\|_{\mathcal{H}}\\
&=\|f_1\|_{\mathcal{H}_{\mathcal{X}}}\|\mathcal{C}_{Y|X} f_2\|_{\mathcal{H}_{ \mathcal{Y}}}\\
 &\leq \|f_1\|_{\mathcal{H}_{\mathcal{X}}} \|f_2\|_{\mathcal{H}_{ \mathcal{X}}}\|\mathcal{C}_{Y|X}\|_{\mathrm{o}}\\
 &=\|f_1\otimes f_2\|_{\mathcal{H} }\|\mathcal{C}_{Y|X}\|_{\mathrm{o}}
\end{align*}
showing that
\begin{equation}\label{eq:r2}
\|\tilde{\mathcal{C}}_{P_{Y|X}}\|_{\mathrm{o}}\leq \|\mathcal{C}_{Y|X}\|_{\mathrm{o}}<\infty.
\end{equation}
 where the last inequality holds by assumption.

Next, remark that for every $f\in\mathcal{H}_{\mathcal{X}}$ the linear operator $(f \otimes \cdot ):\mathcal{H}_{\mathcal{Y}}\rightarrow\mathcal{H}$ is such that
\begin{align}\label{eq:r3}
\|f \otimes \cdot\|_{\mathrm{o}}\leq \|f \|_{\mathcal{H}_{\mathcal{X}}}<\infty.
\end{align}
since
$$
\|f\otimes g\|_{\mathcal{H}}=\|f\|_{\mathcal{H}_{\mathcal{X}}}\|g\|_{\mathcal{H}_{\mathcal{Y}}},\quad\forall f\in\mathcal{H}_{\mathcal{X}},\,\forall g\in\mathcal{H}_{\mathcal{Y}}.
$$

Let $\tilde{\mu}(P_X)=\mathbb{E}_{X\sim P_X}\big[k_{\mathcal{X}}(X,\cdot)\otimes k_{\mathcal{X}}(X,\cdot)]$ be the embedding  of  $P_X\in\mathcal{P}(\mathcal{X})$ in $\mathcal{H}_{\mathcal{X}}\otimes\mathcal{H}_{\mathcal{X}}$

Then, for $P'_X\in\mathcal{P}(\setX)$ and using the shorthand $P'=P'_XP_{Y|\cdot}$, we have
\begin{align*}
\mu(P'):&=\mathbb{E}_{(X,Y)\sim P'}\big[k_{\mathcal{X}}(X,\cdot)\otimes k_{\mathcal{Y}}(Y,\cdot)]\\
&=\mathbb{E}_{X\sim P'_X}\Big[\mathbb{E}_{Y\sim P_{Y|X}}\big[ k_{\mathcal{X}}(X,\cdot)\otimes k_{\mathcal{Y}}(Y,\cdot)\big]\Big]\\
&=\mathbb{E}_{X\sim P'_X}\Big[k_{\mathcal{X}}(X,\cdot)\otimes \mathbb{E}_{Y\sim P_{Y|X}}\big[ k_{\mathcal{Y}}(Y,\cdot)\big]\Big]\\
&=\mathbb{E}_{X\sim P'_X}\Big[k_{\mathcal{X}}(X,\cdot)\otimes \mu(P_{Y|X})\Big]\\
&=\mathbb{E}_{X'\sim P'_X}\Big[k_{\mathcal{X}}(X,\cdot)\otimes \mathcal{C}_{Y|X}k_{\mathcal{X}}(X,\cdot)\Big]\\
&=\mathbb{E}_{X'\sim P'_X}\Big[\tilde{\mathcal{C}}_{P_{Y|X}}\Big(k_{\mathcal{X}}(X',\cdot)\otimes  k_{ \mathcal{X}}(X',\cdot)\Big)\Big]\\
&=\tilde{\mathcal{C}}_{P_{Y|X}}\mathbb{E}_{X\sim P'_X}\big[k_{\mathcal{X}}(X,\cdot)\otimes  k_{\mathcal{X}}(X,\cdot)\big]\\
&= \tilde{\mathcal{C}}_{P_{Y|X}}\tilde{\mu}(P'_X)
\end{align*}
where the  interchange  between expectation and tensor product between the second and the third equality is  justified by Lemma \ref{lemma:swap} and by \eqref{eq:r3}, where the  interchanges  between expectation and tensor product between the fifth and the sixth   equality is justified by Lemma \ref{lemma:swap} and by  \eqref{eq:r2}, while the fifth  equality holds by \eqref{eq:r1}.

Similarly, for   $P''_X\in\mathcal{P}(\setX)$ and with $P''=P'_XP_{Y|\cdot}$, we have 
\begin{align*}
\mu(P''):&=\mathbb{E}_{(X,Y)\sim P''}\big[k_{\mathcal{X}}(X,\cdot)\otimes k_{\mathcal{Y}}(Y,\cdot)]\tilde{\mathcal{C}}_{P_{Y|X}}\tilde{\mu}(P''_X)
\end{align*}
 and thus,  
\begin{equation}\label{eq:final}
\begin{split}
\mathbb{D}_{k}(P', P'')&=\|\mu(P')-\mu(P'')\|_{\mathcal{H}}
\\
&= \Big\|\tilde{\mathcal{C}}_{P_{Y|X}}\Big(\tilde{\mu}(P_X')-\tilde{\mu}(P_X'')\Big)\Big\|_{\mathcal{H}}
\\
&\leq  \|\tilde{\mathcal{C}}_{P_{Y|X}}\|_{\mathrm{o}}\,\big\|\tilde{\mu}(P_X')-\tilde{\mu}(P_X'')\big\|_{\mathcal{H}_{\mathcal{X}} \otimes\mathcal{H}_{{\mathcal{X}}}}
\\
 &\leq \|\mathcal{C}_{Y|X}\|_{\mathrm{o}}\,\,\mathbb{D}_{k^2_{\mathcal{X}}}(P'_X, P''_X)
\end{split}
\end{equation}
where the last inequality holds by \eqref{eq:r2}. The proof is complete.

\end{proof}

\section{Proof of Lemma~\ref{thm:fixed1}}
\label{subsection:proofs:theorems1}

\subsection{Preliminary results}
\label{subsection:proofs:lemma1}

The following lemma is adapted from Lemma 5 in~\cite{cherief2019mmd}. While the proof is quite similar, the statement is more general.

\begin{lemma}
\label{lemma:preliminary:MMD}
 Let $\mathcal{S}$ be a set (equipped with a $\sigma$-algebra). Let $K$ be any symmetric function $\mathcal{S}^2 \rightarrow [-1,1]$ that can be written $K(s,s') = \left<\varphi(s),\varphi(s')\right>_{\mathcal{H}}$ for some Hilbert space $\mathcal{H}$ and some function $\varphi$ (note that we do not assume that $K$ is a characteristic kernel). Let $S_1,\dots,S_n$ be independent random variables on $\mathcal{S}$ with respective distributions $Q_1,\dots,Q_n$. Define $\bar{Q}=(1/n)\sum_{i=1}^n Q_i$ and $\hat{Q}=(1/n)\sum_{i=1}^n \delta_{S_i} $. We define, for any $Q$ and $Q'$ probability distributions on $\mathcal{S}$,
 $$ \mathbb{D}^2_K(Q,Q') = \mathbb{E}_{S\sim Q,S'\sim Q}[K(Z,Z')] - 2  \mathbb{E}_{S\sim Q,S'\sim Q'}[K(Z,Z')] +  \mathbb{E}_{S\sim Q',S'\sim Q'}[K(Z,Z')] $$
 (which is indeed a metric if $K$ is a characteristic kernel). We have:
 $$ \mathbb{E}\left[ \mathbb{D}_K(\bar{Q},\hat{Q}) \right] \leq \frac{1}{\sqrt{n}} \text{ and } \mathbb{E}\left[ \mathbb{D}_K^2(\bar{Q},\hat{Q}) \right] \leq \frac{1}{n}. $$
\end{lemma}

\begin{proof}  Jensen's inequality gives $ \mathbb{E}[ \mathbb{D}_K(\bar{Q},\hat{Q}) ] \leq \sqrt{\mathbb{E}[\mathbb{D}_K^2(\bar{Q},\hat{Q}) ]} $. Put $m_i = \mathbb{E}_{S\sim Q_i}[\varphi(S)]$, then
\begin{align*}
\mathbb{E}\left[ \mathbb{D}_K^2(\bar{Q},\hat{Q}) \right]
 & = \mathbb{E} \left[  \left\| \frac{1}{n}\sum_{i=1}^{n}\left[\varphi(S_i) - m_i \right] \right\|_{\mathcal{H}}^2 \right]
\\
 & = \frac{1}{n^2}\sum_{i=1}^n \mathbb{E} \left[ \left\| \varphi(S_i) - m_i \right\|_{\mathcal{H}}^2 \right] + \frac{1}{n(n-1)} \sum_{i\neq j} \mathbb{E} \left[  \left< \varphi(S_i) - m_i, \varphi(S_j) - m_j \right>_{\mathcal{H}} \right]
\\
& = \frac{1}{n^2}\sum_{i=1}^n \Bigl( \mathbb{E} \left[ \left\| \varphi(S_i) \right\|_{\mathcal{H}}^2 \right] - \left\| m_i\right\|_{\mathcal{H}}^2 \Bigr) + 0
\\
& \leq \frac{1}{n^2}\sum_{i=1}^n \mathbb{E} \left[  \left\| \varphi(S_i) \right\|_{\mathcal{H}}^2 \right]
= \frac{1}{n^2} \sum_{i=1}^n K(S_i,S_i) \leq \frac{1}{n}\text{. }
\end{align*}
\end{proof}

 Our proof strategy to study $\hat{\theta}_n(D_n)$ actually relies on the fact that despite contamination, the performance of $\hat{\theta}_n(D_n)$ remains close to the one of $\hat{\theta}_n(D_n)$. The following lemma will help to formalize this claim.

\begin{lemma}
 \label{lemma:preliminary:MMD:2}
 Let 
 $ \hat{P}^{n,0} = \frac{1}{n}\sum_{i=1}^n \delta_{(X_i^0,Y_i^0)}  $
 be the non-contaminated empirical distribution and
 $ \hat{P}^{n,0}_\theta = \frac{1}{n}\sum_{i=1}^n \delta_{X_i^0} P_{g(\theta,X_i^0)}   $
 be the  uncontaminated counterpart of $\hat{P}^{n}_\theta$. Then, for any probability distribution $Q$ on $\mathcal{X}\times\mathcal{Y}$, we have
 \begin{equation}
\label{equa:proof:adversarial}
\left|  \mathbb{D}_k\left( \hat{P}^{n,0},Q \right) -  \mathbb{D}_k\left( \hat{P}^n,Q \right) \right| < 2\epsilon
\end{equation}
and
\begin{equation}
\label{equa:proof:adversarial:2}
\left|  \mathbb{D}_k\left(\hat{P}^{n,0}_\theta ,Q \right) -  \mathbb{D}_k\left( \hat{P}^{n}_\theta,Q \right) \right| < 2\epsilon.
\end{equation}
\end{lemma}

\begin{proof} 
For the first inequality~\eqref{equa:proof:adversarial},
\begin{align*}
\left|  \mathbb{D}_k\left( \hat{P}^{n,0},Q \right) -  \mathbb{D}_k\left( \hat{P}^n,Q \right) \right|
& \leq 
  \mathbb{D}_k\left( \hat{P}^{n,0},\hat{P}^n\right) 
\\
& = \left\| \frac{1}{n}\sum_{i=1}^n \Bigl[ k((X_i^0,Y_i^0),\cdot) - k((X_i,Y_i),\cdot) \Bigr] \right\|_{\mathcal{H}}
\\
& \leq \frac{1}{n}\sum_{i=1}^n \left\|  \Bigl[ k((X_i^0,Y_i^0),\cdot) - k((X_i,Y_i),\cdot) \Bigr] \right\|_{\mathcal{H}}
\\
& = \frac{1}{n}\sum_{i\in I} \left\|  \Bigl[ k((X_i^0,Y_i^0),\cdot) - k((X_i,Y_i),\cdot) \Bigr] \right\|_{\mathcal{H}}
\\
& \leq \frac{1}{n} \sum_{i\in I} 2 = \frac{2 |I| }{n } < 2\epsilon .
\end{align*}
The proof of~\eqref{lemma:preliminary:MMD:2} is exactly the same.
\end{proof}

\subsection{Proof of the lemma}

\begin{proof} 
Thanks to~\eqref{equa:proof:adversarial} of Lemma~\ref{lemma:preliminary:MMD:2} we have, for any fixed $\theta\in\Theta$,
\begin{align}
 \mathbb{D}_k(\hat{P}^n_{\hat{\theta}_n} ,\bar{P}_n^0)
 & \leq
 \mathbb{D}_k(\hat{P}^n_{\hat{\theta}_n}, \hat{P}^{n}) + \mathbb{D}_k( \hat{P}^{n},\bar{P}_n^0) \text{ (triangle inequality)} \nonumber
 \\
 & \leq
  \mathbb{D}_k(\hat{P}^n_{\hat{\theta}_n}, \hat{P}^{n}) + \mathbb{D}_k( \hat{P}^{n,0},\bar{P}_n^0) +2 \epsilon \text{ where we used~\eqref{equa:proof:adversarial} with }Q=\bar{P}_n^0 \nonumber
\\ 
 & \leq
  \mathbb{D}_k(\hat{P}^n_{\theta}, \hat{P}^{n}) + \mathbb{D}_k( \hat{P}^{n,0},\bar{P}_n^0) +2 \epsilon \text{ (by definition of } \hat{\theta}_n \text{)} \nonumber
\\
& \leq
  \mathbb{D}_k(\hat{P}^n_{\theta}, \hat{P}^{n,0}) + \mathbb{D}_k( \hat{P}^{n,0},\bar{P}_n^0) +4 \epsilon \text{ by~\eqref{equa:proof:adversarial} with }Q=\hat{P}^n_{\theta} \nonumber
\\
& \leq   \mathbb{D}_k(\hat{P}^n_{\theta}, \bar{P}_n^0) + 2 \mathbb{D}_k( \hat{P}^{n,0},\bar{P}_n^0) +4 \epsilon
 \text{ (triangle inequality)}.  \label{eq:1:lemma:preliminary:MMD}
\end{align}
Taking the expectation in~\eqref{eq:1:lemma:preliminary:MMD} gives:
\begin{align}\label{eq:thm:fixed1:e1}
 \mathbb{E}\left[  \mathbb{D}_k(\hat{P}^n_{\hat{\theta}_n} ,\bar{P}^0_n ) \right] \leq 4\epsilon +\mathbb{D}_k(\hat{P}^n_{\theta},\bar{P}^0_n ) + 2 \mathbb{E}\left[ \mathbb{D}_k(\hat{P}^{n,0},\bar{P}^0_n ) \right]. 
 \end{align}
We can control the expectation in the right-hand side by an application of Lemma~\ref{lemma:preliminary:MMD}, where  $S_i=(X_i^0,Y_i^0)\sim Q_i: = \delta_{X_i^0} P^0_{Y|X_i^0} $ that are indeed independent, and where $K=k$. The lemma gives:
\begin{align}\label{eq:thm:fixed1:e}
\mathbb{E}\left[ \mathbb{D}_k(\hat{P}^n,\bar{P}^0_n ) \right] \leq\frac{1}{\sqrt{n}}.
\end{align}
We take the infinimum with respect to $\theta$ to obtain:
\begin{align}\label{eq:thm:fixed1:e2}
 \mathbb{E}\left[  \mathbb{D}_k(\hat{P}^n_{\hat{\theta}_n} ,\bar{P}^0_n ) \right] \leq 4\epsilon +\inf_{\theta\in\Theta} \mathbb{D}_k(\hat{P}^n_{\theta},\bar{P}^0_n ) + \frac{2}{\sqrt{n}}. 
 \end{align}

In order to prove~\eqref{eq:tthm:fixed1:p}, take any $z_i'\in\mathcal{Z}$ and define  
$$
\hat{P}^{n,0}_{(i)}=\frac{1}{n} \Big(\sum_{j\neq i}  \delta_{(X_j^0,Y_j^0)} + \delta_{z_i'}\Big).
$$ We note that:
$$ \left|\mathbb{D}_k(\hat{P}^{n,0},\bar{P}^0_n ) - \mathbb{D}_k(\hat{P}^{n,0}_{(i)},\bar{P}^0_n )\right|
\leq \mathbb{D}_k(\hat{P}^{n,0},\hat{P}^{n,0}_{(i)}) \leq \frac{2}{n} .$$
This allows to use the McDiarmind's bounded difference inequality~\cite{mcdiarmid1989method}, which gives:
\begin{align}\label{eq:thm:fixed1:p}
\mathbb{P}\left\{  \mathbb{D}_k(\hat{P}^n,\bar{P}^0_n ) - \mathbb{E}\left[ \mathbb{D}_k(\hat{P}^n,\bar{P}^0_n) \right]  \geq t \right\} \leq \exp\left(-  \frac{n t^2}{2} \right),\quad\forall t>0.
\end{align} 
Put $\eta=\exp(-n t^2 /2)$ to get
$$ \mathbb{P}\left\{  \mathbb{D}_k(\hat{P}^{n,0},\bar{P}^0_n ) - \mathbb{E}\left[ \mathbb{D}_k(\hat{P}^{n,0},\bar{P}^0_n) \right]  \geq \sqrt{\frac{2 \log(1/\eta)}{n}} \right\} \leq \eta, $$
which, together with \eqref{eq:thm:fixed1:e1}-\eqref{eq:thm:fixed1:e}, gives the statement of the theorem. 
\end{proof}

\section{Proof of Theorem \ref{cor:fixed}}

\subsection{Preliminary result}

\begin{lemma}
 \label{lemma:MMD:euclid}
 Let $\|\cdot\|$ be a semi-norm on $\Theta$.
 Let $M:\Theta\rightarrow [0,2]$ be such that there exists a unique $\theta_\star\in\Theta$ verifying $\inf_{\theta\in\Theta}M(\theta)=M(\theta_\star)$  and such that there exists a neighborhood $U$ of $\theta_\star$ and a constant $\mu>0$ for which
\begin{align*}
M(\theta)-M(\theta_\star)\geq \mu\|\theta-\theta_\star\|,\quad\forall \theta\in U.
\end{align*}
 Let $(\check{\theta}_n)_{n\geq 1}$ be a sequence of random variables taking values in $\Theta$ and such  that there exist  a strictly increasing function $h_1:(0,\infty)\rightarrow (0,\infty)$ with $\lim_{x\rightarrow \infty} h_1(x) = \infty$, a continuous and strictly decreasing function $h_2:(0,1)\rightarrow (0,\infty)$, and a constant $x \geq 0$  such that 
 \begin{equation}
\label{eqn:proof:proba}
\mathbb{P}\Big\{M(\check{\theta}_n) < M(\theta_\star)+x+\frac{h_2(\eta)}{h_1(n)}\Big\}\geq 1-\eta,\quad\forall \eta\in(0,1),\quad\forall n\geq 1.
\end{equation}
Then for any $t>0$,
$$
 \mathbb{P}\Big\{  \|\check{\theta}_n-\theta_\star\| \geq x/\mu + t  \Big\} \leq 2 h_2^{-1}\left[ \left((\mu t) \wedge (\alpha-x)_+ \right)h_1(n) \right],
$$
and
$$
\mathbb{P}\Big\{  \|\check{\theta}_n-\theta_\star\| < \frac{x}{\mu}+   \frac{h_2\left( \frac{\eta}{2}\right)}{\mu h_1(n)}  \Big\} \geq 1-\eta ,\quad\forall n \geq 1,\quad\forall \eta \in \left[ 2 h_2^{-1}((\alpha-x)_+ h_1(n)),1 \right)
$$
where $\alpha = \inf_{\theta\in U^c}   M(\theta)-M(\theta_\star) \in (0,2]$.
\end{lemma}

\begin{remark}
It would also be possible to get a result on $ \mathbb{E}[ \|\check{\theta}_n-\theta_\star\|] $, but at the price of the additional assumption that the parameter space $\Theta$ is bounded: $\sup_{(\theta,\theta')\in\Theta^2}\|\theta-\theta'\|_\Theta <\infty$.
\end{remark}

\begin{proof} 
Note that~\eqref{eqn:proof:proba} is equivalent to
 \begin{equation}
\label{eqn:proof:proba:2}
\mathbb{P}\Big\{M(\check{\theta}_n)-  M(\theta_\star) -x> t\Big\}\leq h_2^{-1}(th_1(n)) ,\quad\forall t>0,\quad\forall n\geq 1.
\end{equation}
Remind that $\alpha = \inf_{\theta\in U^c}   M(\theta)-M(\theta_\star)  $. It is immediate to see that $\alpha \leq 2$. Moreover, $\alpha>0$, otherwise, $U^c$ being a closed set, there would be a $\theta'\in U^c$ such that $M(\theta')-M(\theta_\star)=0$.

Now, for any $t>0$,
\begin{align*}
 \mathbb{P} & \left\{  \|\check{\theta}_n-\theta_\star\| \geq t + x/\mu \right\}
 \\
 & =  \mathbb{P}\left\{  \|\check{\theta}_n-\theta_\star\| \geq t + x/\mu , \check{\theta}_n \in U \right\} + \mathbb{P}\left\{  \|\check{\theta}_n-\theta_\star\| \geq t + x/\mu , \check{\theta}_n \notin U \right\}
 \\
 & \leq \mathbb{P}\left\{  M(\check{\theta}) - M(\theta_\star)  \geq \mu t +x , \check{\theta}_n \in U \right\} + \mathbb{P}\left\{ \check{\theta}_n \notin U \right\}
 \\
 & \leq \mathbb{P}\left\{  M(\check{\theta}) - M(\theta_\star) -x \geq \mu t \right\} + \mathbb{P}\left\{  M(\check{\theta}) - M(\theta_\star)  \geq \alpha \right\}
 \\
 & \leq h_2^{-1}\left( \mu t h_1(n) \right) + h_2^{-1}\left( (\alpha-x)_+ h_1(n) \right)
\end{align*}
where we used~\eqref{eqn:proof:proba:2} for the last inequality. As $h_2^{-1}$ is strictly decreasing, we obtain:
\begin{equation}
\label{eq:lemma5:step}
\mathbb{P}\left\{  \|\check{\theta}_n-\theta_\star\| \geq t  + x/\mu \right\} \leq 2 h_2^{-1}\left[ \left((\mu t) \wedge (\alpha-x)_+ \right)h_1(n) \right] .
\end{equation}

Fix $\eta \in \left[ 2 h_2^{-1}((\alpha-x)_+ h_1(n)),1 \right)$ as in the statement of the lemma, and note that
$$ 2 h_2^{-1}\left[ \left((\mu t) \wedge (\alpha-x)_+ \right)h_1(n) \right] = \eta \Leftrightarrow t = \frac{h_2\left(\frac{\eta}{2}\right)}{\mu h_1(n)} . $$
Plugging these values in~\eqref{eq:lemma5:step}, we obtain:
$$
\mathbb{P}\Big\{  \|\check{\theta}_n-\theta_\star\| < \frac{x}{\mu}+   \frac{h_2\left( \frac{\eta}{2}\right)}{\mu h_1(n)}  \Big\} \geq 1-\eta.
\text{ } $$
\end{proof}

\subsection{Proof of the theorem}

\begin{proof} 

From Lemma~\eqref{thm:fixed1}, \eqref{eqn:proof:proba} in Lemma~\ref{lemma:MMD:euclid} holds with $\theta_\star = \theta_0$, $x=4\epsilon$, $h_1(n)=\sqrt{n}$, $h_2(\eta) = 2+\sqrt{2\log(1/\eta)}$ and $\check{\theta}_n = \hat{\theta}_n$ . Apply  Lemma~\ref{lemma:MMD:euclid} to get:
$$
\sum_{n\geq 1} \mathbb{P}\left\{  \|\check{\theta}_n-\theta_\star\| \geq +\frac{4\epsilon}{\mu}+ t  \right\} \leq 2 \sum_{n\geq 1} \exp\left[- \frac{\left[\left((\mu t) \wedge (\alpha-x)_+ \right)\sqrt{n} - 2\right]^2}{2} \right]<\infty,\quad\forall t>0
$$
showing that  $\mathbb{P}\big(\limsup_{n\rightarrow\infty}\|\check{\theta}_n-\theta_\star\|\leq 4\epsilon/\mu\big)=1$.
Lemma~\ref{lemma:MMD:euclid} also states
$$
\mathbb{P}\Big\{  \|\check{\theta}_n-\theta_\star\| <  \frac{h_2\left( \frac{\eta}{2}\right)}{\mu h_1(n)}  \Big\} \geq 1-\eta ,\quad\forall n \geq 1,\quad\forall \eta \in \left[ 2 h_2^{-1}((\alpha-x)_+ h_1(n)),1 \right).
$$
Note that
$$\frac{h_2\left( \frac{\eta}{2}\right)}{\mu h_1(n)}  = \frac{1}{\mu \sqrt{n}} \left(2+\sqrt{2\log(2/\eta)}\right) $$
and $2 h_2^{-1}((\alpha-x)_+ h_1(n)) = 2 \exp(- ((\alpha-x)_+ \sqrt{n}-2)^2 / 2)$. For the sake of simplicity, we only consider $n\geq 16/(\alpha-x)_+^2 $, in this case, we have $(\alpha-x)_+ \sqrt{n}-2 \geq (\alpha-x)_+ \sqrt{n}/2$ and thus the result holds in particular for any $\eta\in [ 2\exp(-n(\alpha-x)_+^2 /8),1)$. Finally, remind that $x=4\epsilon< \alpha/8$ so it holds in particular for $n\geq 64/\alpha^2 $ and $\eta\in [ 2\exp(-n \alpha^2 /32),1)$.
\end{proof}

\section{Proof of Lemma~\ref{thm:random1}}\label{p-thm:random1}

\subsection{Preliminary result}

We start with a result that will be an essential tool in the proof of Lemma \ref{thm:random1}. Essentially, it quantifies how well $ \hat{P}^{n,0}_{\hat{\theta}_n} = (1/n)\sum_{i=1}^n \delta_{X_i^0} P_{g(\hat{\theta}_n,X_i^0)} $ approximates $P^0$. Usually, in regression literature, we focus mostly on the estimation of the distribution of $Y|X$ rather than on the estimation of the distribution of the pair $(X,Y)$. Still, we believe that this result has in interpretation on its own, so we state is as a theorem.

\begin{theorem}\label{thm:random0}
 Under Assumption~\ref{asm:k} we have
 $$ \mathbb{E}\left[  \mathbb{D}_k(\hat{P}^{n,0}_{\hat{\theta}_n},P^0 ) \right] \leq 8\epsilon + \inf_{\theta\in\Theta}  \mathbb{D}_k(P_{\theta},P^0) + \frac{3}{\sqrt{n}} $$
 and, for any $\eta\in(0,1)$,
  $$ \mathbb{P}\left\{  \mathbb{D}_k(\hat{P}^{n,0}_{\hat{\theta}_n},P^0 ) \leq 8\epsilon + \inf_{\theta\in\Theta}  \mathbb{D}_k(P_{\theta},P^0) + \frac{3}{\sqrt{n}}\left(1+\sqrt{2\log(2/\eta)} \right) \right\} \geq 1-\eta. $$
\end{theorem}

\begin{proof}  The proof is quite similar to the proof of Lemma~\ref{thm:fixed1}, but requires some adaptations, in particular in the application of Lemma~\ref{lemma:preliminary:MMD}.

First,
\begin{equation}
   \mathbb{D}_k(\hat{P}^{n,0}_{\hat{\theta}_n},P^0 ) \leq \mathbb{D}_k(\hat{P}^{n,0}_{\hat{\theta}_n},\hat{P}^{n,0} ) + \mathbb{D}_k(\hat{P}^{n,0},P^0 ) \label{veryfirststep}.
\end{equation}
Let us deal with the first term of this upper bound in a first time. Here, we will use both~\eqref{equa:proof:adversarial} and~\eqref{equa:proof:adversarial:2} of Lemma~\ref{lemma:preliminary:MMD:2}. We have:
\begin{align*}
 \mathbb{D}_k(\hat{P}^{n,0}_{\hat{\theta}_n},\hat{P}^{n,0} )
   \leq  \mathbb{D}_k\left( \hat{P}^{n,0}_{\hat{\theta}_n},\hat{P}^{n} \right) + 2\epsilon   \leq  \mathbb{D}_k\left(  \hat{P}^{n}_{\hat{\theta}_n},\hat{P}^{n} \right) + 4\epsilon 
 & \leq \mathbb{D}_k\left( \hat{P}^{n}_{\hat{\theta}_n(D_n^0)},\hat{P}^{n} \right) + 4\epsilon
 \\
 & \leq \mathbb{D}_k\left(  \hat{P}^{n,0}_{\hat{\theta}_n(D_n^0)},\hat{P}^{n} \right) + 6 \epsilon 
 \\
 & \leq \mathbb{D}_k\left( \hat{P}^{n,0}_{\hat{\theta}_n(D_n^0)},\hat{P}^{n,0} \right) + 8 \epsilon 
 \\
 & = \inf_{\theta\in\Theta} \mathbb{D}_k\left( \hat{P}^{n,0}_{\theta},\hat{P}^{n,0} \right)+ 8 \epsilon.
\end{align*}
where the first inequality uses \eqref{equa:proof:adversarial}, the second   \eqref{equa:proof:adversarial:2}, the third   the definition of $\hat{\theta}_n$, the fourth \eqref{equa:proof:adversarial:2}, the fifth  \eqref{equa:proof:adversarial} and the sixth the definition of $\hat{\theta}_n$.

Together with \eqref{veryfirststep}, this shows that
\begin{align}\label{eq:proof:thm0:1}
  \mathbb{D}_k(\hat{P}^{n,0}_{\hat{\theta}_n},P^0 ) 
& \leq \inf_{\theta\in\Theta} \mathbb{D}_k\left( \hat{P}^{n,0}_{\theta},\hat{P}^{n,0} \right)+ \mathbb{D}_k(\hat{P}^{n,0},P^0 )  + 8 \epsilon 
 \nonumber
 \\ & \leq \inf_{\theta\in\Theta}  \mathbb{D}_k\left( \hat{P}^{n,0}_{\theta},P^0\right) + 2 \mathbb{D}_k(\hat{P}^{n,0},P^0 ) + 8\epsilon
 \nonumber
 \\ & \leq \inf_{\theta\in\Theta} \left[   \mathbb{D}_k\left( \hat{P}^{n,0}_{\theta},P_\theta \right)+
 \mathbb{D}_k(P_\theta,P^0 )  \right]+ 2 \mathbb{D}_k(\hat{P}^{n,0},P^0 ) + 8\epsilon
\end{align}
and so, taking expectations on both sides,
\begin{equation}
\begin{split}
\label{multline-step1}
\mathbb{E}\left[ \mathbb{D}_k(\hat{P}^{n,0}_{\hat{\theta}_n},P^0 ) \right]
&\leq
\inf_{\theta\in\Theta} \left\{ \mathbb{E}\left[\mathbb{D}_k(\hat{P}^{n,0}_{\theta},P_\theta  )\right] +  \mathbb{D}_k(P_\theta,P^0 )  \right\} 
 2\mathbb{E}\left[\mathbb{D}_k(\hat{P}^{n,0},P^0 ) \right] + 8\epsilon.
\end{split}
\end{equation}
We tackle the term $ \mathbb{E}\left[\mathbb{D}_k(\hat{P}^{n,0}_{\theta},P_\theta  )\right] $.
Letting $\Phi$ denote the function such that $k((x,y),(x',y'))=\left<\Phi(x,y),\Phi(x',y')\right>_\mathcal{H}$, we have
\begin{align*}
 \mathbb{D}_k( \hat{P}^{n,0}_{\theta},P_\theta  )
  &
  = \sqrt{ \mathbb{D}_k^2(\hat{P}^{n,0}_{\theta},P_\theta  )}
  \\
  &
  = \Biggl(\mathbb{E}_{(X,Y)\sim \hat{P}^{n,0}_{\theta},(X',Y')\sim \hat{P}^{n,0}_{\theta} } \left< \Phi(X,Y),\Phi(X',Y')\right>_{\mathcal{H}}
  \\
  & \quad -2 \mathbb{E}_{(X,Y)\sim \hat{P}^{n,0}_{\theta},(X',Y')\sim P_\theta} \left< \Phi(X,Y),\Phi(X',Y')\right>_{\mathcal{H}}
  \\
  & \quad +\mathbb{E}_{(X,Y)\sim P_\theta,(X',Y')\sim P_\theta} \left< \Phi(X,Y),\Phi(X',Y')\right>_{\mathcal{H}} \Biggl)^{\frac{1}{2}}
  \\
  &   = \Biggl(\mathbb{E}_{X\sim \frac{1}{n}\sum_{i=1}^n \delta_{X_i^0} ,X'\sim \frac{1}{n}\sum_{i=1}^n \delta_{X_i^0}} \left< \mathbb{E}_{Y\sim P_{g(\theta,X)}} [\Phi(X,Y)] ,  \mathbb{E}_{Y'\sim P_{g(\theta,X')}} [\Phi(X',Y')]\right>_{\mathcal{H}}
  \\
  & \quad -2 \mathbb{E}_{X\sim \frac{1}{n}\sum_{i=1}^n \delta_{X_i^0} ,X'\sim P_X^0} \left< \mathbb{E}_{Y\sim P_{g(\theta,X)}} [\Phi(X,Y)] ,  \mathbb{E}_{Y'\sim P_{g(\theta,X')}} [\Phi(X',Y')]\right>_{\mathcal{H}}
  \\
  & \quad +\mathbb{E}_{X\sim P_X^0,X'\sim P_X^0} \left< \mathbb{E}_{Y\sim P_{g(\theta,X)}} [\Phi(X,Y)] ,  \mathbb{E}_{Y'\sim P_{g(\theta,X')}} [\Phi(X',Y')]\right>_{\mathcal{H}} \Biggl)^{\frac{1}{2}}
  \\
  & = \sqrt{\mathbb{D}_{\bar{k}}^2\left( \frac{1}{n}\sum_{i=1}^n \delta_{X_i^0} ,P_X^0 \right)} = \mathbb{D}_{\bar{k}}\left( \frac{1}{n}\sum_{i=1}^n \delta_{X_i^0} ,P_X^0 \right)
\end{align*}
where the function $\bar{k}$ is given by:
$$ \bar{k}(x,x') = \left< \mathbb{E}_{Y\sim P_{g(\theta,x)}} [\Phi(x,Y)] ,  \mathbb{E}_{Y'\sim P_{g(\theta,x')}} [\Phi(x',Y')]\right>_{\mathcal{H}} .$$
Note that $-1\leq \bar{k} \leq 1$ so we can apply Lemma~\ref{lemma:preliminary:MMD} to $S_i=X_i^0\sim Q_i = P_X^0$ and $K=\bar{k}$ to get:
$$ \mathbb{E} \left[ \mathbb{D}_{\bar{k}}\left( \frac{1}{n}\sum_{i=1}^n \delta_{X_i^0} ,P_X^0 \right)\right] \leq \frac{1}{\sqrt{n}} .$$
Combining this last result with  \eqref{multline-step1}, and applying Lemma~\ref{lemma:preliminary:MMD} with $S_i=(X_i^0,Y_i^0)\sim Q_i = P^0 $ and $K=k$ that gives  $\mathbb{E}[ \mathbb{D}_k(\hat{P}^{n,0},P^0 )]\leq 1/\sqrt{n}$, we finally obtain:
\begin{align*}
\mathbb{E}\left[\mathbb{D}_k(\hat{P}^{n,0}_{\hat{\theta}_n},P^0 ) \right]
& \leq
\inf_{\theta\in\Theta} \left\{ \frac{1}{\sqrt{n}} +  \mathbb{D}_k(P_\theta,P^0 )  \right\} + \frac{2}{\sqrt{n}} + 8\epsilon
\\
& = \inf_{\theta\in\Theta}  \mathbb{D}_k(P_\theta,P^0 )  + \frac{3}{\sqrt{n}} + 8\epsilon,
\end{align*}
that is the first inequality of the theorem.

In order to prove the second inequality let $\theta_0\in\argmin_{\theta\in\Theta} \mathbb{D}_k(P_\theta,P^0 )  $. Then~\eqref{eq:proof:thm0:1} implies
\begin{equation*}
\begin{split}
\mathbb{D}_k & (\hat{P}^{n,0}_{\hat{\theta}_n},P^0 )
\\
&
\leq
\mathbb{D}_k( \hat{P}^{n,0}_{\theta_0} , P_{\theta_0} )+  \mathbb{D}_k(P_{\theta_0},P^0 )   + 2 \mathbb{D}_k(\hat{P}^{n,0},P^0 ) + 8\epsilon
\\
& = \mathbb{D}_k( \hat{P}^{n,0}_{\theta_0} , P_{\theta_0} ) + \inf_{\theta\in\Theta} \mathbb{D}_k(P_{\theta},P^0 )    + 2\mathbb{D}_k(\hat{P}^{n,0},P^0 ) + 8\epsilon.
\end{split}
\end{equation*}
McDiarmid's bounded difference inequality leads to
$$\mathbb{P}\left\{  \mathbb{D}_k(\hat{P}^{n,0},P^0 ) - \mathbb{E}\big[ \mathbb{D}_k(\hat{P}^{n,0},P^0 ) \big]  \geq t \right\} \leq \exp\left(-  \frac{n t^2}{2} \right)$$
and to
$$\mathbb{P}\left\{   \mathbb{D}_k(\hat{P}^{n,0},P^0 )- \mathbb{E}\left(  \mathbb{D}_k(\hat{P}^{n,0},P^0 ) \right)  \geq t \right\} \leq \exp\left(-  \frac{n t^2}{2} \right).$$
By a union bound, the probability that one of the two events hold is smaller or equal to $2 \exp(-nt^2/2)$, which leads to
  $$ \mathbb{P}\left\{  \mathbb{D}_k(\hat{P}^{n,0}_{\hat{\theta}_n},P^0 ) \leq \inf_{\theta\in\Theta}  \mathbb{D}_k(P_\theta,P^0) + \frac{3}{\sqrt{n}}\left(1+\sqrt{2\log(2/\eta)} \right) + 8\epsilon \right\} \geq 1-\eta \text{. }  $$
This ends the proof. 

\end{proof}

  \subsection{Proof of the lemma}

\begin{proof} 

By Lemma \ref{lemma:constant:kernel} applied to $P_X' = \hat{P}^{n,0}_X$ and $P_X'' = P^0_X$, we have 
\begin{align}
\label{eq:lemma:kernel:1}
\mathbb{D}_k(\hat{P}^n_{\hat{\theta}_n} ,P_{\hat{\theta}_n} )\leq \mathfrak{C}\, \mathbb{D}_{k^2_X}(\hat{P}_X^{n,0}, P^0_X)
\end{align}
 and thus
\begin{align*}
\mathbb{E}\Big[ \mathbb{D}_k(\hat{P}^{n,0}_{\hat{\theta}_n} ,P_{\hat{\theta}_n} )\Big] \leq   \mathfrak{C}\mathbb{E}\Big[\mathbb{D}_{k^2_X}(\hat{P}_X^{n,0}, P^0_X) \Big].
\end{align*}
Applying Lemma~\ref{lemma:preliminary:MMD}  with $Z_i=X_i\sim Q_i = P^0_X$ and $K=k^2_{\mathcal{X}}$, we obtain
\begin{equation}
\label{eq:lemma:kernel:2}
\mathbb{E}\Big[ \mathbb{D}_k(\hat{P}^{n,0}_{\hat{\theta}_n} ,P_{\hat{\theta}_n} )\Big] \leq   \frac{ \mathfrak{C}} {\sqrt{n}}.
\end{equation}
Now:
\begin{align*}
 \mathbb{E}\left[  \mathbb{D}_k(P_{\hat{\theta}_n},P^0 ) \right]
 & \leq  \mathbb{E}\left[  \mathbb{D}_k(P_{\hat{\theta}_n},\hat{P}_{\hat{\theta}_n}) \right] +  \mathbb{E}\left[  \mathbb{D}_k(\hat{P}_{\hat{\theta}_n},P^0 ) \right]
 \\
 & \leq  \frac{ \mathfrak{C}} {\sqrt{n}}+ \left( \inf_{\theta\in\Theta}  \mathbb{D}_k(P_\theta,P^0) +8 \epsilon +  \frac{3}{\sqrt{n}}  \right) 
\end{align*}
where we used~\eqref{eq:lemma:kernel:2} to upper bound the first term, and Theorem~\ref{thm:random0} for the second term. This ends the proof of the bound in expectation.

Let us now prove the inequality in probability. Let $\eta\in(0,1)$ and use the bounded difference inequality to get
$$ \mathbb{P}\bigg\{ \mathbb{D}_{k_{\mathcal{X}}^2}(\hat{P}_X^{n,0}, P^0_X) - \mathbb{E}\Big[ \mathbb{D}_{k^2_X}(\hat{P}_X^{n,0}, P^0_X) \Big]  \leq \sqrt{\frac{ 2\log(2/\eta)}{n}} \bigg\} \geq 1-\frac{\eta}{2} $$
while, by Theorem \ref{thm:random0},
$$
\mathbb{P}\bigg\{ \mathbb{D}_k(\hat{P}^{n,0}_{\hat{\theta}_n},P^0 ) \leq 8 \epsilon + \inf_{\theta\in\Theta}  \mathbb{D}_k(P_{\theta},P^0) + \frac{3}{\sqrt{n}}\left(1+\sqrt{2\log(4/\eta)} \right) \bigg\} \leq 1-\frac{\eta}{2}. $$
Together with~\eqref{eq:lemma:kernel:1}, and using a union bound, we obtain
\begin{multline*}
\mathbb{P}\bigg\{  \mathbb{D}_k(\hat{P}_{\hat{\theta}_n},P^0 ) \leq  \inf_{\theta\in\Theta}  \mathbb{D}_k(P_{\theta},P^0) + \frac{ 3\big(1+\sqrt{2\log(4/\eta)}\big)+\mathfrak{C}\big(1+\sqrt{2\log(2/\eta)}\big)}{\sqrt{n}}  \bigg\}\\
\geq 1-\eta.
\end{multline*}
\end{proof}

\section{Proof of Theorem~\ref{cor:random}}
\begin{proof} 
From Lemma~\ref{thm:random1}, \eqref{eqn:proof:proba} in Lemma~\ref{lemma:MMD:euclid} holds with $h_1(n)=\sqrt{n}$, $h_2(\eta) = (\mathfrak{C}  + 3 ) (1+\sqrt{2\log(4/\eta)})$ and $\check{\theta}_n = \hat{\theta}_n$. Then, the result is proved following the computations done in the proof of Theorem \ref{cor:fixed}.
\end{proof}

\section{Proof of Proposition \ref{prop:Taylor}}

\begin{proof} 

 Let $f:\Theta\rightarrow[0,4]$ be defined by 
$$
f(\theta)=\big(\mathbb{D}_k(P_{\theta},\tilde{P}^0  )- \mathbb{D}_k(P_{\theta_0},\tilde{P}^0  )\big)^2,\quad \theta\in\Theta
$$
and let $U$ be an open set containing $\theta_0$ such that $f$ is twice continuously differentiable on $U$. Let $H_\theta$ be the Hessian matrix of $f$ evaluated at $\theta\in U$.

Then, using Taylor's theorem, for every $\theta\in U$ we have, for some $\tau\in[0,1]$
\begin{align*}
f(\theta)&=f(\theta_0)+(\theta-\theta_0)^\top\nabla f(\theta_0)+\frac{1}{2}(\theta-\theta_0)^\top H_{\theta_0+\tau(\theta-\theta_0)} (\theta-\theta_0)\\
&=(\theta-\theta_0)^\top H_{\theta_0+\tau(\theta-\theta_0)} (\theta-\theta_0)\\
&\geq \|\theta-\theta_0\|^2\frac{\lambda_{\mathrm{min}}\big(H_{\theta_0+\tau(\theta-\theta_0)}\big)}{2}\\
&\geq \|\theta-\theta_0\|^2\frac{\inf_{\theta\in U,\tau\in[0,1]} \lambda_{\mathrm{min}}\big(H_{\theta_0+\tau(\theta-\theta_0)}\big)}{2}
\end{align*}
where for every $\theta\in U$ we denote by $\lambda_{\mathrm{min}}(H_\theta)$ the minimum eigenvalue of $H_\theta$. Under the assumptions of the proposition, we can take $U$ sufficiency small so that
$
c:=\inf_{\theta\in U,\tau\in[0,1]} \lambda_{\mathrm{min}}\big(H_{\theta_0+\tau(\theta-\theta_0)}\big)
>0$. Then,
\begin{align*}
\mathbb{D}_k(P_{\theta},\tilde{P}^0 )- \mathbb{D}_k(P_{\theta_0},\tilde{P}^0)=\sqrt{f(\theta)}\geq \sqrt{c/2}\, \|\theta-\theta_0\|
\end{align*}
showing that \eqref{eq:loweDk} holds for $\mu=\sqrt{c/2}$.
\end{proof}

\section{Proof of Proposition \ref{prop:random2}}
\label{subsection:proofs:theorems2}

\begin{proof} 

 For all $(\theta,x,y)\in\Theta\times\setX\times\setY$, let
$$
m_\theta(x,y)=\mathbb{E}_{Y, Y'\iid P_{g(\theta,X)} } \big[k_{\mathcal{Y}}(Y,Y') -2 k_{\mathcal{Y}}(Y,y)\big] +\E_{X\sim P_X^0}\Big[\mathbb{E}_{Y, Y'\iid P^0_{Y|X} } \big[k_{\mathcal{Y}}(Y,Y')\big]\Big]
$$
and remark that
\begin{align*}
\E_{(X,Y)\sim P^0}&[m_\theta(X,Y)]=\E_{X \sim P_X^0}\big[ \mathbb{D}_{k_{\mathcal{Y}}}(P_{g(\theta,X)}, P^0_{Y|X})^2\big],\quad\forall\theta\in\Theta.
\end{align*}
Under the assumptions of the theorem, the  mapping $\theta\mapsto m_\theta(x,y)$ is continuous on the compact set $\Theta$ and is such that $|m_\theta(x,y)|\leq  4$ for all $(\theta,x,y)\in\Theta\times\setX\times\setY$. Then  \citep[see e.g][page 46]{Vaart2000}
$$
\sup_{\theta\in\Theta}\Big|\frac{1}{n}\sum_{i=1}^n  m_\theta(X_i,Y_i)-\E_{X \sim P_X^0}\big[ \mathbb{D}_{k_{\mathcal{Y}}}(P_{g(\theta,X)}, P^0_{Y|X})^2\big]\Big|\rightarrow 0,\quad \text{in $\mathbb{P}$-probability}
$$
and therefore, noting that $
\tilde{\theta}_n\in\argmin_{\theta\in\Theta} \frac{1}{n}\sum_{i=1}^n m_\theta(X_i,Y_i)$,
the result follows by \citet[][Theorem 5.7]{Vaart2000}.
\end{proof}

\section{Proof of Theorem~\ref{thm:random2}}

\begin{proof} 
 Let $\epsilon\in[0,1)$ and, for  all   $x\in\setX$, let $\tilde{P}^0_{Y|x}=(1-\epsilon)P^0_{Y|x}+\epsilon Q_{Y|x}$ and $\tilde{P}^0_X=(1-\epsilon) P^0_X+\epsilon Q_X$ where $Q_X$ denotes the distribution of $X$ under $Q$. 

Then, for all $\theta\in\Theta$ we have
\begin{equation}\label{eq:decompos}
\begin{split}
\E_{X \sim P_X^0}\big[ \mathbb{D}_{k_{\mathcal{Y}}} (P_{g(\theta,X)},& P^0_{Y|X})^2\big]\\
&\leq  \E_{X \sim P_X^0}\big[ \big(\mathbb{D}_{k_{\mathcal{Y}}}(P_{g(\theta,X)},  \tilde{P}^0_{Y|X})+\mathbb{D}_{k_{\mathcal{Y}}}(P^0_{Y|X}, \tilde{P}^0_{Y|X})\big)^2\big]\\
&\leq  \E_{X \sim P_X^0}\big[ \big(\mathbb{D}_{k_{\mathcal{Y}}}(P_{g(\theta,X)},  \tilde{P}^0_{Y|X})+2\epsilon\big)^2\big]\\
&\leq \E_{X \sim P_X^0}\big[ \big(\mathbb{D}_{k_{\mathcal{Y}}}(P_{g(\theta,X)},  \tilde{P}^0_{Y|X})^2\big]+8 \epsilon +4\epsilon^2\\
&\leq \E_{X \sim P_X^0}\big[ \big(\mathbb{D}_{k_{\mathcal{Y}}}(P_{g(\theta,X)}, \tilde{P}^0_{Y|X})^2\big]+12\epsilon\\
&\leq  \frac{1}{1-\epsilon} \E_{X \sim \tilde{P}_X^0}\big[ \big(\mathbb{D}_{k_{\mathcal{Y}}}(P_{g(\theta,X)}, \tilde{P}^0_{Y|X})^2\big]+12\epsilon
\end{split}
\end{equation}
where the third inequality the fact uses the  that, since $|k_{\mathcal{Y}}|\leq 1$,  $\mathbb{P}(\mathbb{D}_{k_{\mathcal{Y}}}(P_{g(\theta,X)}, P^0_{Y|X})\leq 2)=1$, the penultimate inequality holds since $\epsilon< 1$ and the last inequality holds since $ \E_{X \sim Q_X}\big[ \big(\mathbb{D}_{k_{\mathcal{Y}}}(P_{g(\theta,X)}, \tilde{P}^0_{Y|X})^2\big]\geq 0$ for all $\theta\in\Theta$.

Then, applying  \eqref{eq:decompos} with $\theta=\tilde{\theta}_{Q,\epsilon}$ yields
\begin{equation}\label{eq:decompos2}
\begin{split}
\E_{X \sim P_X^0}\big[ &\mathbb{D}_{k_{\mathcal{Y}}}(P_{g(\tilde{\theta}_{Q,\epsilon},X)} , P^0_{Y|X})^2\big]\\
& \leq \frac{1}{1-\epsilon} \inf_{\theta\in\Theta}\E_{X \sim \tilde{P}_X^0}\big[ \mathbb{D}_{k_{\mathcal{Y}}}(P_{g(\theta,X)},  \tilde{P}^0_{Y|X})^2\big]+12\epsilon\\
&\leq \frac{1}{1-\epsilon} \inf_{\theta\in\Theta} \E_{X \sim \tilde{P}_X^0}\big[ \mathbb{D}_{k_{\mathcal{Y}}}(P_{g(\theta,X)}, P^0_{Y|X})^2\big]+\frac{12 \epsilon}{1-\epsilon}+12\epsilon\\
&\leq   \inf_{\theta\in\Theta} \E_{X \sim P_X^0}\big[ \mathbb{D}_{k_{\mathcal{Y}}}(P_{g(\theta,X)}, P^0_{Y|X})^2\big]+\frac{24 \epsilon}{1-\epsilon}+ \frac{12 \epsilon}{1-\epsilon}+16\epsilon\\
&\leq   \E_{X \sim P_X^0}\big[ \mathbb{D}_{k_{\mathcal{Y}}}(P_{g(\tilde{\theta}_{0},X)}, P^0_{Y|X})^2\big]+\frac{52 \epsilon}{1-\epsilon} 
\end{split}
\end{equation}
where the second inequality follows by swapping $\tilde{P}_{Y|X}^0$ and $P_{Y|X}^0$ in \eqref{eq:decompos} and the third one uses the fact that, since $|k_{\mathcal{Y}}|\leq 1$,
$$
 \E_{X \sim Q_X}\big[ \mathbb{D}_{k_{\mathcal{Y}}}(P_{g(\theta,X)}, P^0_{Y|X})^2\big]\leq 4,\quad\forall\theta\in\Theta.
$$

By assumption, $\tilde{\theta}_0$ is the unique minimizer of the function $\theta\mapsto \E_{X \sim P_X^0}\big[ \mathbb{D}_{k_{\mathcal{Y}}}(P_{g(\theta,X)}, P^0_{Y|X})^2\big]$ and therefore (see the proof of Lemma \ref{lemma:MMD:euclid})
\begin{align*}
\alpha=\inf_{\theta\in U^c}\Big(&\E_{X \sim P_X^0}\big[ \mathbb{D}_{k_{\mathcal{Y}}}(P_{g(\theta,X)},P^0_{Y|X})^2\big]-\E_{X \sim P_X^0}\big[ \mathbb{D}_{k_{\mathcal{Y}}}(P_{g(\tilde{\theta}_0,X)}, P^0_{Y|X})^2\big]\Big)>0.
\end{align*}
Together with \eqref{eq:decompos2}, this shows that if 
\begin{align*}
\frac{52 \epsilon}{1-\epsilon} <\alpha\Leftrightarrow \epsilon<\frac{\alpha}{52+\alpha}
\end{align*}
then
\begin{align*}
\E_{X \sim P_X^0}\big[ \mathbb{D}_{k_{\mathcal{Y}}}(P_{g(\tilde{\theta}_{Q,\epsilon},X)}, & P^0_{Y|X})^2\big] -\E_{X \sim P_X^0}\big[ \mathbb{D}_{k_{\mathcal{Y}}}(P_{g(\tilde{\theta}_0,X)}, P^0_{Y|X})^2\big]\\
&<\inf_{\theta\in U^c}\Big(\E_{X \sim P_X^0}\big[ \mathbb{D}_{k_{\mathcal{Y}}}(P_{g(\theta,X)}, \tilde{P}^0)^2\big]-\E_{X \sim P_X^0}\big[ \mathbb{D}_{k_{\mathcal{Y}}}(P_{g(\tilde{\theta}_0,X)}, \tilde{P}^0)^2\big]\Big)
\end{align*}
implying that  $\tilde{\theta}_{Q,\epsilon}\in U$. Consequently, using again  \eqref{eq:decompos2},
\begin{align*}
\frac{52 \epsilon}{1-\epsilon}\geq \E_{X \sim P_X^0}\big[ \mathbb{D}_{k_{\mathcal{Y}}}(P_{g(\tilde{\theta}_{Q,\epsilon},X)}, P^0_{Y|X})^2\big] -\E_{X \sim P_X^0}\big[ \mathbb{D}_{k_{\mathcal{Y}}}(P_{g(\tilde{\theta}_0,X)}, P^0_{Y|X})^2\big]\geq \mu\|\tilde{\theta}_{Q,\epsilon}-\theta_0\|
\end{align*}
and the result follows.
\end{proof}

\end{document}